\documentclass[10pt,a4paper]{amsart} 
\usepackage[utf8]{inputenc}
\usepackage[english]{babel}

\usepackage{amsmath}
\usepackage{amsfonts}
\usepackage{amssymb}

\usepackage{color}
\usepackage{ifthen}

\usepackage{graphics}
\usepackage{epsfig}

\usepackage{url}
\usepackage{hyperref}

\usepackage[a4paper,top=1.7cm,left=2cm,right=2cm,textheight=717.00946pt]{geometry}
\markleft{\hfill \textsc{Catalog of Properties of the First Isodynamic Point} \hfill}
\newcommand{\degrees}{^\circ}

\newsavebox{\figbox}

\newboolean{showIsoNumber}
\newboolean{showRef}
\long\def\void#1{}

\setboolean{showIsoNumber}{false}
\setboolean{showRef}{true}

\newcommand{\newcode}{\color{blue}$^*$\color{black} }
\newboolean{discovery}
\setboolean{discovery}{false}
\newcommand{\new}{\setboolean{discovery}{true}}

\definecolor{brown}{rgb}{0.6,0.3,0}

\newcommand{\propertyNumber}{
\textcolor{red}{\textbf{\Large Property \thesubsubsection.%
\ifthenelse{\boolean{discovery}}{\newcode}{}}}
}

\newcommand{\then}{$\blacktriangleright$\hspace{4pt}}

\newcommand{\pro}[5][]{
\bigskip
\refstepcounter{subsubsection}
\vbox{
\propertyNumber
\ifthenelse{\equal{#1}{}}{}{\textcolor{blue}{\textbf{(#1)}}}
\ifthenelse{\boolean{showIsoNumber}}{\hfill #3}{}
\ifthenelse{\boolean{showRef}}{\hfill #2}{}

\nopagebreak
\parbox{3.19in}{
\begin{center}
\includegraphics[width=0.3\textwidth]{#3.pdf}\par
\ifthenelse{\equal{#4}{}}{}
{\parbox{2in}{\begin{center}\textcolor{brown}{#4}\end{center}}\medskip}
\end{center}
\raggedright
\then #5
}
\smallskip
\label{property-#3}
\bigskip
\vfill
}
\setboolean{discovery}{false}
}

\newcommand{\prowide}[5][]{
\bigskip
\refstepcounter{subsubsection}
\vbox{
\propertyNumber
\ifthenelse{\equal{#1}{}}{}{\textcolor{blue}{\textbf{(#1)}}}
\ifthenelse{\boolean{showIsoNumber}}{\hfill #3}{}
\ifthenelse{\boolean{showRef}}{\hfill #2}{}

\nopagebreak
\parbox{3.19in}{
\begin{center}
\includegraphics[width=0.4\textwidth]{#3.pdf}\par
\ifthenelse{\equal{#4}{}}{}
{\parbox{2in}{\begin{center}\textcolor{brown}{#4}\end{center}}\medskip}
\end{center}
\raggedright
\then #5
}
\smallskip
\label{property-#3}
\bigskip
\vfill
}
\setboolean{discovery}{false}
}

\newlength{\innertextlength}
\newlength{\picturewidth}
\newlength{\exactpicturewidth}
\newsavebox{\myPicture}

\newcommand{\titledBox}[3]{
  \savebox{\myPicture}{\scalebox{#3}{\includegraphics{#1}}}%
  \settowidth{\exactpicturewidth}{\usebox{\myPicture}}%
  \setlength{\fboxrule}{1pt}%
  \savebox{\myPicture}{\framebox{\usebox{\myPicture}}}%
  \settowidth{\picturewidth}{\usebox{\myPicture}}%
  \begin{minipage}{\picturewidth}
  \lineskip=-1pt
  \setlength{\innertextlength}{\exactpicturewidth}
  \addtolength{\innertextlength}{-0.1in}
  \framebox[\picturewidth]{\parbox{\exactpicturewidth}{\centerline{\parbox{\innertextlength}{\centering #2}}}}
  \usebox{\myPicture}
  \end{minipage}
}

\newcommand{\mysubsection}[1]
{
\medskip
\fbox{\refstepcounter{subsection}%
\textbf{\thesubsection\ #1}}%
\addcontentsline{toc}{subsection}{\hspace{36pt}\thesubsection. #1}%
}

\newcommand{\Hrule}{\rule{\linewidth}{1.5pt}\\}

\newcommand{\mysection}[1]{%
\refstepcounter{section}%
\Hrule
\vspace{-3pt}
{\Large\textbf{\thesection.\hfil #1\\}}
\vspace{-6pt}
\Hrule
\addcontentsline{toc}{section}{\hspace{18pt}\thesection. #1}
}

\newcommand{\extra}[1]{\vspace{-15pt}{\Large $\bullet$}\ \ {\small #1}}

\begin{document}
% ===================
International Journal of  Computer Discovered Mathematics (IJCDM) \\
ISSN 2367-7775 \copyright IJCDM \\
Volume 6, 2021, pp. 108--136  \\
Received 27 Nov. 2021. Published on-line 31 Dec. 2021. Revised 1 March 2022. \\ 
web: \url{http://www.journal-1.eu/} \\
\copyright The Author(s) This article is published 
with open access.\footnote{This article is distributed under the terms of the Creative Commons Attribution License which permits any use, distribution, and reproduction in any medium, provided the original author(s) and the source are credited.}\\
% ===========================   

\bigskip
\medskip

\begin{center}
	{\Large \textbf{Catalog of Properties of the First}} \\
	\medskip
	{\Large \textbf{Isodynamic Point of a Triangle}} \\
	\bigskip

	\textsc{Stanley Rabinowitz} \\

	545 Elm St Unit 1,  Milford, New Hampshire 03055, USA \\
	e-mail: \href{mailto:stan.rabinowitz@comcast.net}{stan.rabinowitz@comcast.net}\\
	web: \url{http://www.StanleyRabinowitz.com/} \\
	
\end{center}
\bigskip

% ==============================
\textbf{Abstract.} The first isodynamic point of a triangle is one of many notable points
associated with a triangle. It is named X(15) in the Encyclopedia of Triangle Centers.
This paper surveys known results about this point and gives additional properties that were
discovered by computer.

\medskip
\textbf{Keywords.} triangle geometry, first isodynamic point, computer-discovered mathematics, GeometricExplorer.

\medskip
\textbf{Mathematics Subject Classification (2020).} 51M04, 51-08.

\DeclareGraphicsExtensions{.pdf}

\newcommand{\fig}[1]
{
\begin{center}
\includegraphics[width=0.4\linewidth]{#1}
\label{fig:#1}
\end{center}
}

\newcommand{\blueurl}[1]{\\\textcolor{blue}{\url{#1}}}

% ================================
% 1 Introduction 
% ================================

%\addtocontents{toc}{\hfil and classification scheme}

\section*{Introduction}
\label{section:intro}

\parbox{3.5in}{
The first isodynamic point of a triangle is one of many notable points
associated with a triangle. It is the point with trilinear coordinates

\bigskip
\hspace{0.5in}
$\displaystyle \left(\sin(A+\frac{\pi}{3}): \sin(B+\frac{\pi}{3}): \sin(C+\frac{\pi}{3})\right)$
\bigskip

with respect to the triangle. It is named $X_{15}$ in the Encyclopedia of Triangle Centers
\cite{ETC15}.

\bigskip

\textbf{Geometrical definition.}
The interior and exterior angle bisectors of angle $A$ of $\triangle ABC$
intersects side $BC$ of the triangle (or its extension) in two points, $A_1$ and $A_2$.
The circle with diameter $A_1A_2$ is called the $A$-Apollonius circle and is named $C_A$.
Circles $C_B$ and $C_C$ are defined similarly.
The points in which the three Apollonius circles intersect are the isodynamic points of the triangle. The one inside $\odot ABC$ is the 1st isodynamic point and is named $S$.
The other point of intersection is the 2nd isodynamic point and is named $S'$.
}
\hfill
\raisebox{-1.3in}{\scalebox{0.57}{\includegraphics{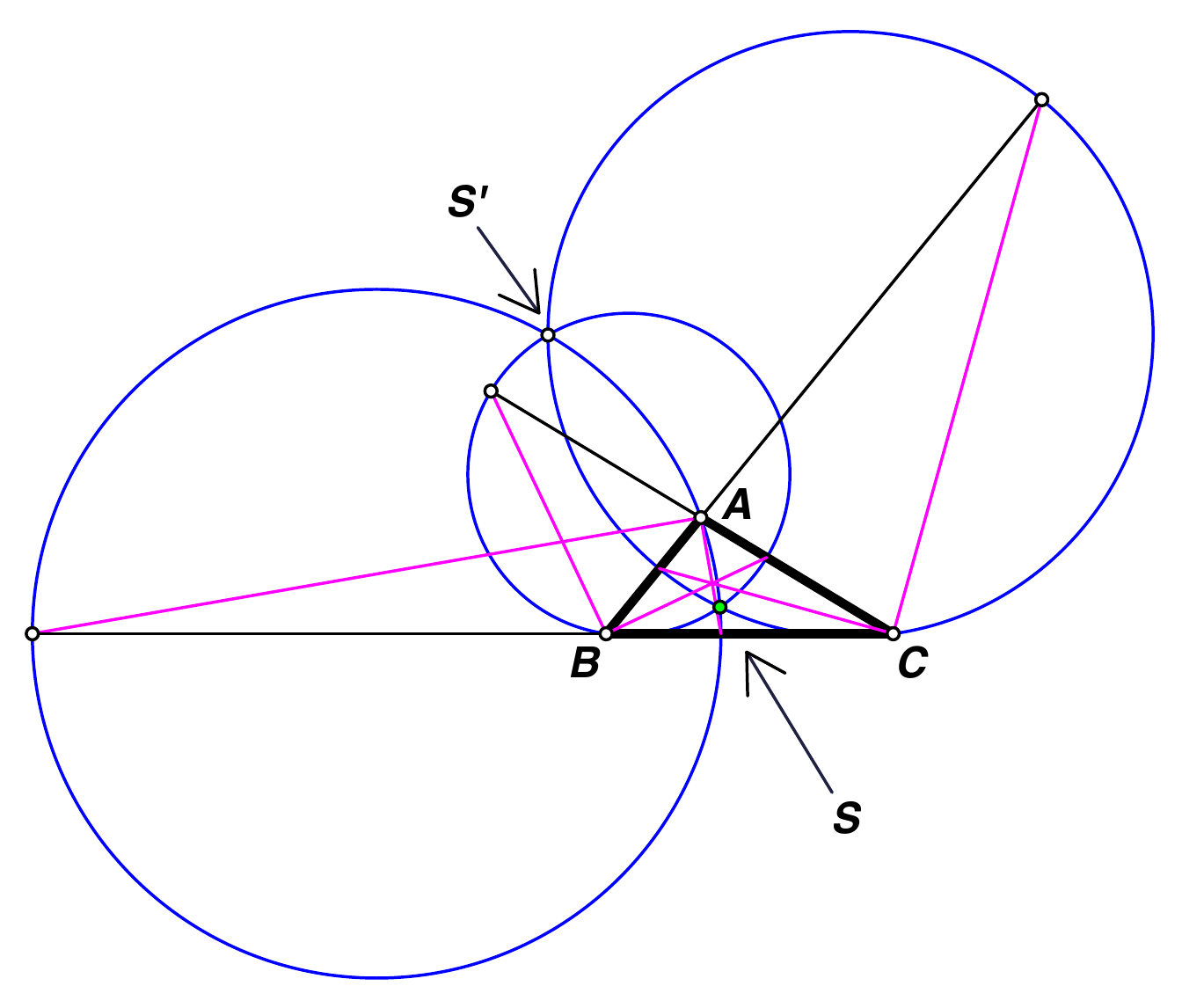}}}

\bigskip
\textbf{Scope of this catalog.}
The mathematical literature is vast. We do not attempt to catalog every property involving
the 1st isodynamic point that appears somewhere in print or on the internet.
We do try to catalog any property that is simple or elegant or that can be obtained
from the configuration associated with one of our top-level classifications
(Triangle plus $S$, Triangle with $S$ and other points, Triangle plus lines through $S$, Quadrilateral plus $S$, etc.)
by applying at most one common geometrical construction (drop a perpendicular, draw an angle bisector,
construct a centroid, etc.) When analyzing triangle centers, we only look at the common ones, $X_1$ through $X_{20}$.

\void{
topmargin=\the\topmargin, 
head=\the\headheight, 
sep=\the\headsep, 
odd=\the\oddsidemargin, 
even=\the\evensidemargin, 
height=\the\textheight, 
width=\the\textwidth, 
foot=\the\footskip, 
paper=\the\paperheight, \the\paperwidth, 
margin=\the\marginparwidth, \the\marginparsep
}

\twocolumn[\huge\hfil \textbf{Classification Scheme}\\]
%\begin{center}
%\parbox{4in}{
%\small
%\setlength{parskip}{-2pt}
\tableofcontents
%}
%\end{center}
\onecolumn

\textbf{Figures.}
The figures in this paper are \emph{not} decorative. In order to reduce verbiage and declutter the catalog, given information
appearing in the figure is not repeated in words. For example, if a figure shows line segments
$WX$ and $YZ$ meeting at a point $P$, we do not state in words that $P$ is the intersection of $WX$ and $YZ$.
This should make it easier for people who don't read English to use this catalog.

A solid line or circle through multiple points means that it is given that these points lie on the same line
or circle.
A dashed line or circle through multiple points means that the conclusion of the theorem or result is that
these points lie on the same line or circle.

\begin{center}
\titledBox{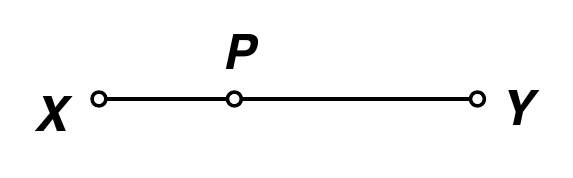}{$P$ is a point on segment $XY$.}{0.6}
\quad
\titledBox{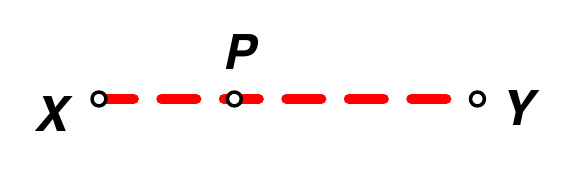}{We can conclude that $X$, $P$, and $Y$ colline.}{0.62}
\qquad
\titledBox{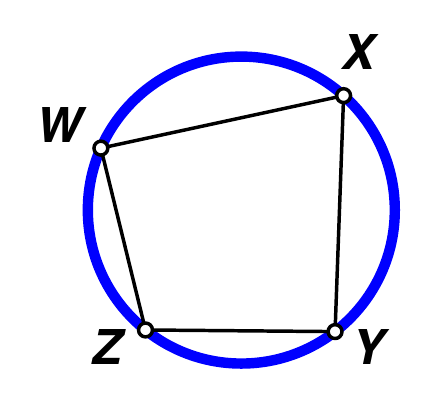}{It is given that $WXYZ$ is a cyclic quadrilateral.}{0.68}
\quad
\titledBox{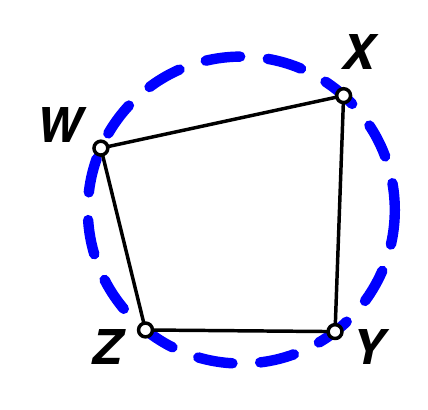}{We can conclude that $W$, $X$, $Y$, and $Z$ are concyclic.}{0.7}
\end{center}

\goodbreak
Angles that are given to be equal are marked with the same filled circle.
Angles that are concluded to be equal are shaded with the same color.

\begin{center}
\titledBox{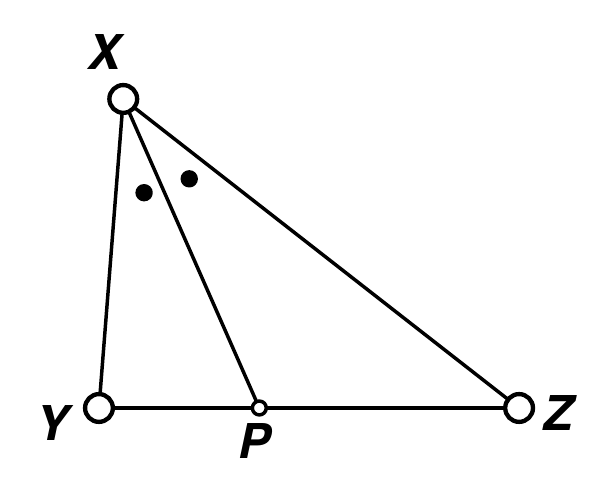}{It is given that $\angle YXP=\angle PXZ$.}{0.6}
\quad
\titledBox{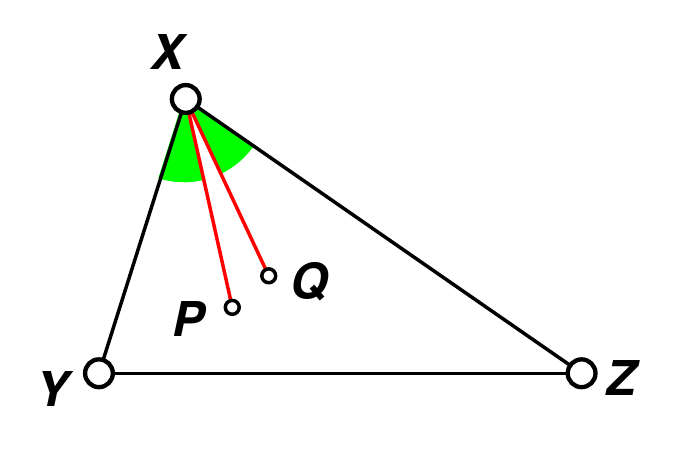}{We can conclude that $\angle YXP=\angle QXZ$.}{0.66}
\end{center}

Two perpendicular brown lines at the point of intersection of two circles
means that we can conclude that the circles are orthogonal (have perpendicular
tangents at that point).

\begin{center}
\titledBox{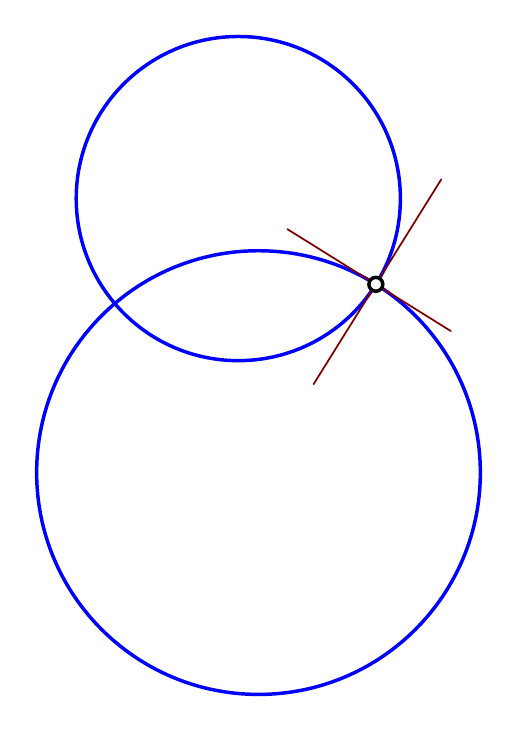}{We can conclude that the circles are orthogonal.}{0.8}
\end{center}

If the title of a section or subsection describes a feature of a figure,
then we do not repeat this description if it is obvious from the figure.
For example, in a subsection entitled ``equilateral triangles", a triangle
highlighted in yellow that looks equilateral can be assumed to be an equilateral triangle.

A right-angle marker is used to indicate two lines that are given to be perpendicular.
All angles are directed angles.
The 1st isodynamic point of $\triangle ABC$ is always colored green.
Given information that is not obvious from the associated figure is shown in brown
text directly beneath the figure.

\newpage
\textbf{Parts of a triangle.}
In order to help with the classification process, we give names for various line segments
associated with a triangle.
A line segment from a vertex of a triangle to a non-vertex point on the opposite side is called a \emph{cevian}.
A line segment joining points on two sides of the triangle is called a \emph{chord}.
A chord parallel to a side of the triangle is called a \emph{parachord}.
If the endpoints of a chord joining points on two sides of a triangle forms a cyclic quadrilateral
with the endpoints of the third side, the chord is called an \emph{antiparallel}.

\begin{center}
\titledBox{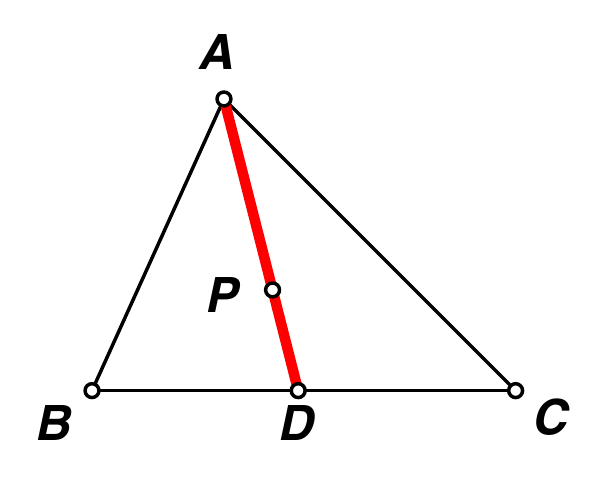}{cevian}{0.5}
\quad
\titledBox{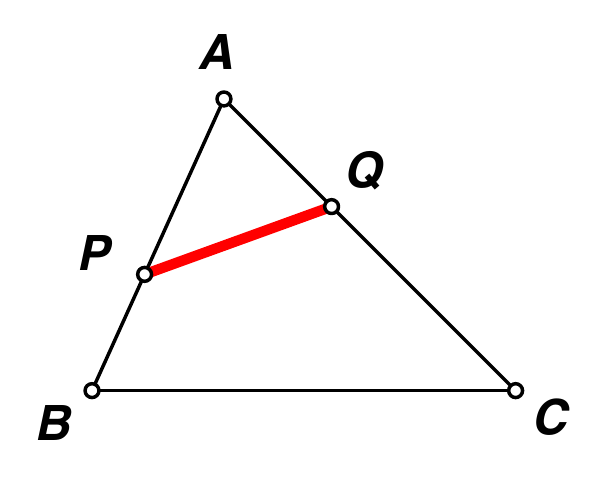}{chord}{0.5}
\quad
\titledBox{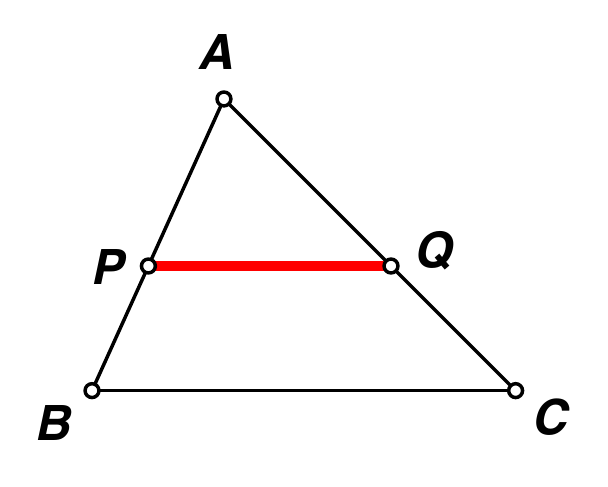}{parachord}{0.5}
\quad
\titledBox{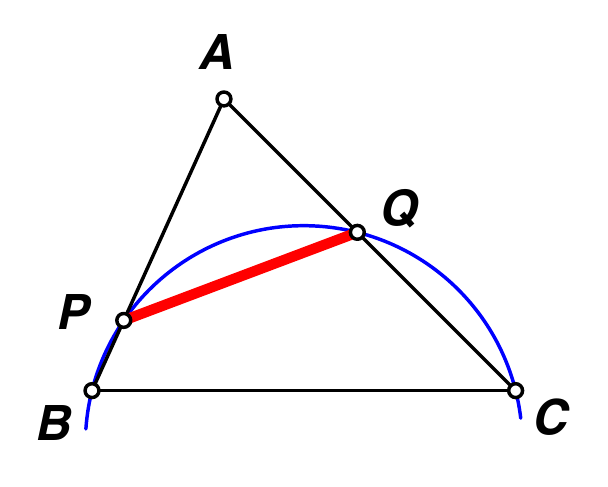}{antiparallel}{0.5}
\end{center}

If $P$ is a point inside a triangle, the line segment from $P$ to a vertex is called a \emph{spoke}.
The line segment from $P$ to the foot of the perpendicular from $P$ to a side of the triangle
is called an \emph{apothem}. A line segment from $P$ parallel to a side of the triangle ending
on another side of the triangle is called a \emph{pararadius}.
A line segment from $P$ to a side of the triangle that forms an angle of $n\degrees$ with that side
is called an \emph{$n\degrees$-incline}.

\begin{center}
\titledBox{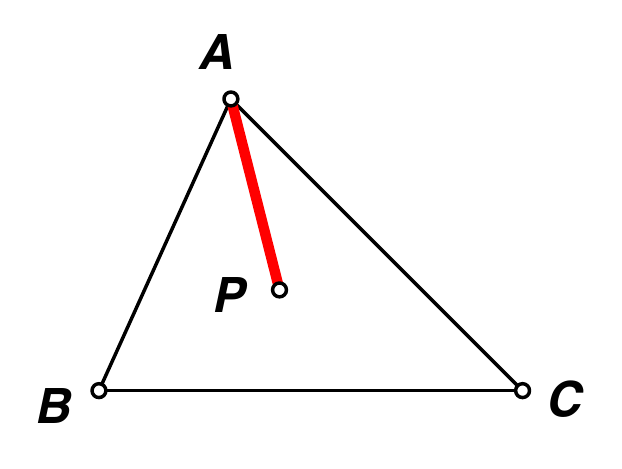}{spoke}{0.5}
\quad
\titledBox{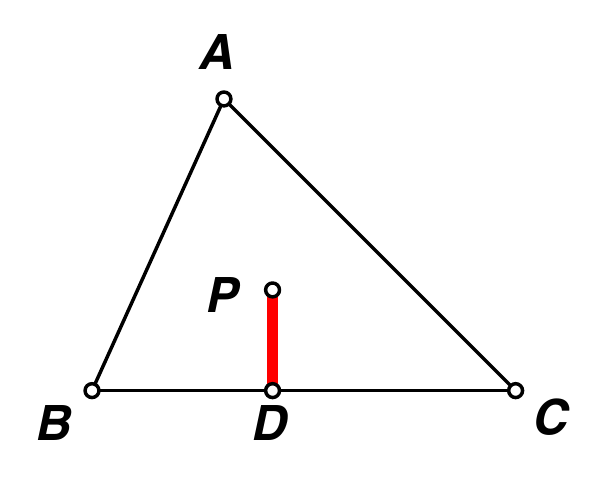}{apothem}{0.5}
\quad
\titledBox{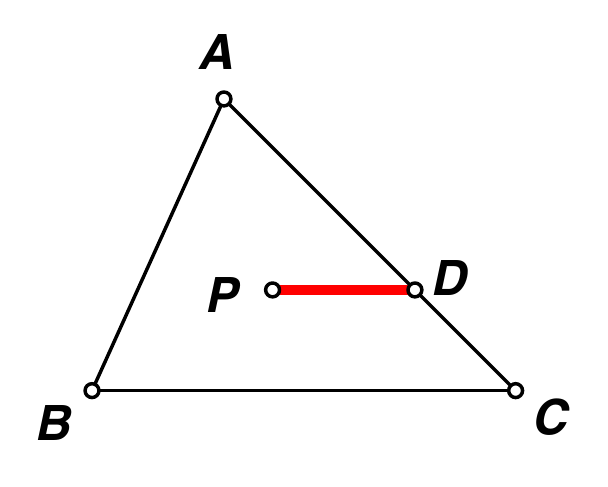}{pararadius}{0.5}
\quad
\titledBox{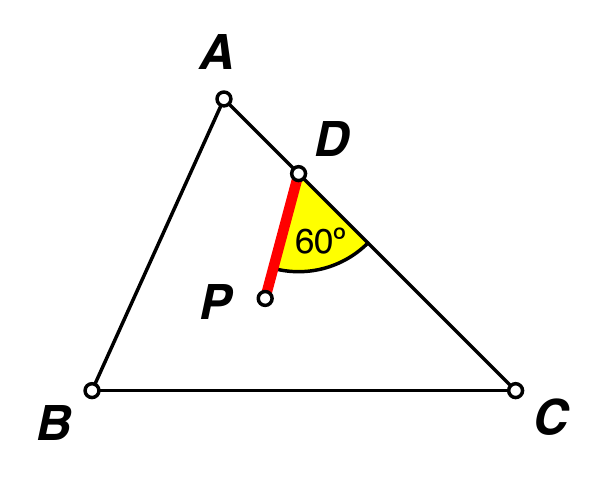}{$60\degrees$-incline}{0.5}
\end{center}

If $P$ is a point inside a triangle, the line segment from a vertex passing through $P$
and extending to the circumcircle of $\triangle ABC$ is called a \emph{circumcevian}.
The line segment from a vertex passing through $P$
and extending to the circumcircle of $\triangle BPC$ is called a \emph{circlecevian}.
The endpoint of a cevian through $P$ (other than a vertex) is called a \emph{trace}
and the line segment along a side from that trace to a vertex of the triangle is called a \emph{trace segment}.
The line segment from the midpoint of the side of a triangle
extending outward to the furthest point on the circumcircle is called a \emph{sagitta}.

\begin{center}
\titledBox{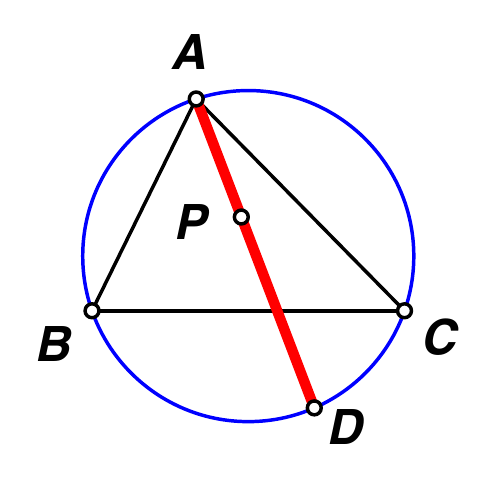}{circumcevian}{0.5}
\quad
\titledBox{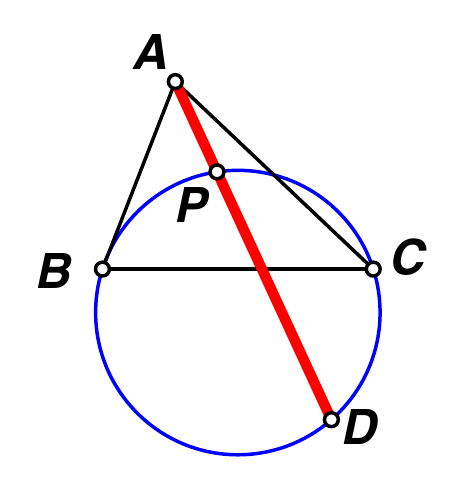}{circlecevian}{0.5}
\quad
\titledBox{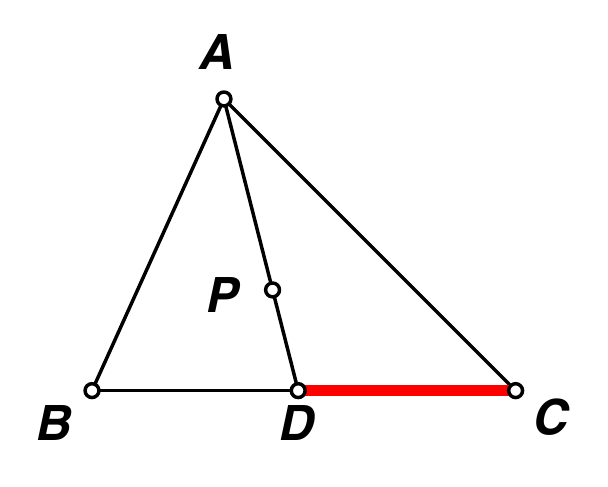}{trace segment}{0.5}
\quad
\titledBox{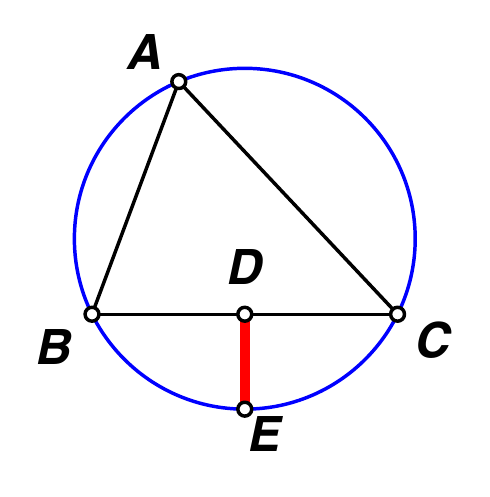}{sagitta}{0.5}
\end{center}

%In the categories used in this paper, cevians, spokes, apothems, etc. shall always be ones through the point $S$.

\bigskip
\textbf{Discoveries.}
An asterisk after a property number indicates that the property was discovered by computer,
either by using GeometricExplorer, Mathematica, or Geometer's Sketchpad.
If a reference is given, this means that the result was posted to an online forum in the
hope that some forum member might find a geometrical proof of the property.

\bigskip
\textbf{Acknowledgement.}

I would like to thank the following people who did thorough reviews of this paper and who provided
additional references and results about the isodynamic points of a triangle:
Kadir Alt{\i}nta\c{s}, Dao Thanh Oai, Tran Quang Hung, Kousik Sett.

\newpage

\textbf{Notation.}

\begin{center}
\begin{tabular}{|l|l|}
\hline
\multicolumn{2}{|c|}{\large \strut Notation used when describing properties}\\
\hline
\textbf{Notation}&\textbf{Description}\\
\hline
$\triangle XYZ$&Triangle $XYZ$\\
$a$, $b$, $c$&The lengths of the sides of $\triangle ABC$\\
$s$&$(a+b+c)/2$\\
$\odot XYZ$&The circle through points $X$, $Y$, and $Z$\\
$H$&The orthocenter of $\triangle ABC$\\
$K$&The symmedian point of $\triangle ABC$\\
$K$&When used in an expression, $K$ denotes the area of $\triangle ABC$.\\
%$I_a$&The $A$-excenter of $\triangle ABC$\\
$I$&The incenter of $\triangle ABC$\\
$M$&The centroid of $\triangle ABC$\\
$O$&The circumcenter of $\triangle ABC$\\
$r$&The inradius of $\triangle ABC$\\
$R$&The circumradius of $\triangle ABC$\\
$R(XYZ)$&The circumradius of $\triangle XYZ$\\
$S$&The 1st isodynamic point of $\triangle ABC$\\
$S(XYZ)$&The 1st isodynamic point of $\triangle XYZ$\\
$S'$ or $T$&The 2nd isodynamic point of $\triangle ABC$\\$\Omega$&The 1st Brocard point of $\triangle ABC$\\
$\Omega'$&The 2nd Brocard point of $\triangle ABC$\\
%$WX\perp YZ$&Lines $WX$ and $YZ$ are perpendicular.\\
%$WX\parallel YZ$&Lines $WX$ and $YZ$ are parallel.\\
$[F]$&The area of figure $F$\\
$\phi$&The golden ratio, $(1+\sqrt 5)/2$\\
$X$ -- $Y$ -- $Z$&The points $X$, $Y$, and $Z$ colline (i.e. are collinear).\\
$\angle XYZ$&Directed angle XYZ. This is the angle through which ray $\overrightarrow{YX}$\\
&must be rotated counterclockwise in order to coincide with ray $\overrightarrow{YZ}$.\\
$X_n$&The nth Kimberling center of $\triangle ABC$ (see \cite{ETC})\\
$X_n(XYZ)$&The nth Kimberling center of $\triangle XYZ$\\
%$\blacktriangleright$ statement&The statement follows from the given information.\\
\hline
\end{tabular}
\end{center}

\bigskip

\void{
\fbox{
\parbox{3.5in}
{
\setcounter{section}{45}
\setcounter{subsection}{6}
\setcounter{subsubsection}{2}
\new
\pro[Stan's Sample Theorem]{[123]}{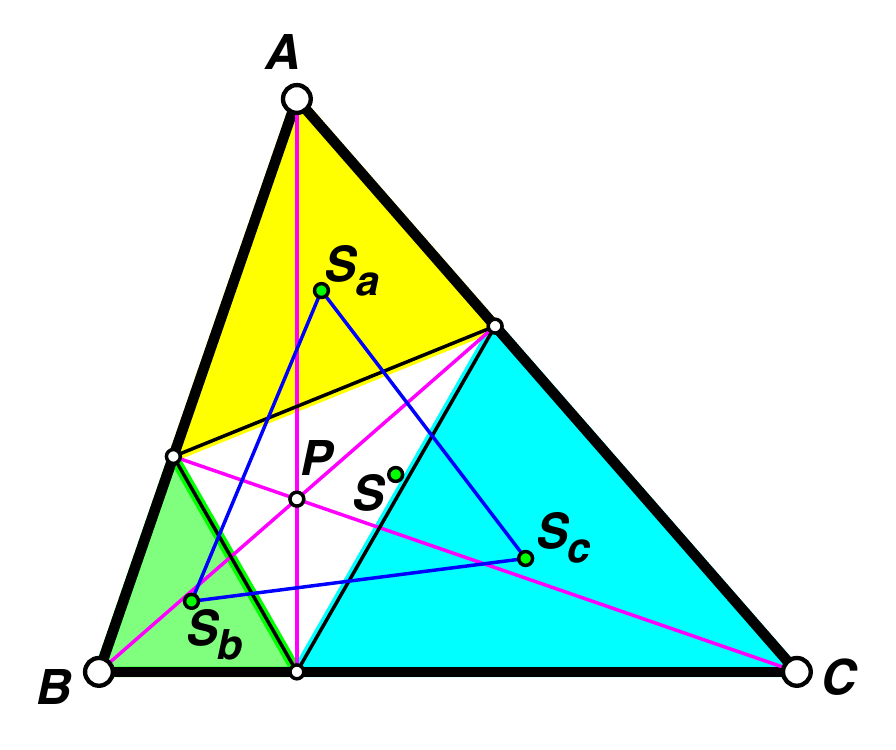}
{$P$ is the awesome point of $\triangle ABC$}
{$PS$ is ultraparallel to the foobar line of $\triangle ABC$.}
\setboolean{discovery}{false}
\extra{If we replace $P$ by the abysmal point of $\triangle ABC$,
then $PS$ is antiperpendicular to the foobar line.}
}
}
}

\textbf{Key to the property listings.}

\begin{center}
\scalebox{0.8}{\includegraphics{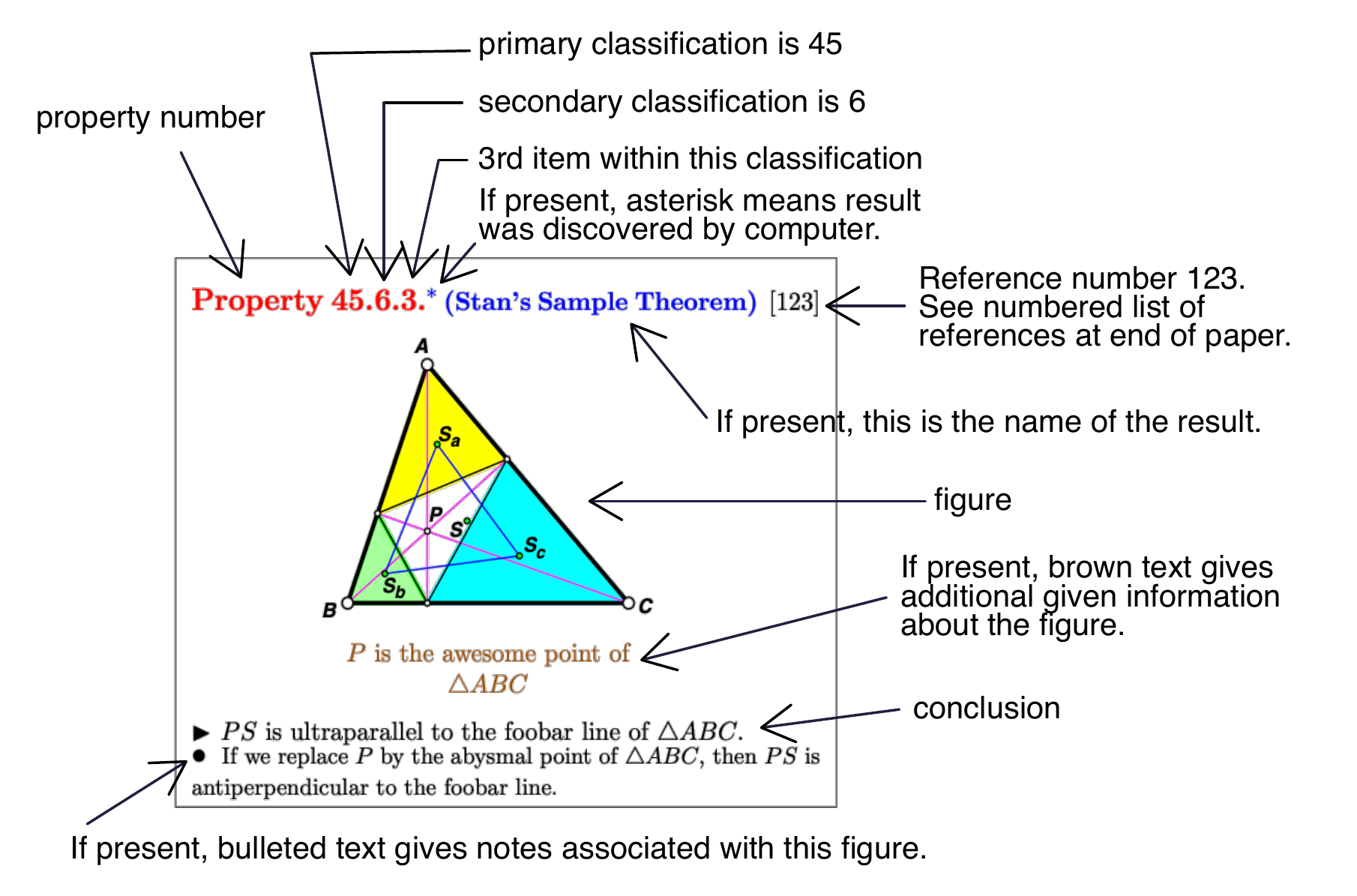}}
\end{center}

Since there are three shaded triangles, this means that the 1st isodynamic points
of these triangles will be named $S_a$, $S_b$, and $S_c$ [or $S_1$, $S_2$, and $S_3$].
The 2nd isodynamic points will be named $T_a$, $T_b$, and $T_c$.

\twocolumn[\huge\hfil \textbf{1st Isodynamic Point $S$}\\]
\raggedright
\raggedbottom
\vspace*{-22pt}
\mysection{Triangle plus $S$}

\void{
This section gives properties of a configuration consisting of a triangle with its
isodynamic point $S$.
}

\mysubsection{location of $S$}

\new
\pro{}{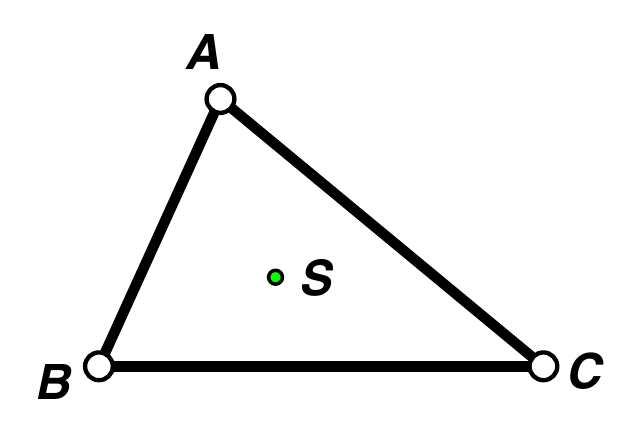}
{Each angle of $\triangle ABC$ has measure less than $120\degrees$.}
{$S$ lies in the interior of $\triangle ABC$.}

\new
\pro{}{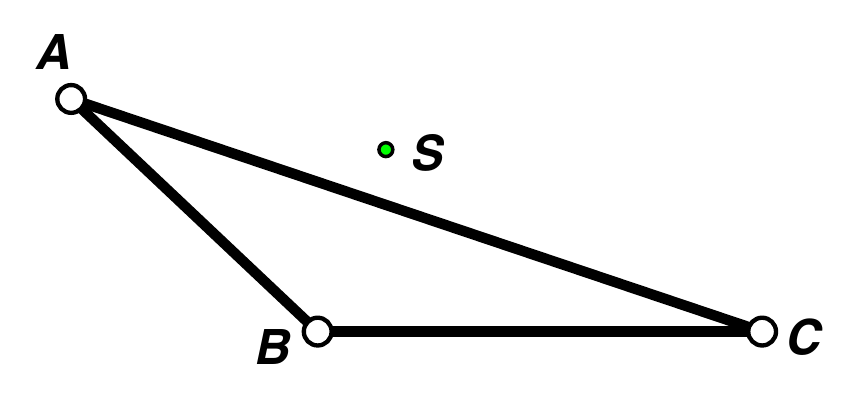}
{$\angle B>120\degrees$.}
{$S$ lies outside $\triangle ABC$.}

\pro{\cite{Wiki}}{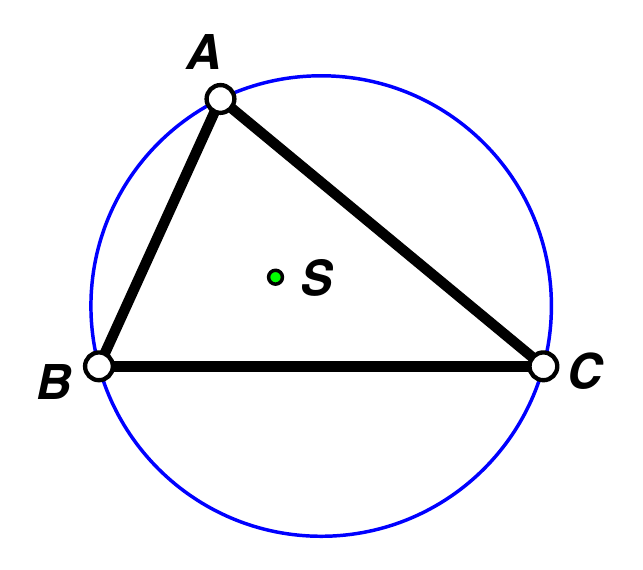}
{}
{$S$ lies inside the circumcircle of $\triangle ABC$.}
\extra{Equivalently, $SO<R$.}

\mysubsection{metric properties}

\nopagebreak
\pro[Tripolar Coordinates]{\cite[p.~106]{Gallatly}}{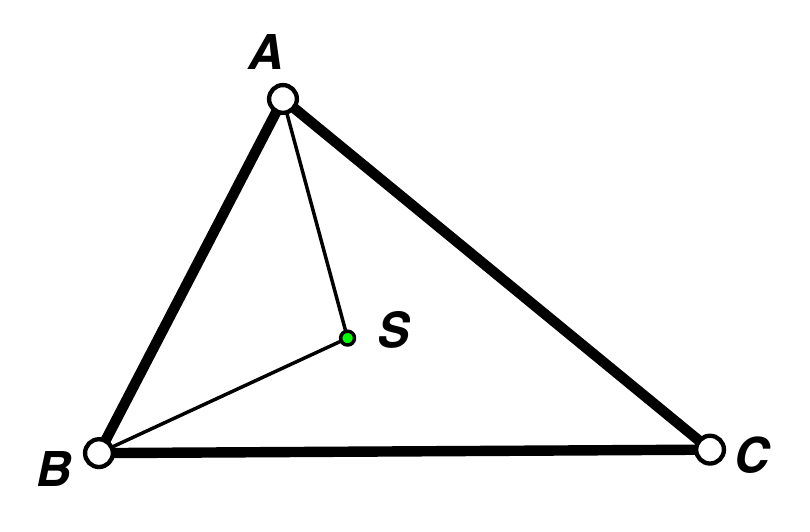}
{}
{$\displaystyle \frac{AS}{BS}=\frac{AC}{BC}$}
\extra{A similar property holds for $S'$. By symmetry, we have $AS\cdot BC=BS\cdot AC=CS\cdot AB$.}

\pro[Spoke Length]{\cite{Fettis}}{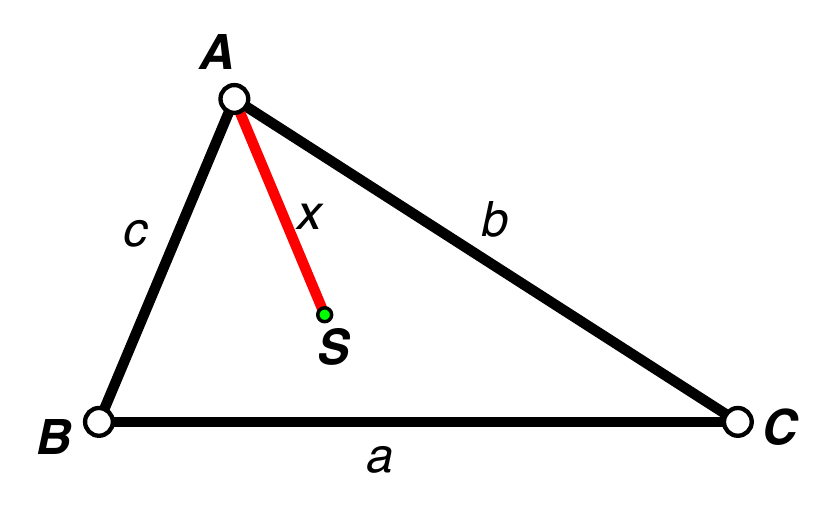}
{\void{Let $x=AS$}}
{$\displaystyle x=\frac{bc\sqrt2}{\sqrt{a^2+b^2+c^2+4K\sqrt3}}$}

\mysubsection{angle properties}

\pro{\cite{RG4850}}{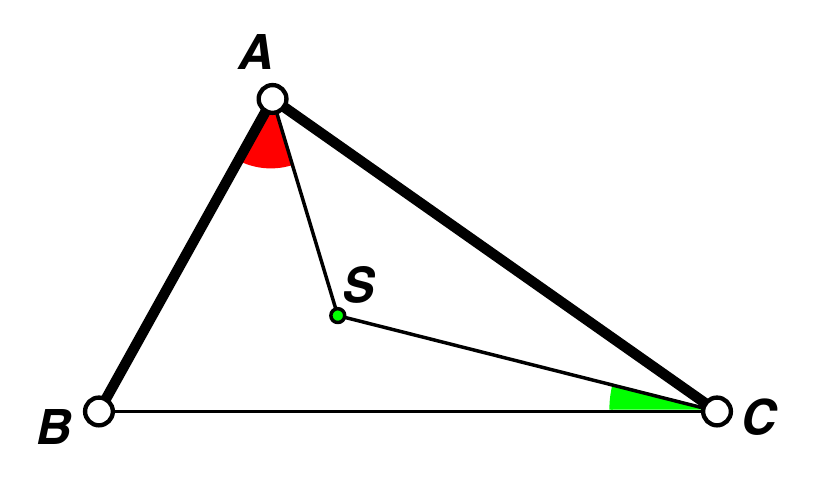}
{}
{$\angle BAS+\angle SCB=60\degrees$}

\pro{\cite[p.~296]{Johnson}}{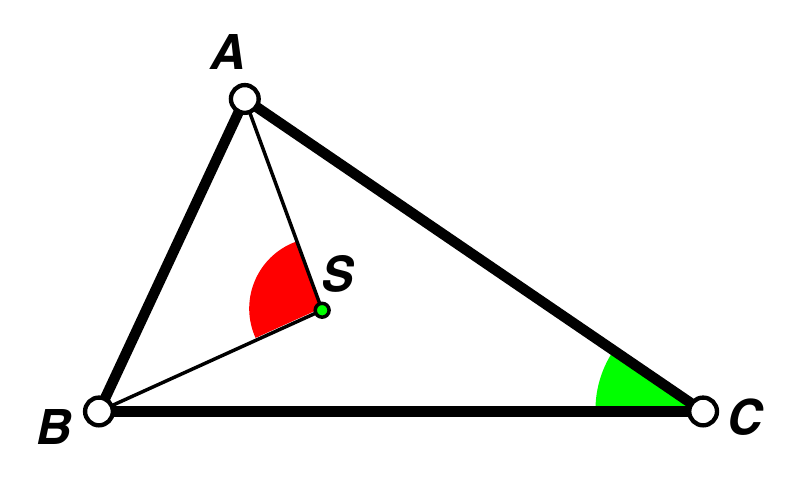}
{}
{$\angle ASB=\angle ACB+60\degrees$}
\extra{A similar result holds for $S'$: $\angle ASB=\angle ACB-60\degrees$}

\mysubsection{geometrical properties}

\new
\pro{\cite{RG7294}}{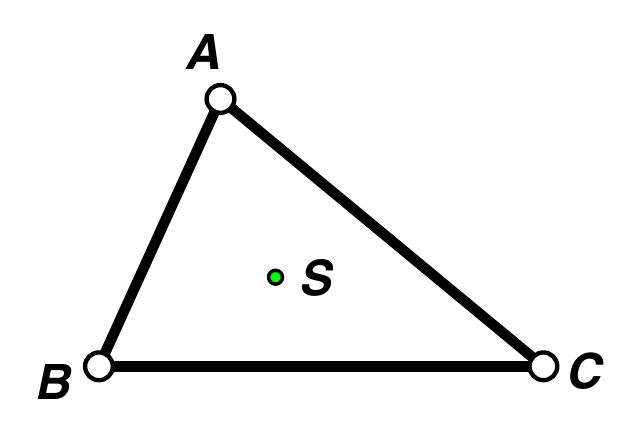}
{}
{$A$ is the 2nd isodynamic point of $\triangle BCS$.}
\setboolean{discovery}{false}

\mysection{Triangle plus $S$ and constructions}

\mysubsection{angle bisectors}

\new
\pro{\cite{RG4846}}{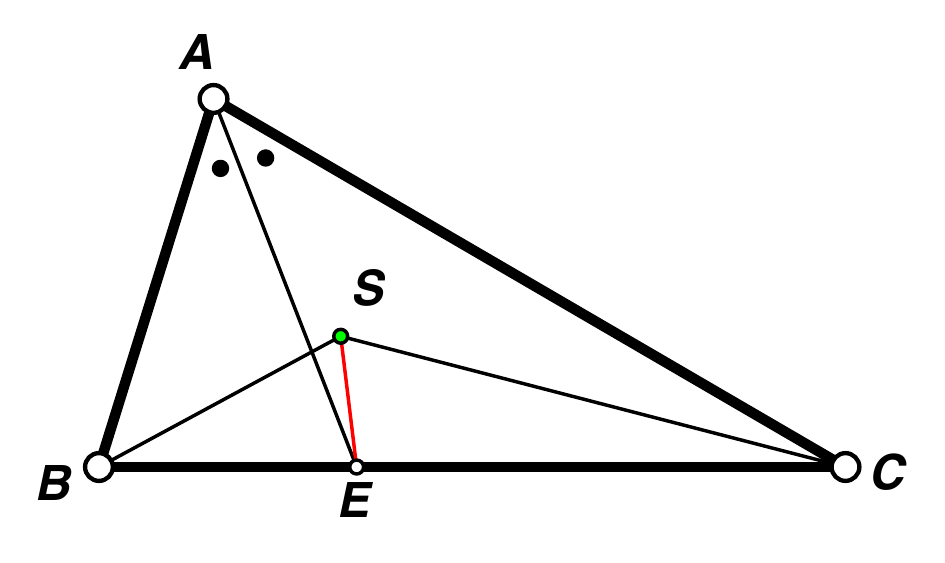}
{\void{Construct $AE$, the angle bisector of $\angle BAC$.}}
{$\angle BSE=\angle ESC$}

\pro{\cite{RG4846}}{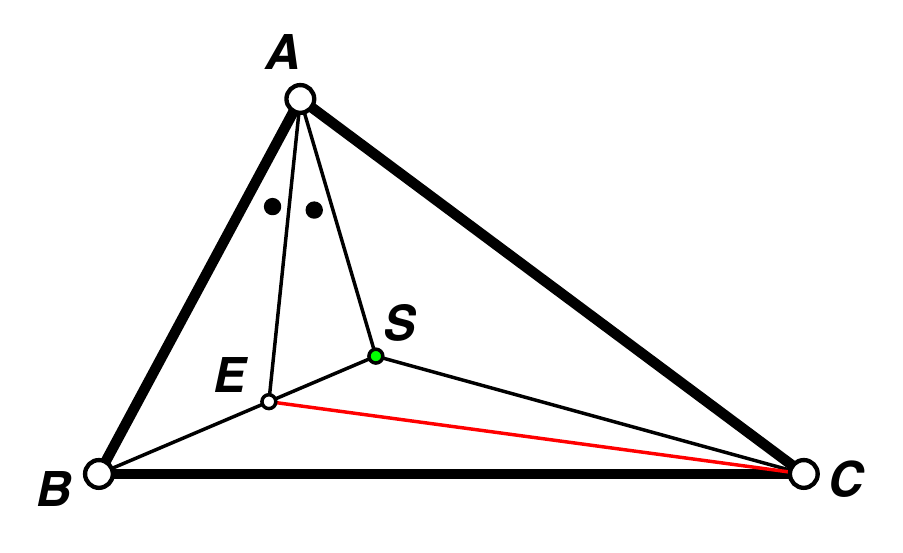}
{\void{Construct $AE$, the angle bisector of $\angle BAS$.}}
{$\angle SCE=\angle ECB$}

\new
\pro{\cite{RG4904}}{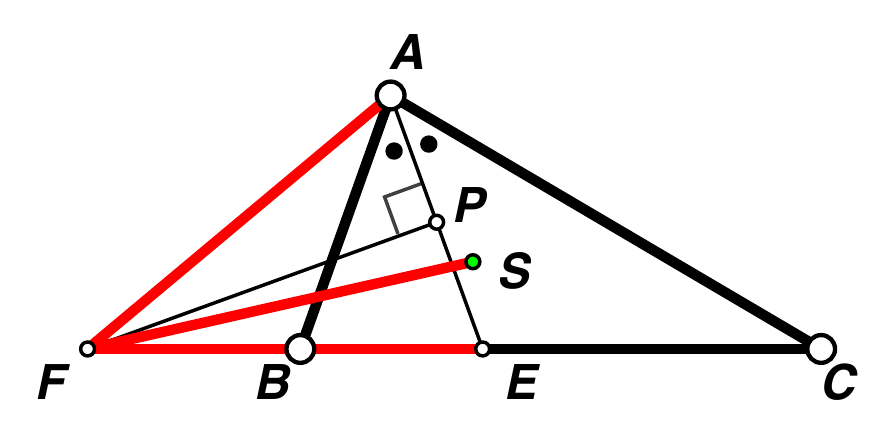}
{$AP=PE$}
{$E$ -- $S$ -- $A'$}

\pro{\cite{Wiki}}{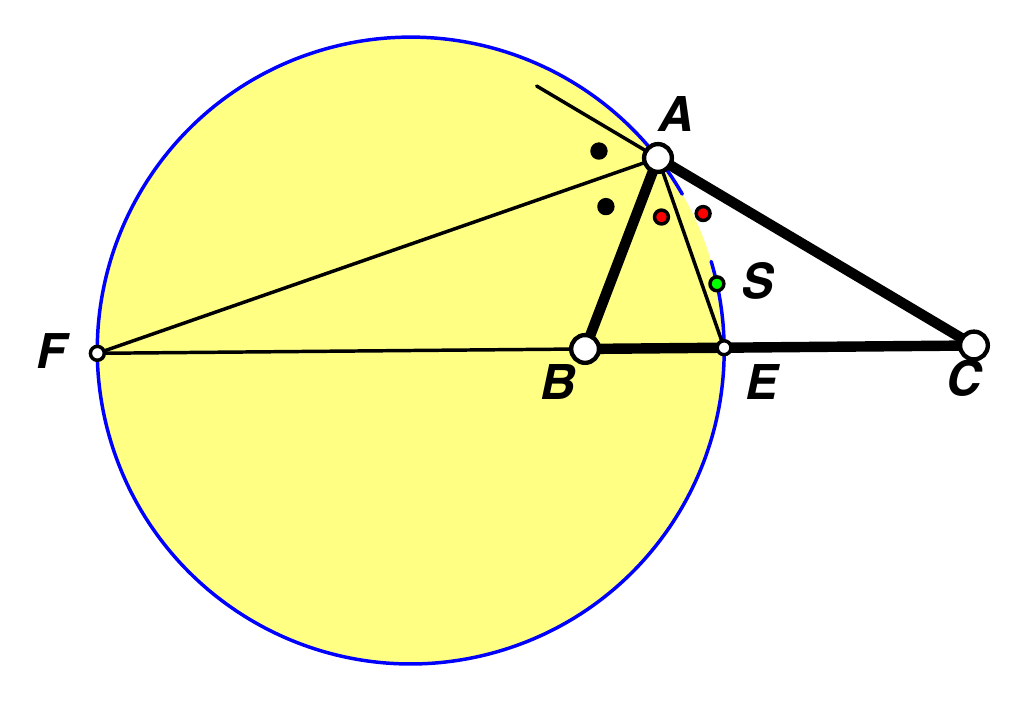}
{\void{Construct interior and exterior angle bisectors of $\angle BAC$.}}
{Circle with diameter $EF$ passes through $S$.}

\void{
\mysubsection{six angle bisectors}
% does not involve S
See\cite{RG9037}
See also \cite{Euclid3072}.
}

\mysubsection{circles}

\pro{\cite{Wiki}}{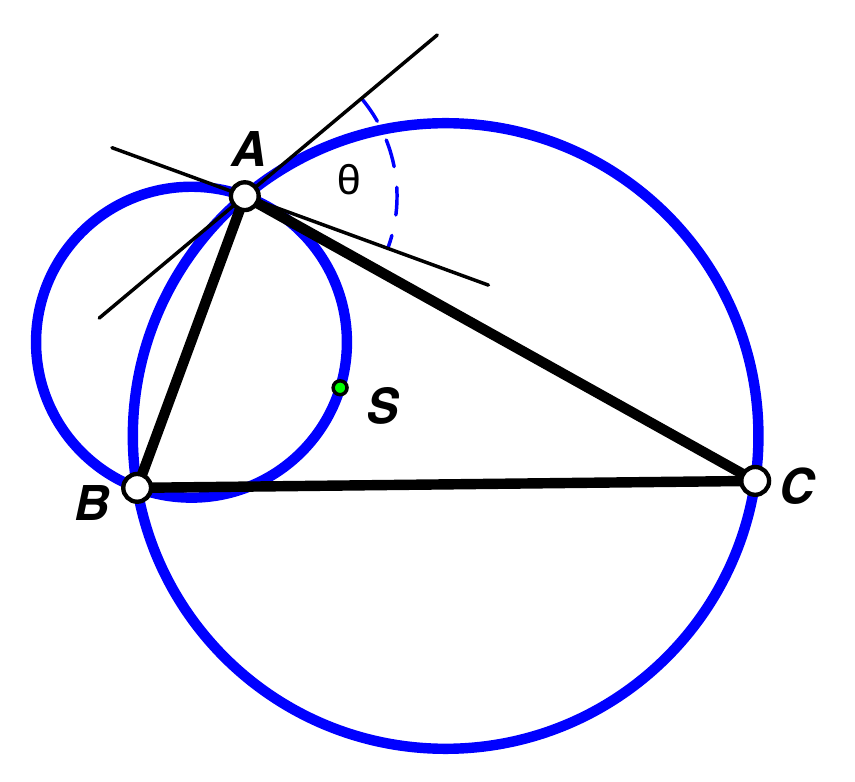}
{\void{Construct $\odot ABC$ and $\odot ABS$.}}
{$\theta=60\degrees$}
\extra{In other words, the circles meet at an angle of $60\degrees$.}

\pro{\cite{PGR14-Nov-2021}}{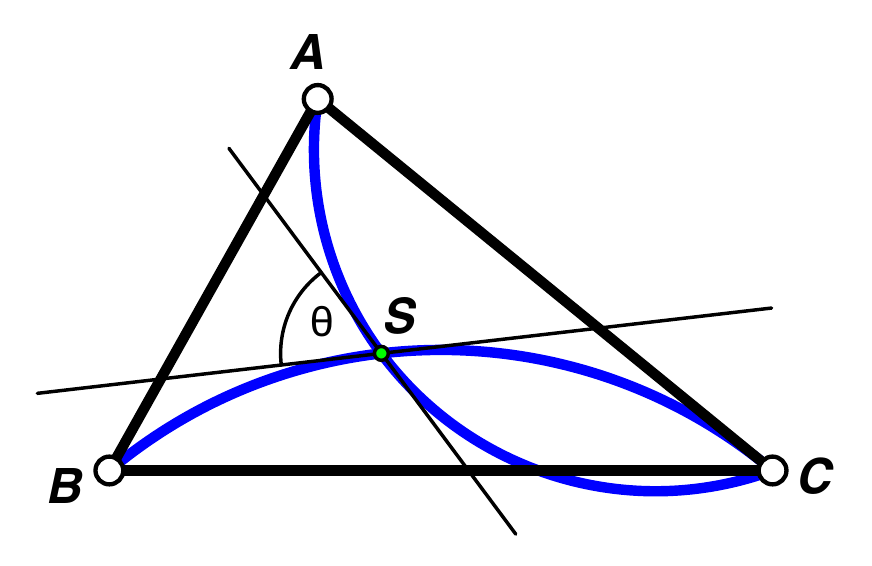}
{}
{$\theta=60\degrees$.}
\extra{This is a special case of Property \ref{property-iso-98}.}

\pro{\cite{PGR14-Nov-2021}}{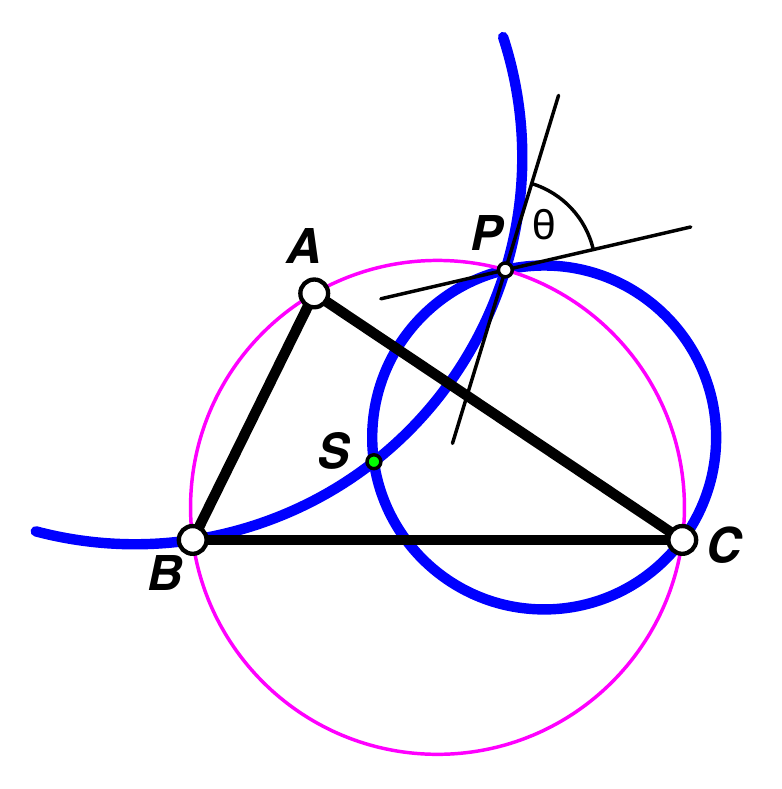}
{}
{$\theta=60\degrees$.}

\pro{\cite{PGR14-Nov-2021}}{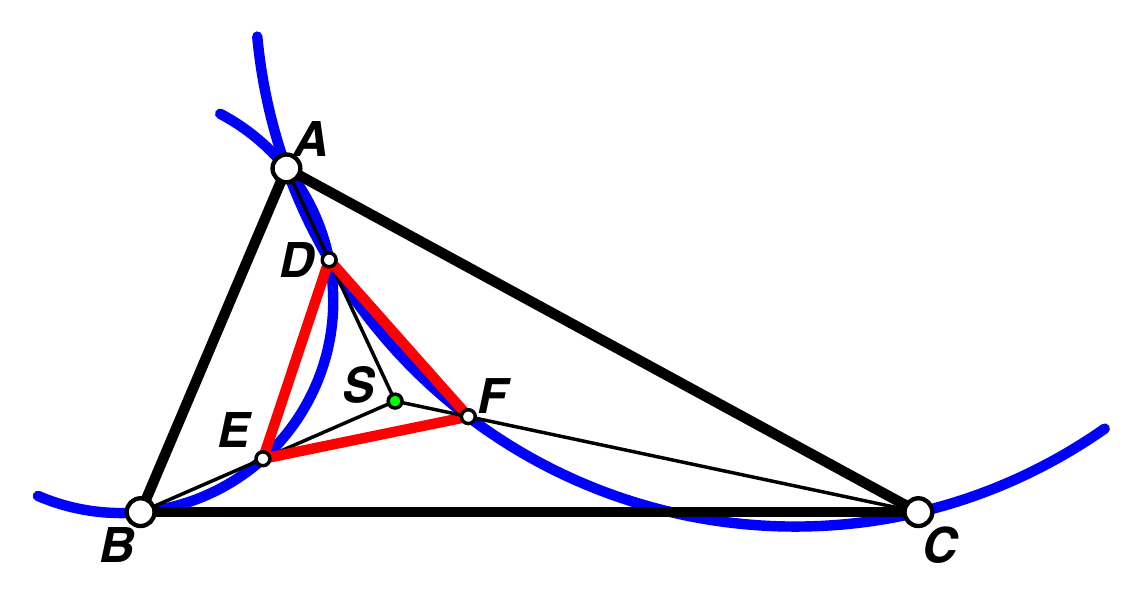}
{$D$ is any point on $AS$.}
{$\triangle DEF$ is equilateral.\\
\then $B$, $E$, $F$, $C$ are concyclic.}

\new
\pro{}{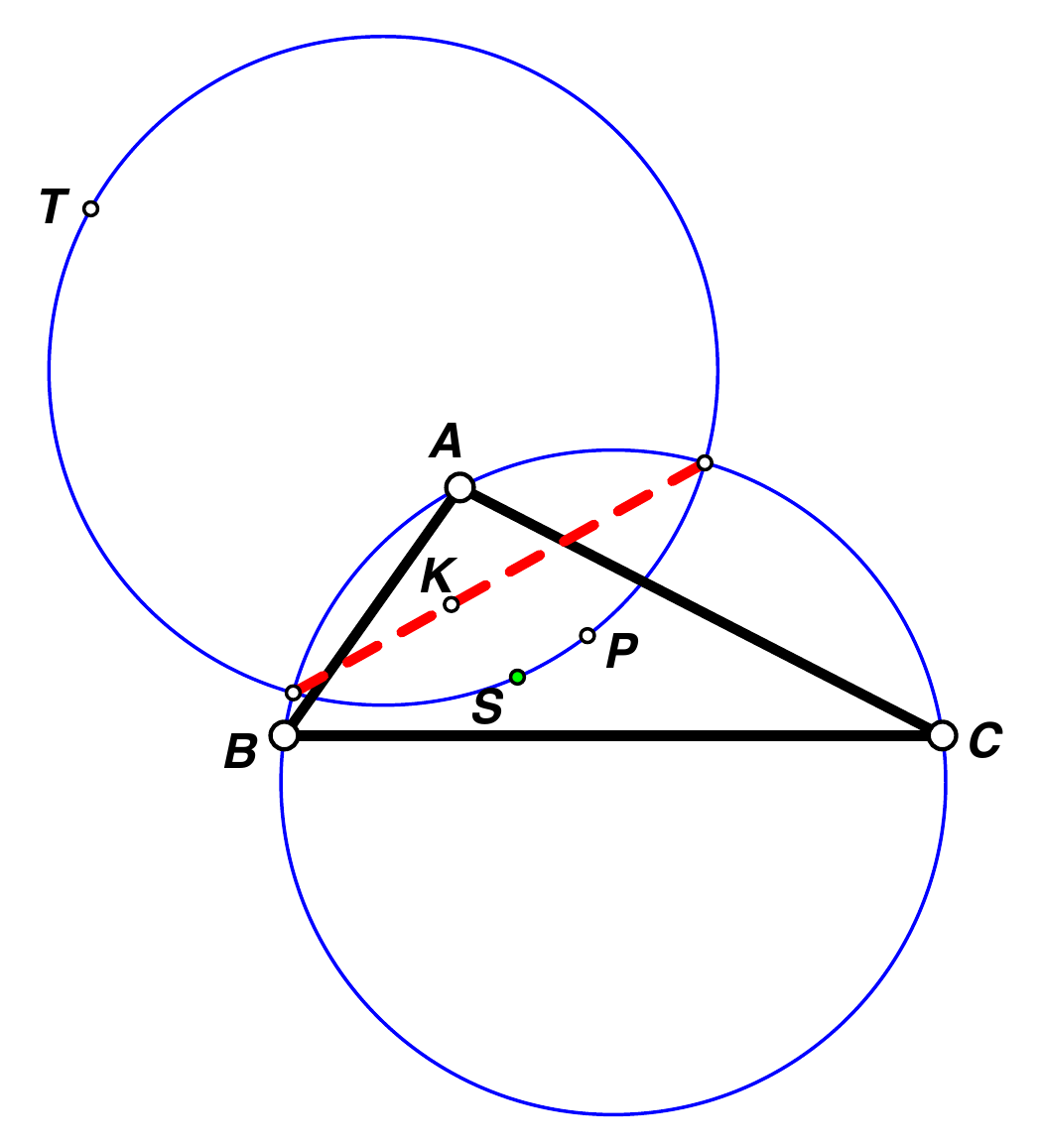}
{$P$ is any point}
{Common chord of $\odot ABC$ and $\odot PST$ passes through $K$.}

\pro{\cite{Jkh-corr}}{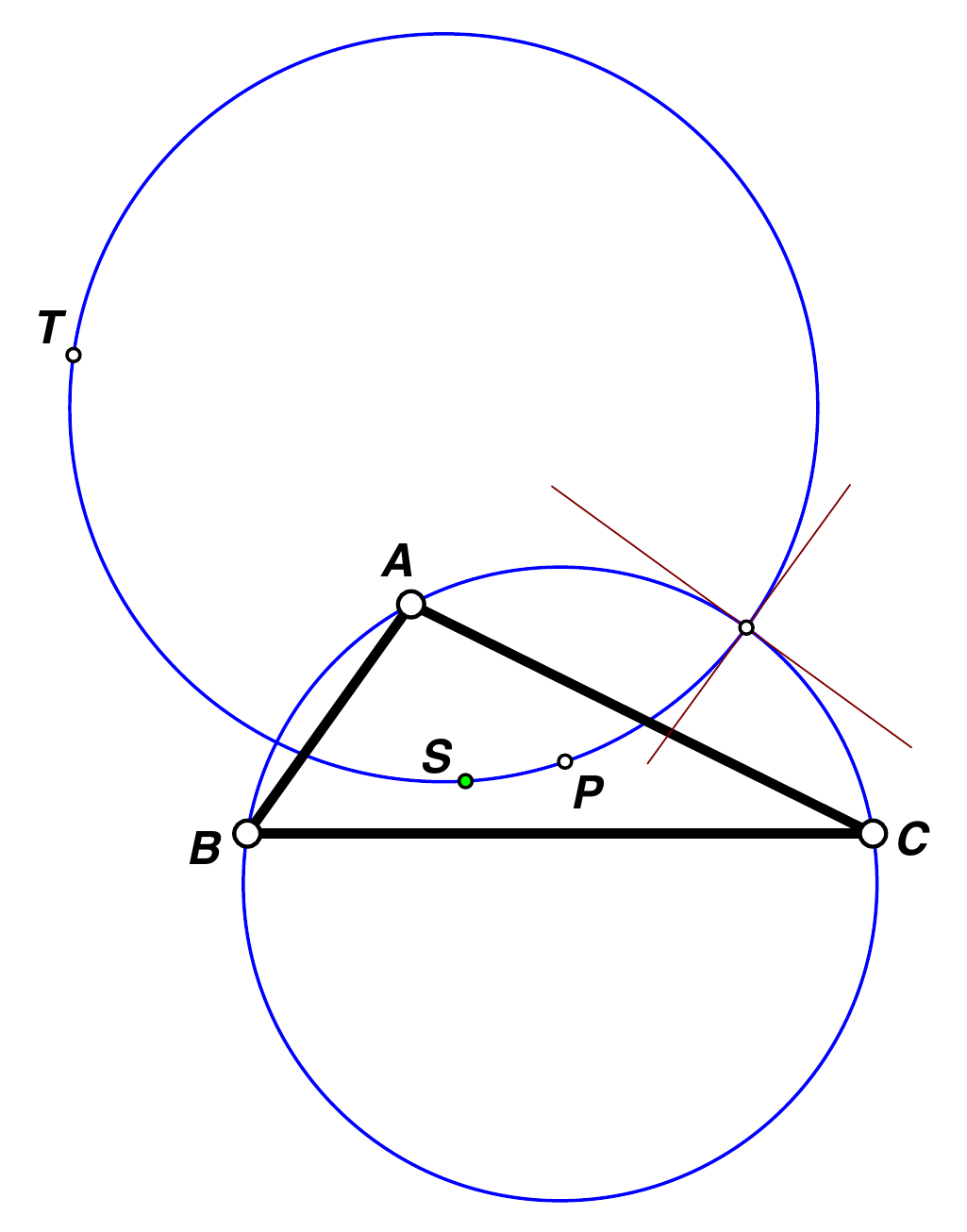}
{$P$ is any point}
{$\odot ABC$ and $\odot PST$ are orthogonal.}
\extra{In particular, the Parry circle is orthogonal to the circumcircle.}

\pro{\cite{Jkh-corr}}{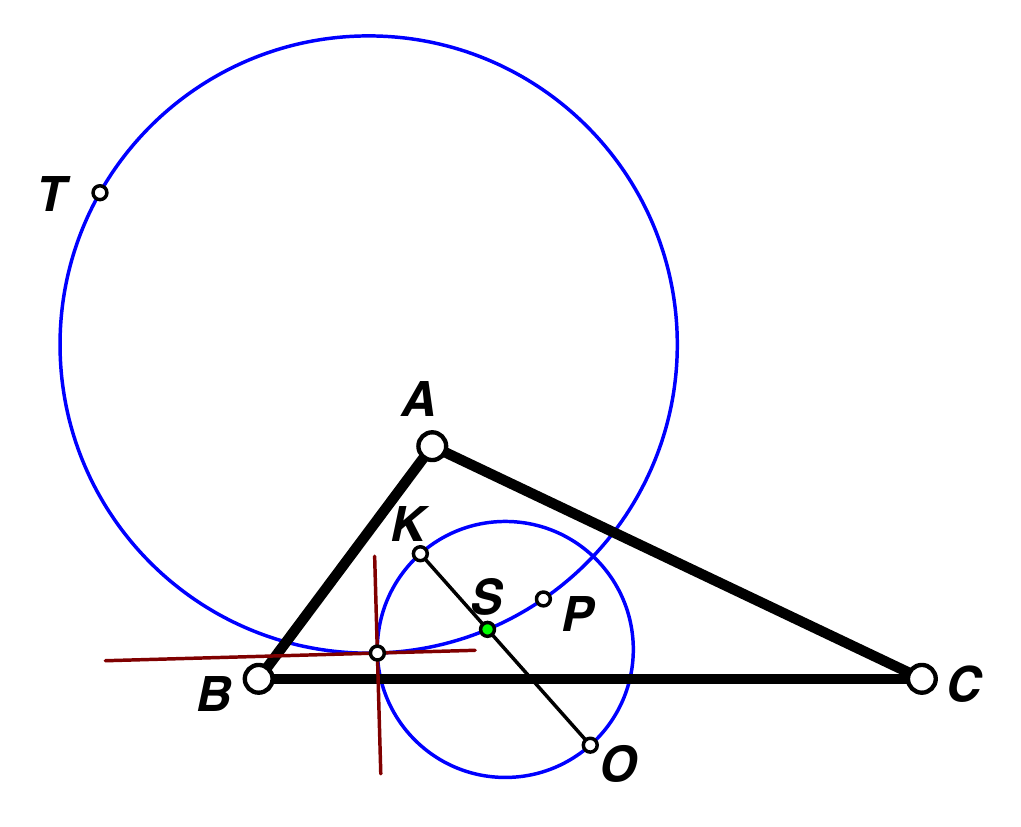}
{$P$ is any point}
{Circle with diameter $KO$ is orthogonal to $\odot PST$.}
\extra{In particular, the Parry circle is orthogonal to the Brocard circle.}

\pro[1st Isodynamic-Dao Triangle]{\cite{ETC16802}}{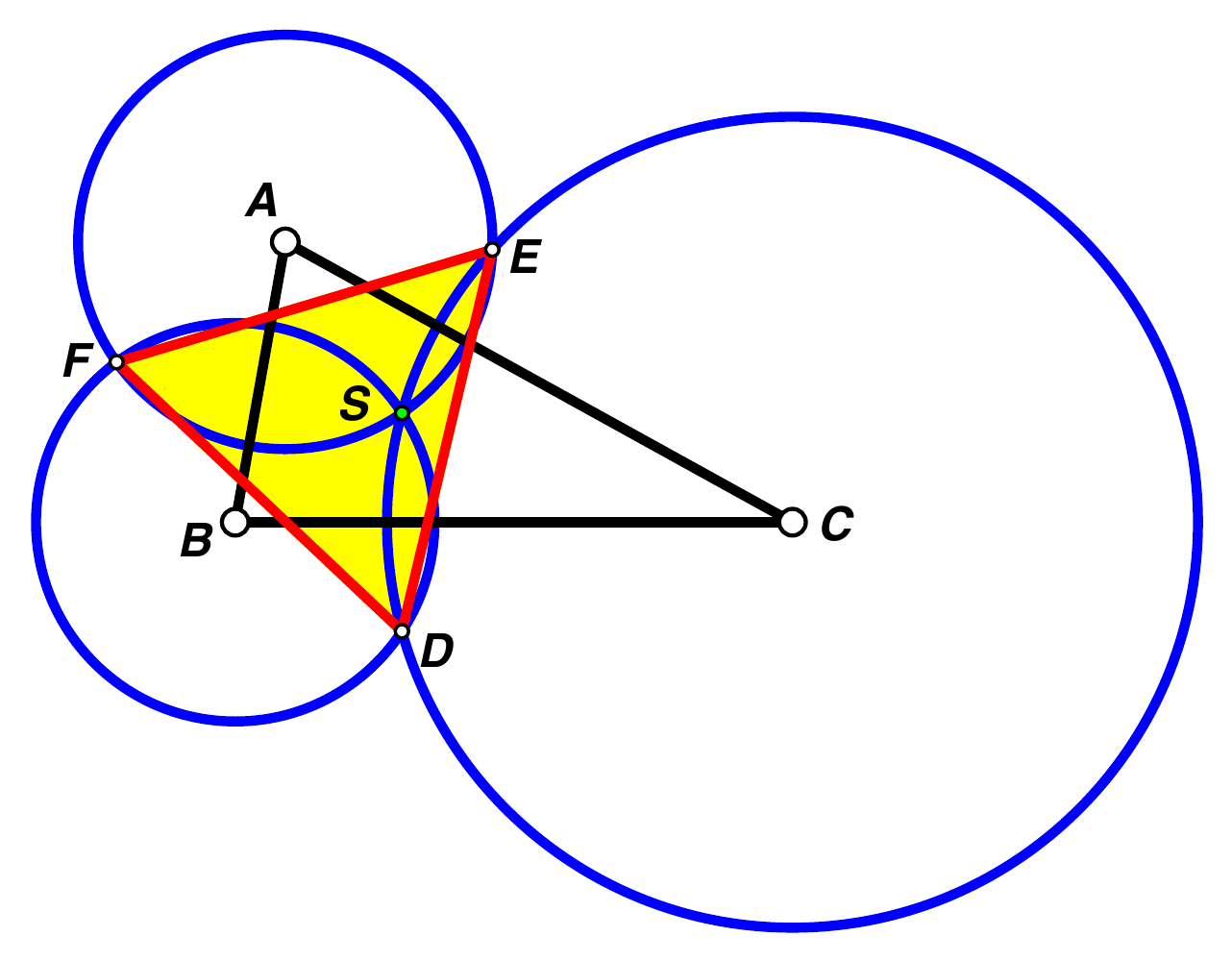}
{}
{$DEF$ is an equilateral triangle.}

\mysubsection{Euler lines}

\pro{\cite{PGR40}}{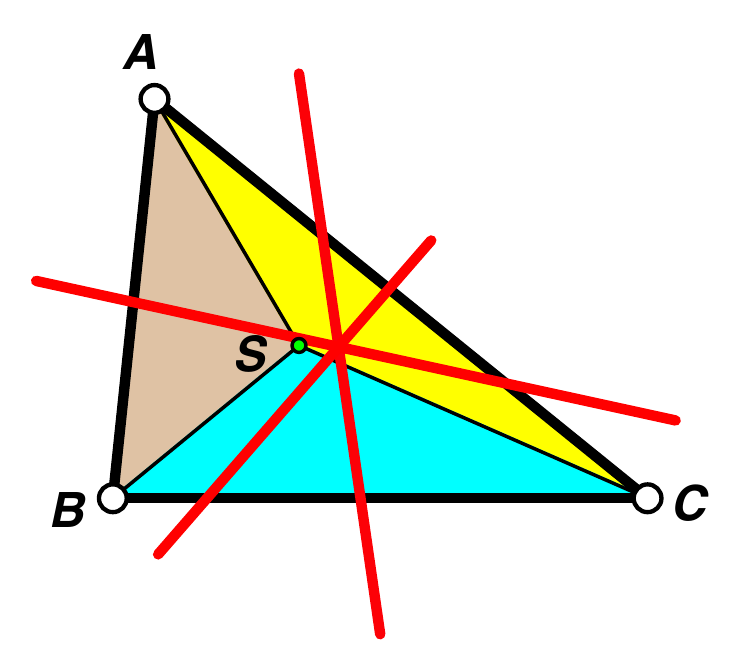}
{}
{The Euler lines of the three colored triangles concur.}
\extra{A computer analysis found that the point of intersection is $X_{61}$.
The three Euler lines will still concur if $S$ is replaced by any point on the circumcircle
or on the Neuberg cubic.
See \cite{MathWorld-Neuberg}.
}

\mysubsection{equilateral triangles}

\nopagebreak
\pro{\cite{MathWorld}}{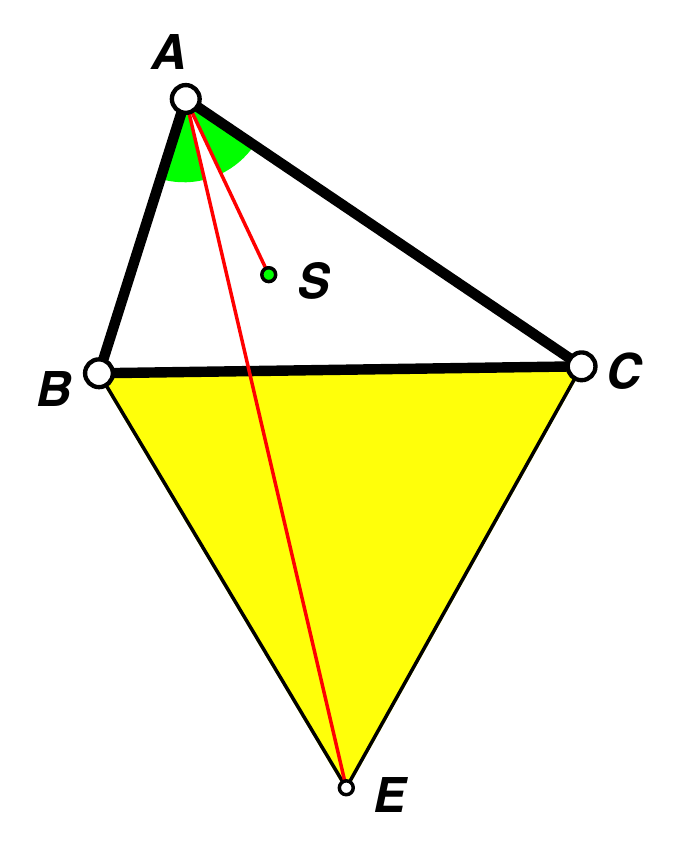}
{\void{equilateral triangle $BCE$}}
{$\angle BAE=\angle SAC$}
\extra{This follows from Property \ref{property-iso-20}.}

\pro{\cite{RG3786}}{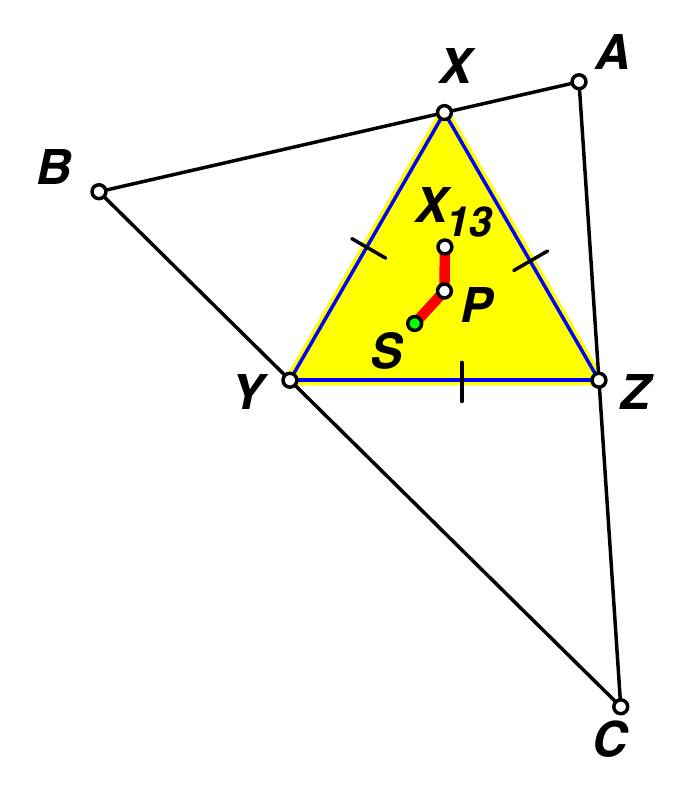}
{$P$ is the center of $\triangle XYZ$.}
{$PS=PX_{13}$}

\pro{\cite{Wiki}}{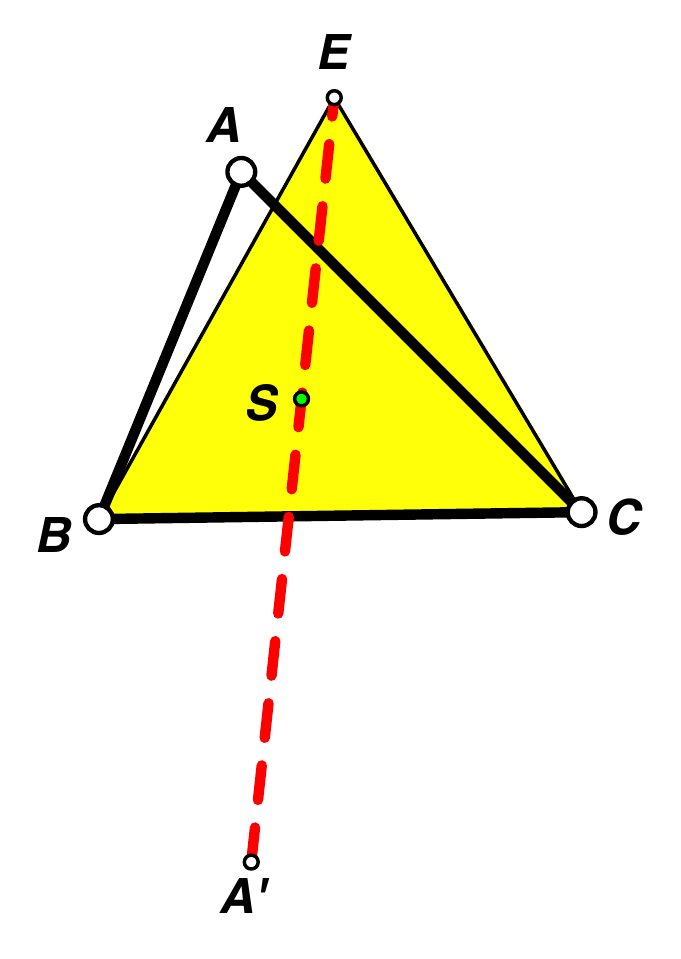}
{$A'$ is the reflection of $A$ about $BC$.}
{$E$ -- $S$ -- $A'$}

\vspace{-12pt}
\mysubsection{reflection}

\nopagebreak
\new
\pro{\cite{RG8847}}{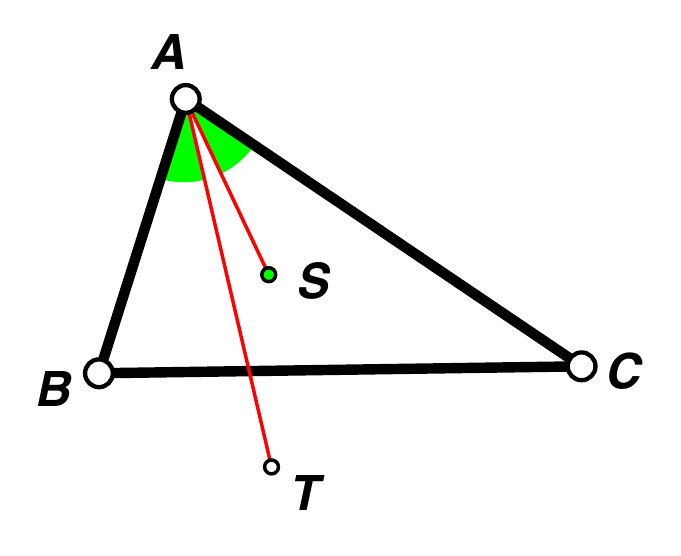}
{$T$ is the reflection of $S$ about $BC$}
{$\angle BAT=\angle SAC$}

\pro{\cite{RG8848}}{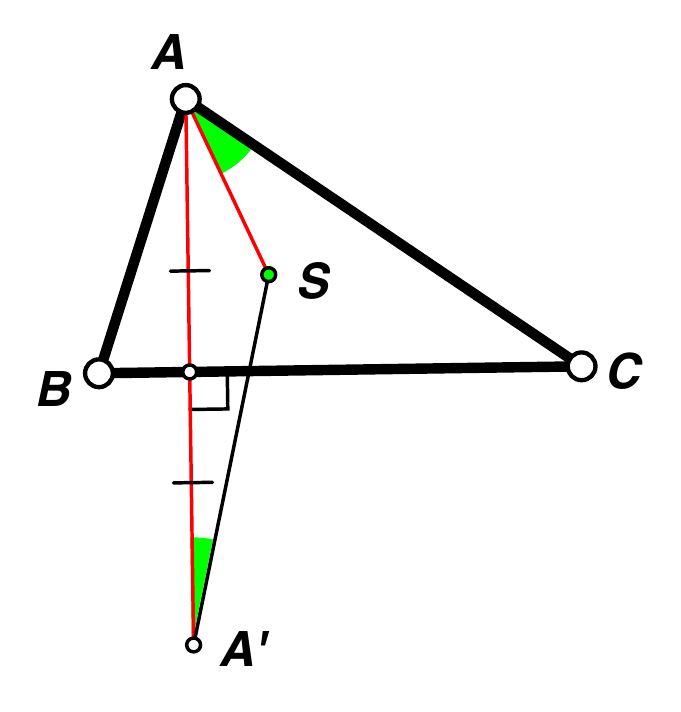}
{$A'$ is the reflection of $A$ about $BC$}
{$\angle SA'A=\angle SAC$}
\extra{The property is true if $S$ is replaced by any point on the $C$-Apollonian circle.}

\mysubsection{sagitta}

\new
\pro{\cite{RG9044}}{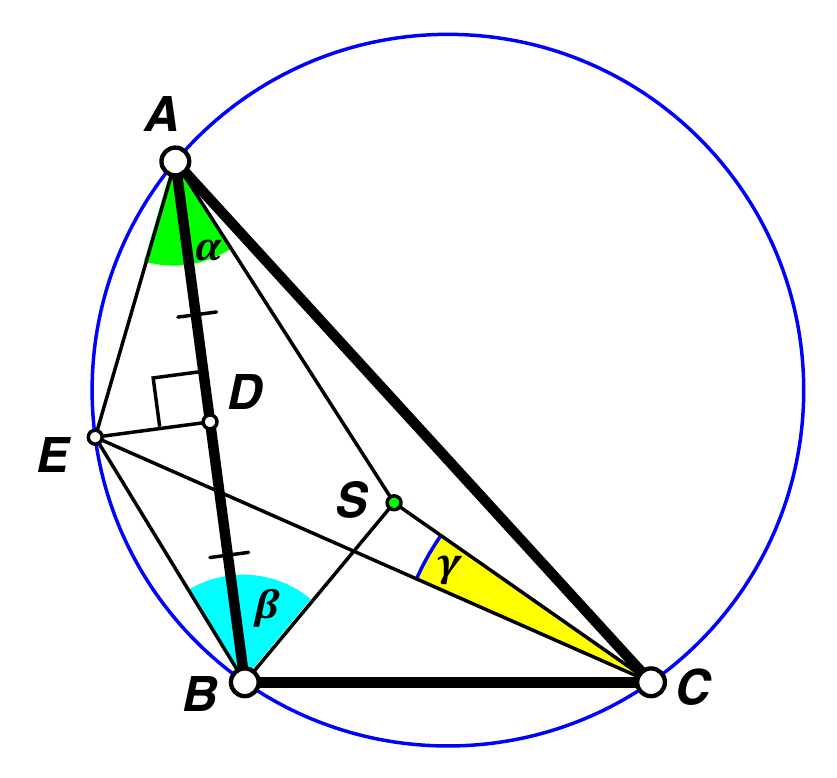}
{}
{$\alpha+\gamma=60\degrees$\\
\then $\beta-\gamma=60\degrees$}

%\newpage
%\vspace*{-24pt}
\mysection{Special triangle plus $S$}

\mysubsection{isosceles triangle}

\new
\pro{\cite{RG4850}}{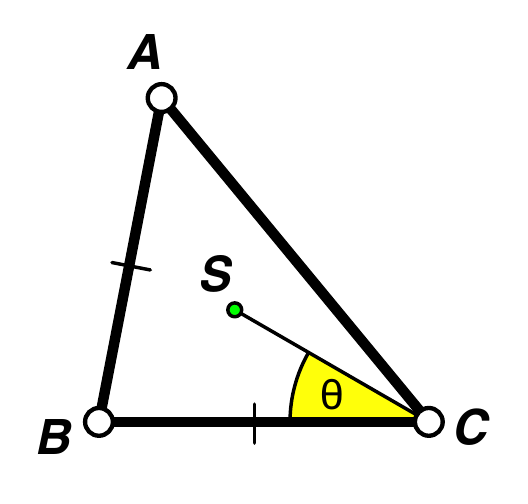}
{}
{$\theta=30\degrees$.}

\mysubsection{isosceles right triangle}

\new
\pro{}{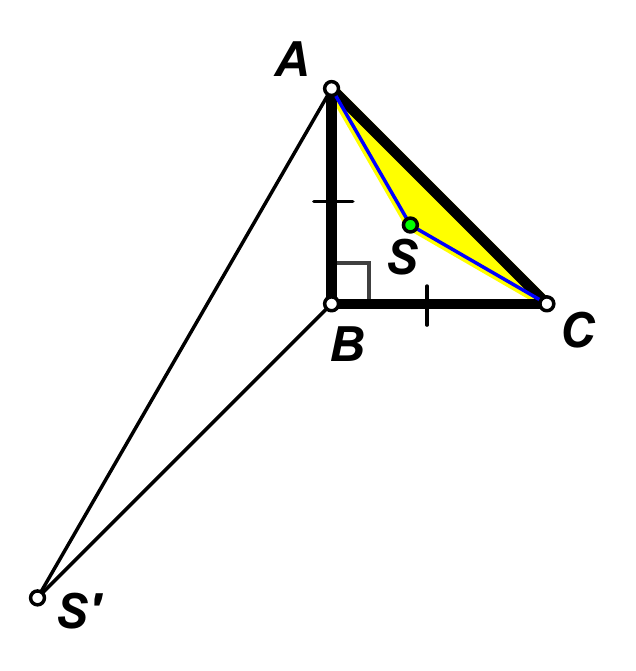}
{}
{$S'$ is the orthocenter of $\triangle ACS$.\\
\then $S$ is the orthocenter of $\triangle ACS'$.}

\mysubsection{$30\degrees$ triangle}

\nopagebreak
\new
\pro{}{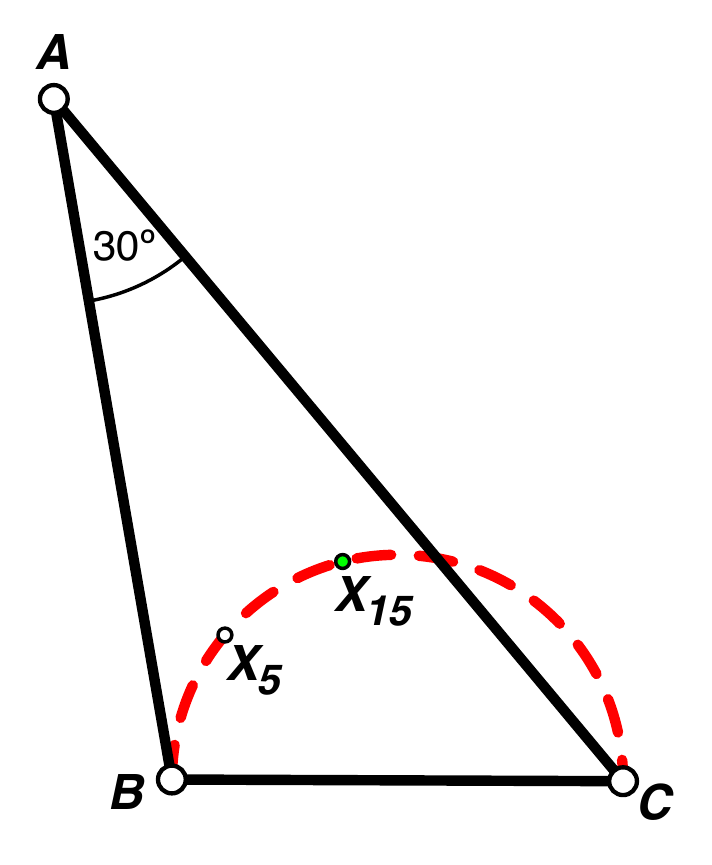}
{}
{$B$, $X_5$, $X_{15}$, and $C$ are concyclic.}

%\mysection{$60\degrees$ triangle plus $S$}

%\mysubsection{area properties}
\mysubsection{$60\degrees$ triangle}

\new
\pro{}{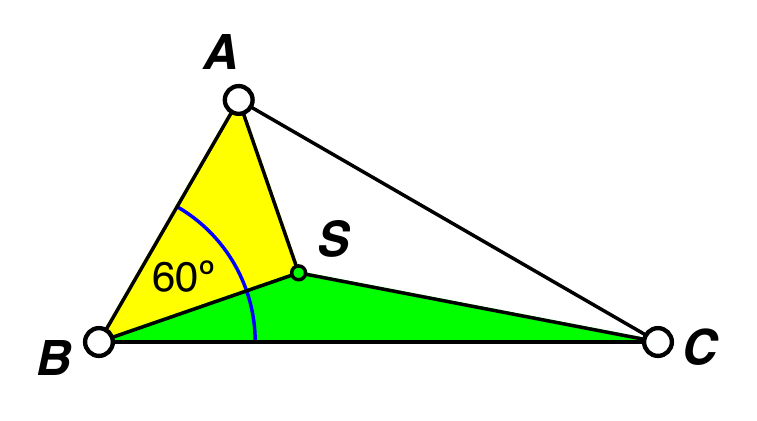}
{}
{$[BAS]=[BSC]$.}

%\mysection{$60\degrees$ triangle plus $S$ and other centers}

%\mysubsection{$S$ and $H$}

\new
\pro{}{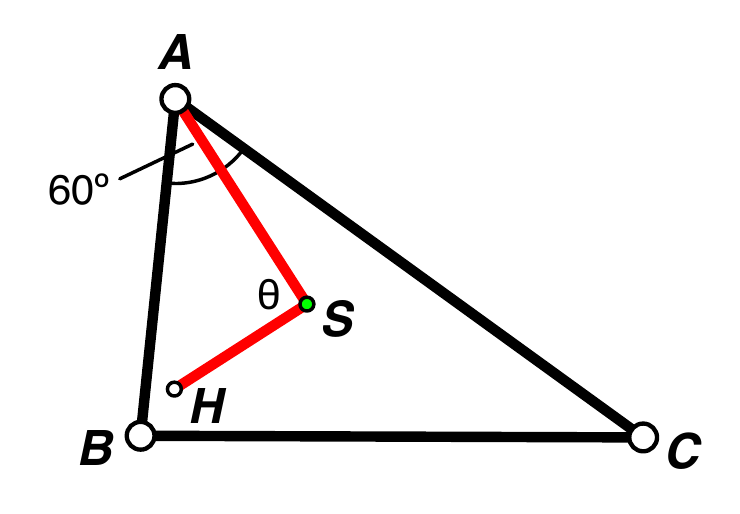}
{}
{$\theta=90\degrees$.}

%\mysubsection{$S$ and $K$}

\new
\pro{}{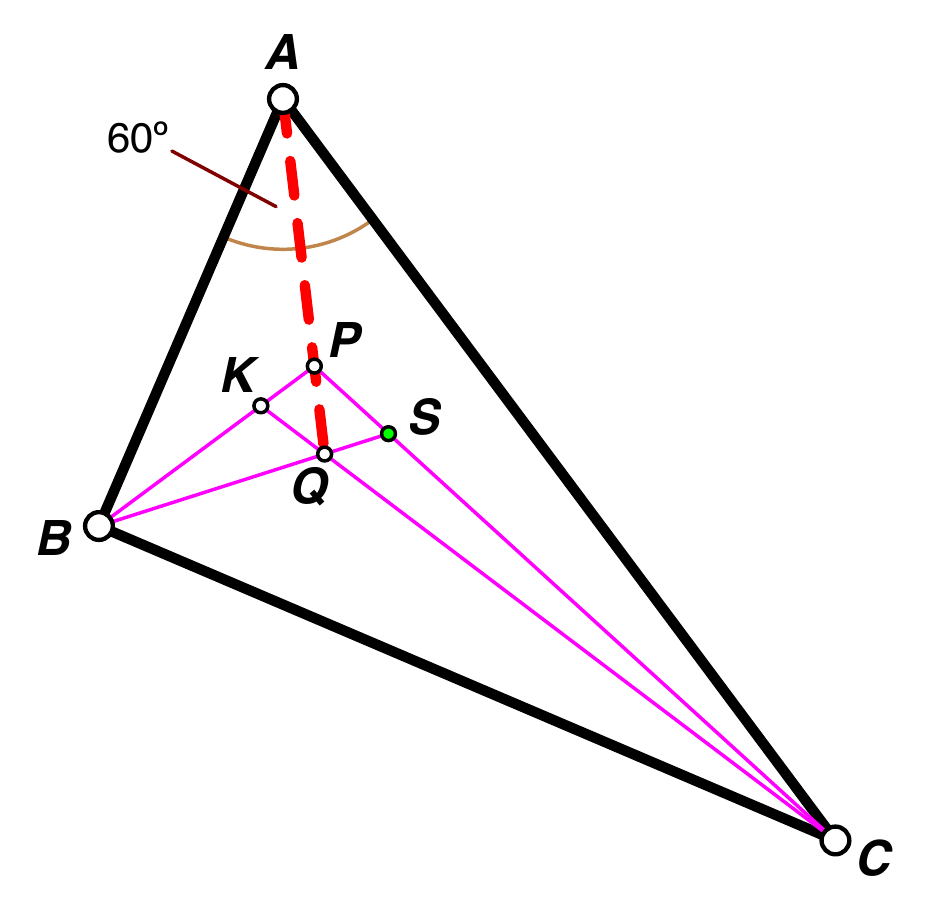}
{}
{$A$ -- $P$ -- $Q$.\\
\then $\angle PAC=30\degrees$.}

%\mysubsection{$S$ and $M$}

\new
\pro{\cite{RG9107}}{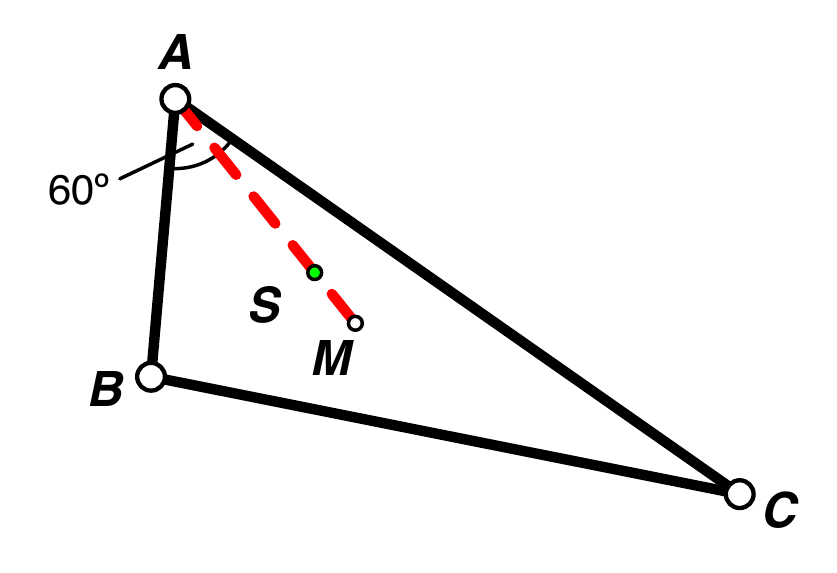}
{}
{$A$ -- $S$ -- $M$.}

%\mysubsection{$S$ and $T$}

\new
\pro{}{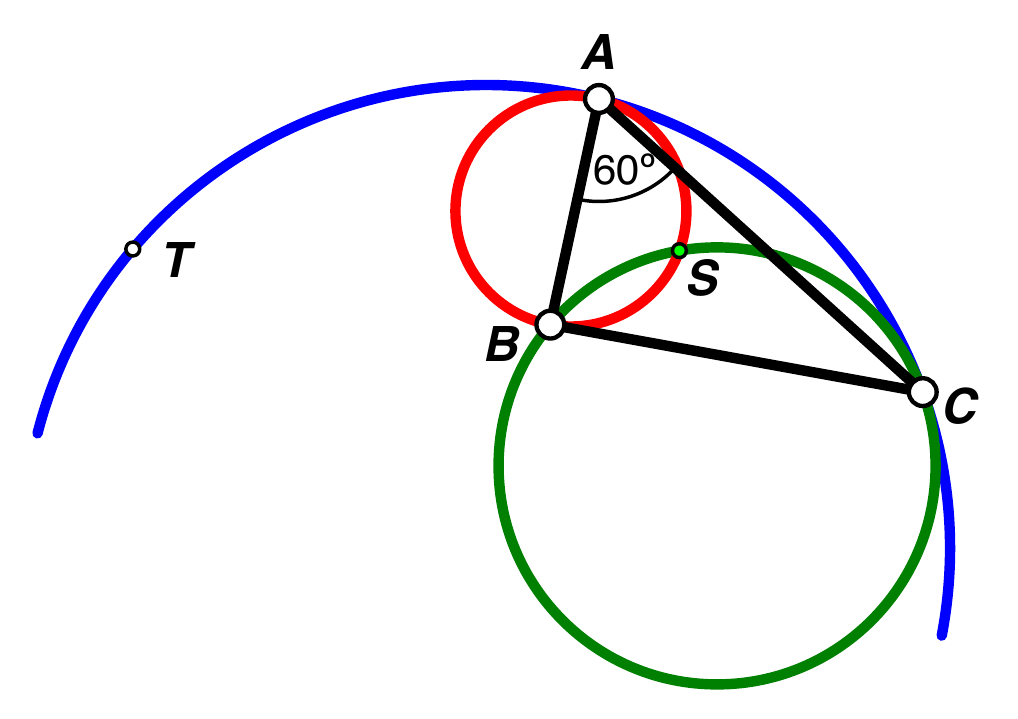}
{}
{$\odot CAT$ and $\odot ABS$ are tangent.\\
\then $\odot CAT$ and $\odot BSC$ are tangent.}

\new
\pro{}{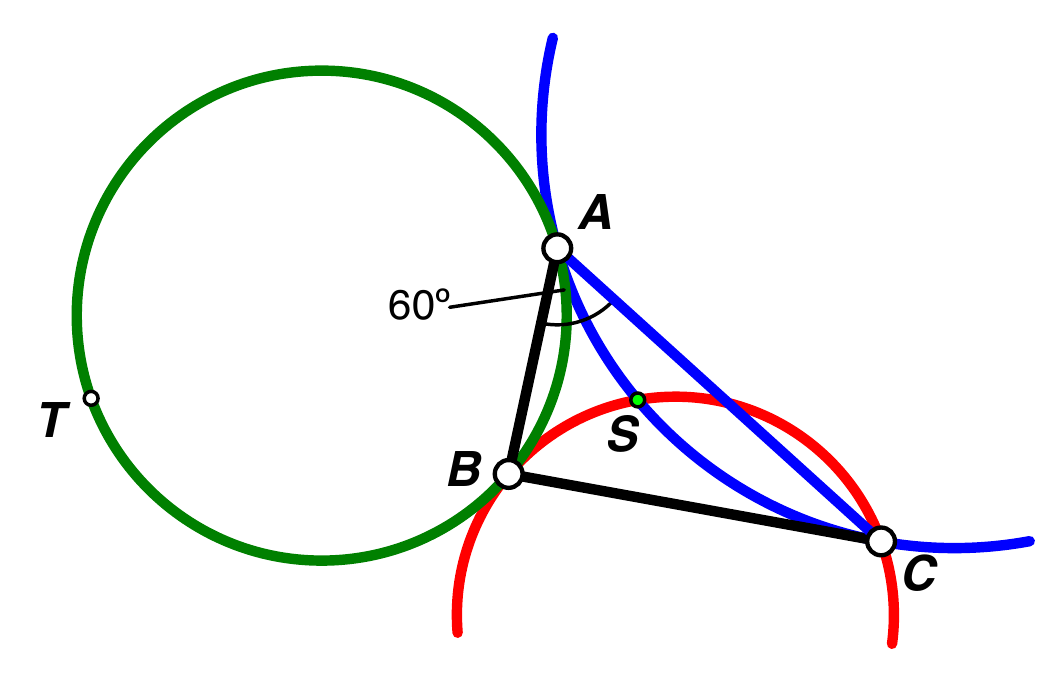}
{}
{$\odot TAB$ and $\odot ASC$ are tangent.\\
\then $\odot TAB$ and $\odot BSC$ are tangent.}

\new
\pro{}{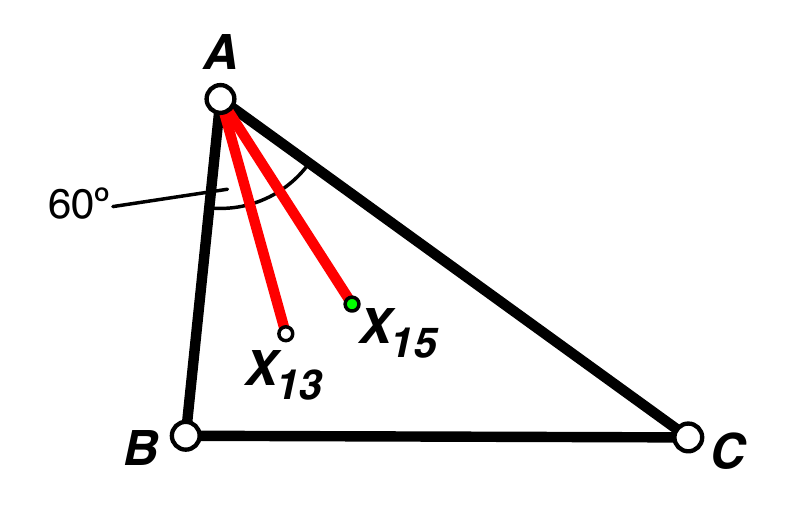}
{}
{$AX_{13}=AX_{15}$.}

\pro{\cite{Titu}}{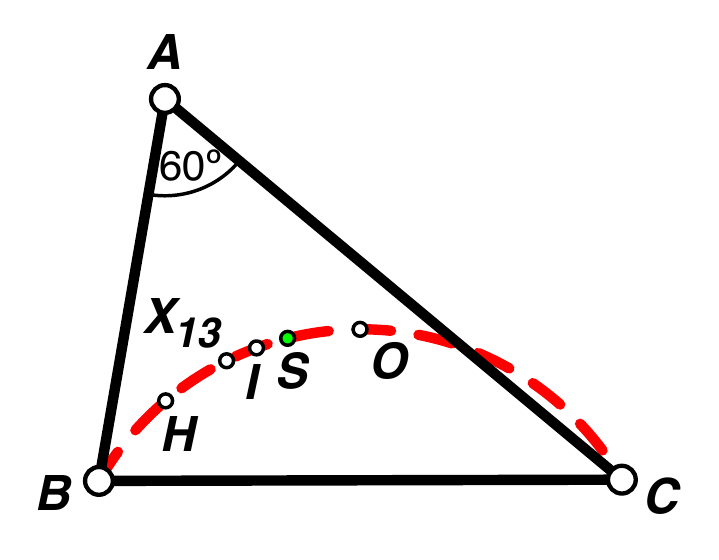}
{}
{$B$, $H$, $X_{13}$, $I$, $S$, $O$, and $C$ are concyclic.}
\extra{See also \cite{Moon}. Other points on this circle (discovered by computer) are $X_{399}$, $X_{616}$, and $X_{617}$.}

\mysubsection{$120\degrees$ triangle}

\new
\pro{}{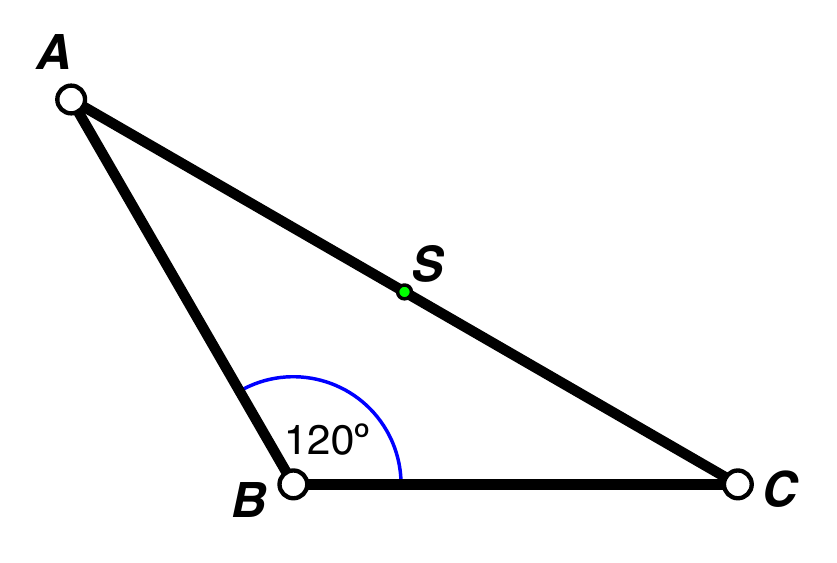}
{\void{$\angle B=120\degrees$.}}
{$S$ lies on $AC$.}

\new
\pro{}{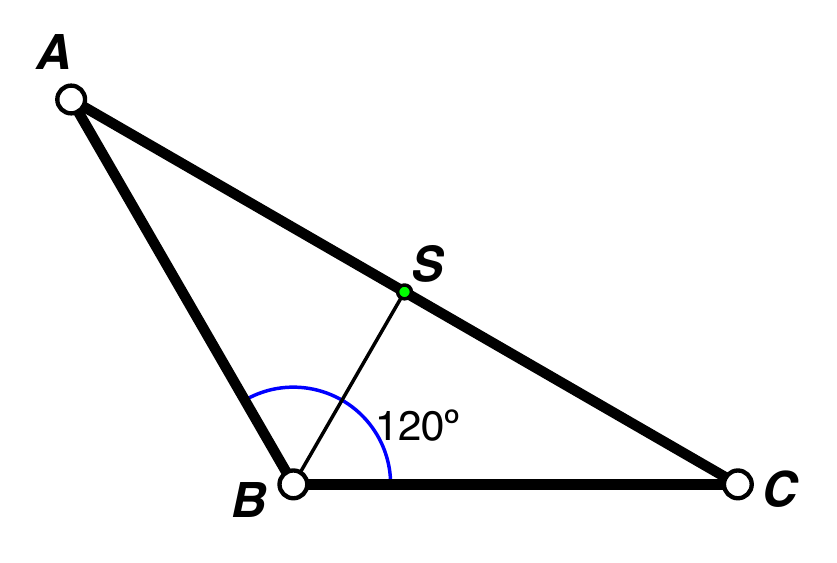}
{\void{$\angle B=120\degrees$.}}
{$BS$ bisects $\angle CBA$.}

\mysubsection{$150\degrees$ triangle}

\nopagebreak
\new
\pro{}{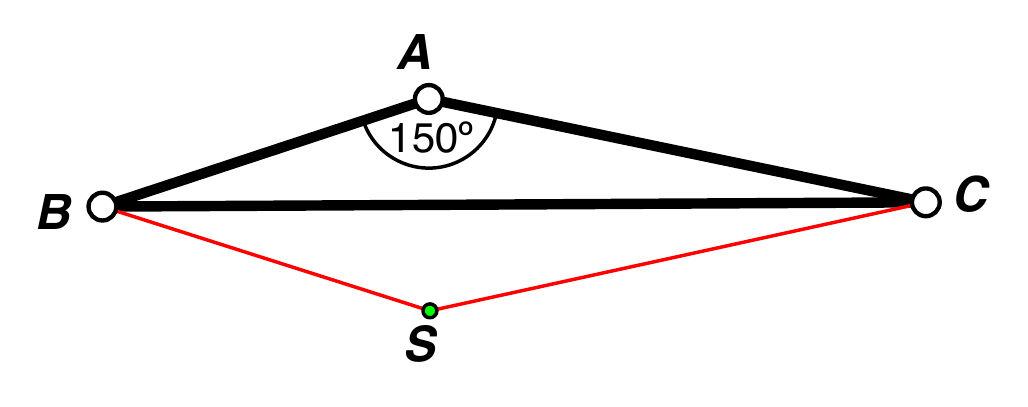}
{}
{$ABSC$ is a kite,.}

\new
\pro{}{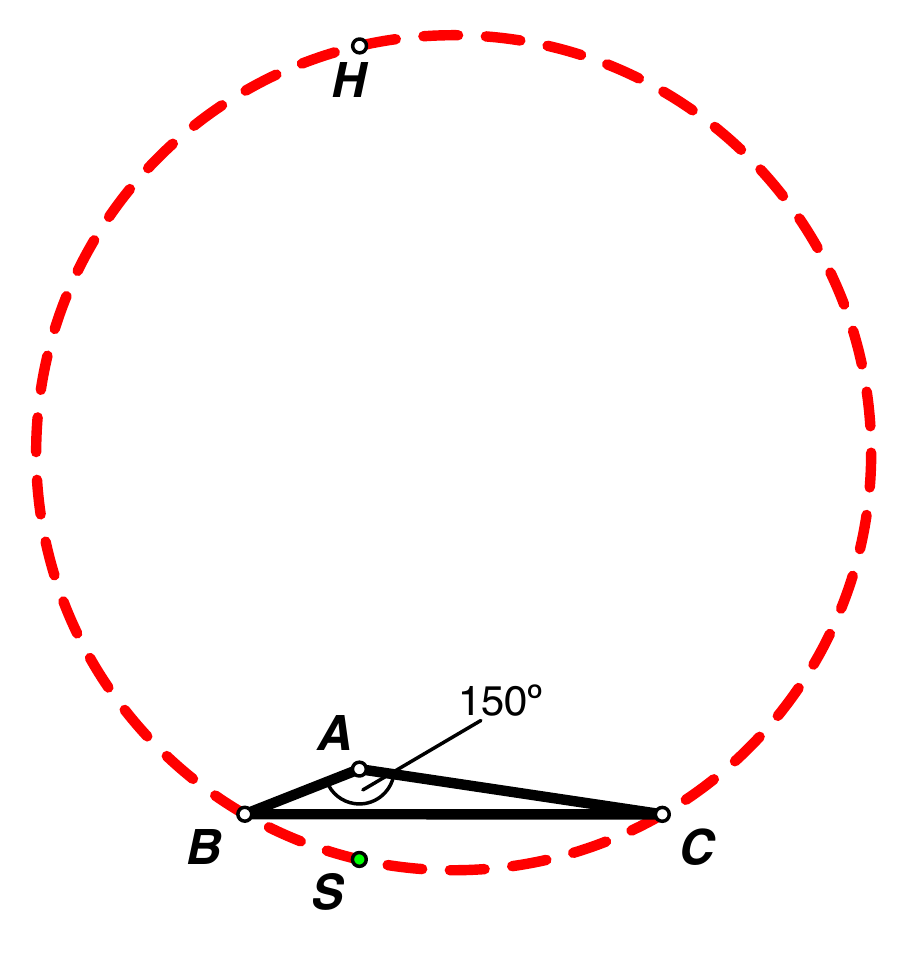}
{}
{$B$, $S$, $C$, and $H$ are concyclic.}

\mysection{Triangle with $S$ and other points}

\mysubsection{$S$ and $T$}

\pro{\cite{Fettis}}{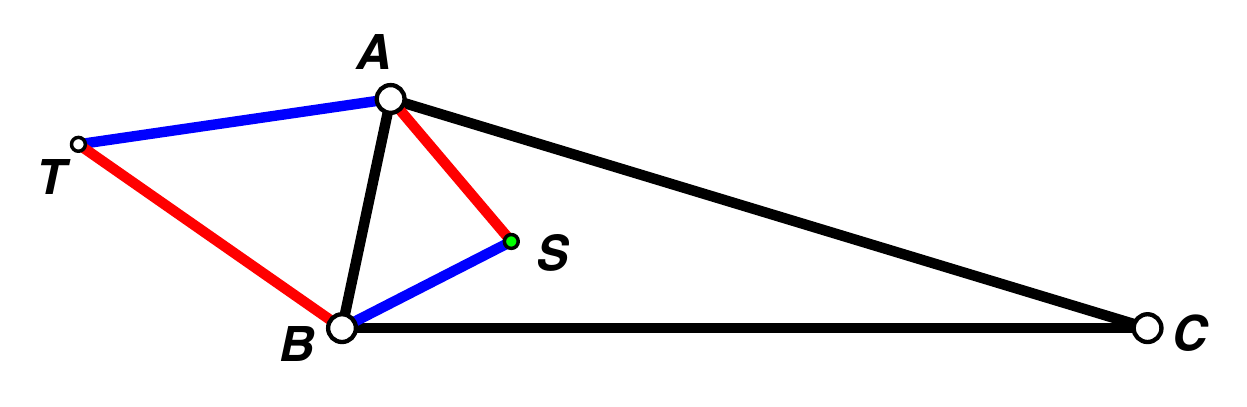}
{}
{$\displaystyle \frac{SA}{TA}=\frac{SB}{TB}$}

\new
\pro{}{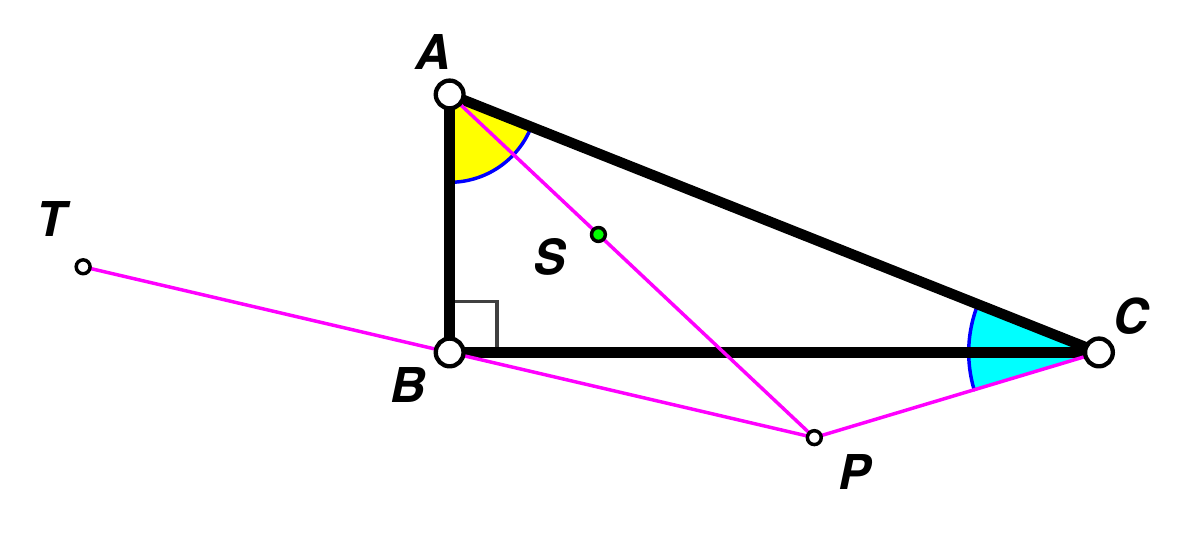}
{}
{$\angle BAC-\angle ACP=30\degrees$}

\pro[Distance Between Isodynamic Points]{\cite{Fettis}}{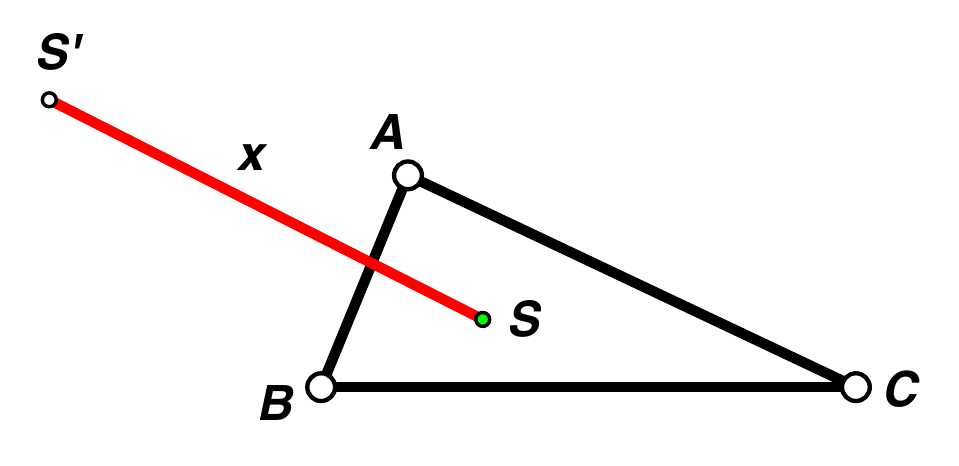}
{}
{$\displaystyle x=\frac{2\sqrt3abc}{\sqrt{(a^2+b^2+c^2)^2-48K^2}}$}

\mysubsection{$S$ and one notable point}

\new
\pro{}{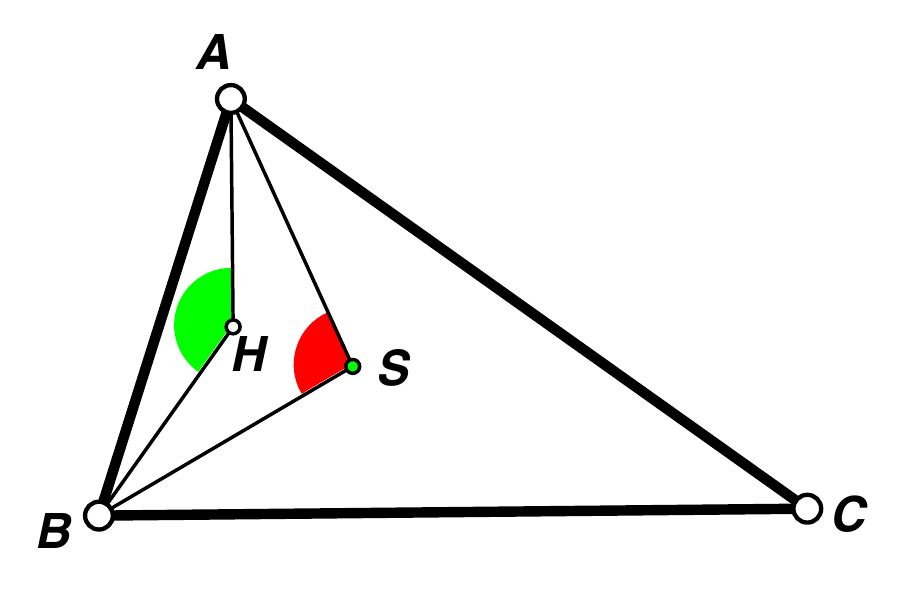}
{$S$ and $H$ lie inside $\triangle ABC$.}
{$\angle AHB+\angle ASB=240\degrees$}

\new
\pro{}{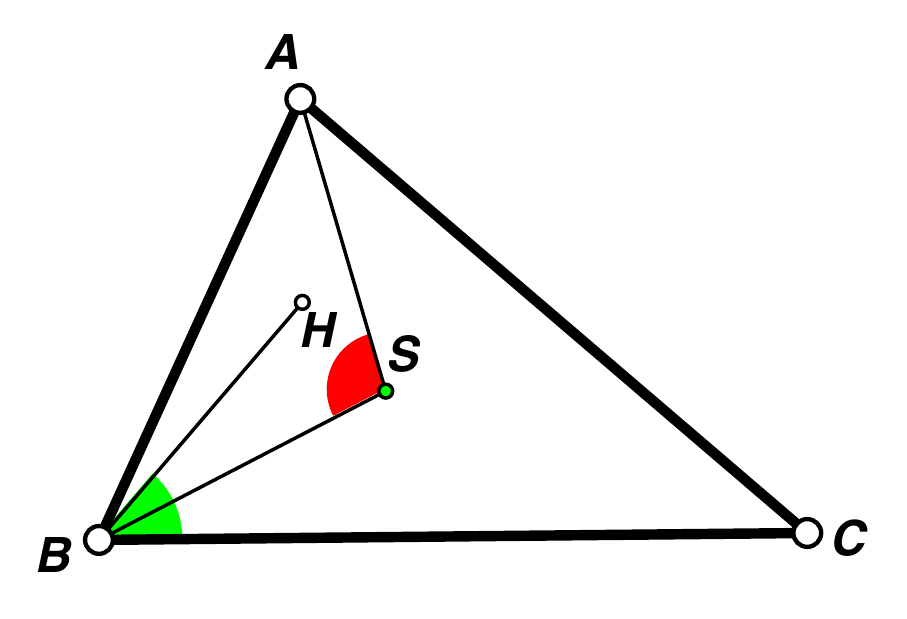}
{$S$ and $H$ lie inside $\triangle ABC$.}
{$\angle CBH+\angle ASB=150\degrees$}

\new
\pro{}{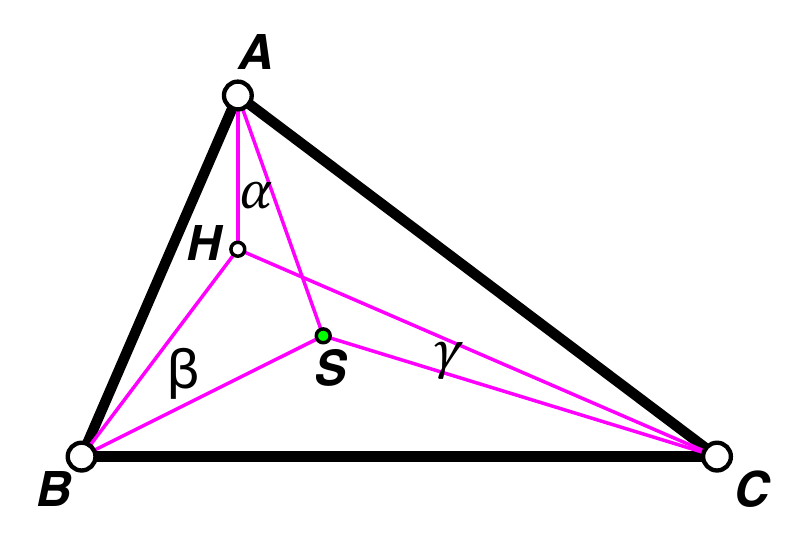}
{$\angle HAS=\alpha$, $\angle HBS=\beta$, and $\angle HCS=\gamma$}
{Using directed angles,
$$\sin\alpha+\sin\beta+\sin\gamma=0.$$}

\new
\pro{}{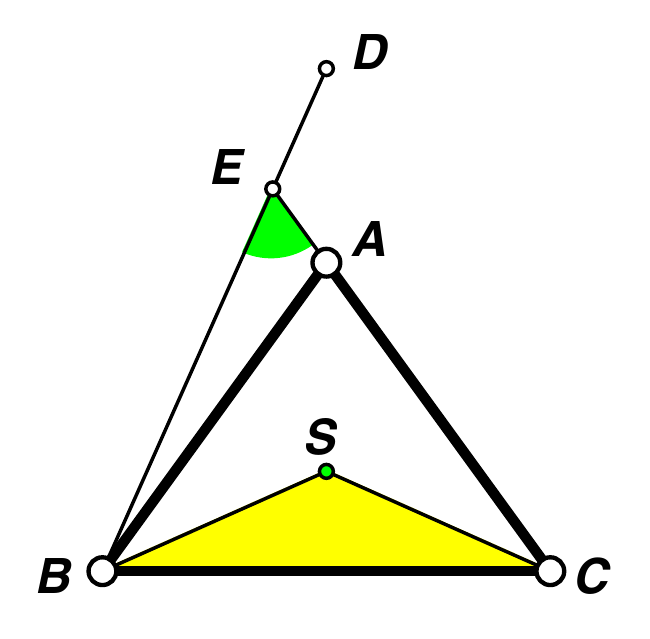}
{$D$ is the orthocenter of $\triangle BCS$.}
{$\angle BEC=60\degrees$}

\new
\pro{}{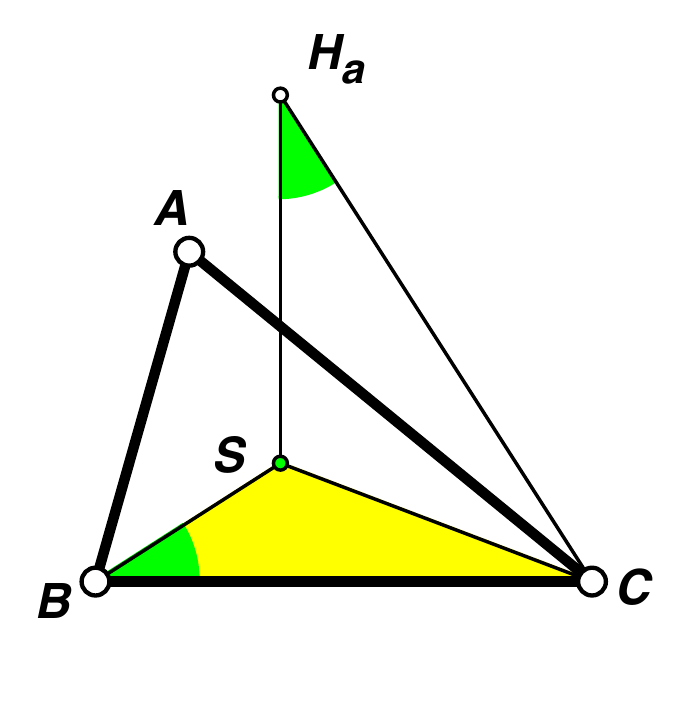}
{$H_a=X_4(BCS)$}
{$\angle CBS=\angle SH_aC$}

\new
\pro{}{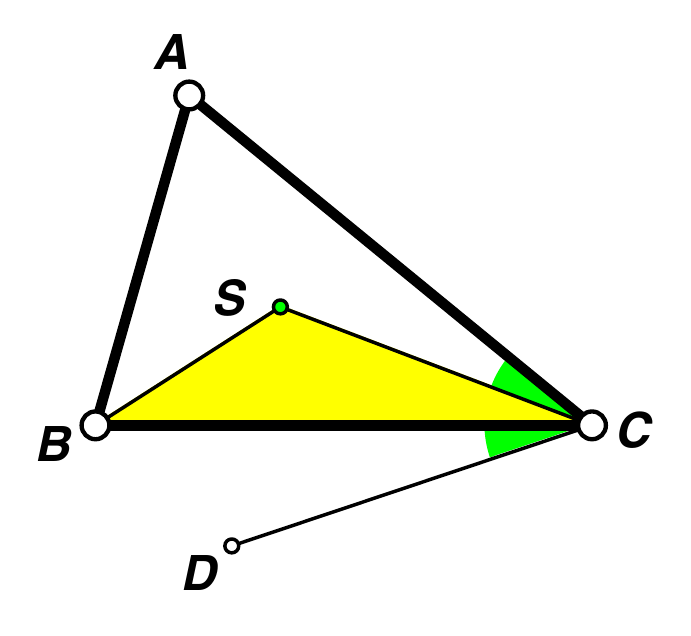}
{$D=X_{14}(BCS)$}
{$\angle BCD=\angle ACS$}

\new
\pro{}{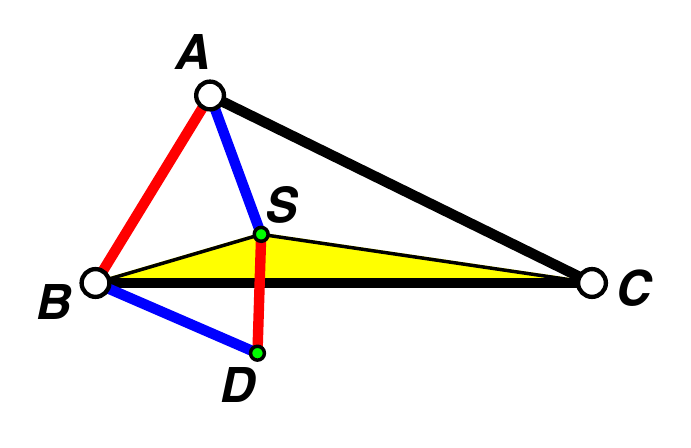}
{$D=X_{15}(BCS)$}
{$AB\cdot SD=AS\cdot BD$}

\new
\pro{}{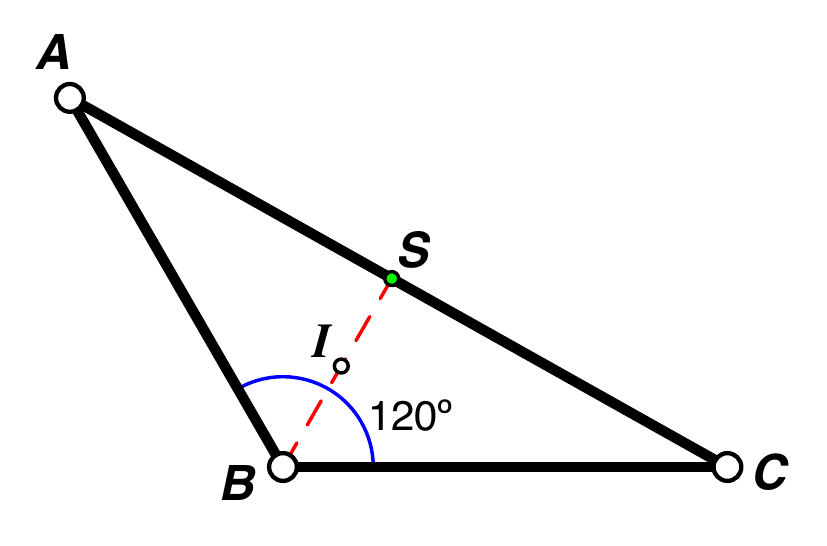}
{\void{$\angle B=120\degrees$.}}
{$B$ -- $I$ -- $S$}
\extra{This result follows from Property \ref{property-iso-30}.}

\new
\pro{\cite{Rabinowitz-A}}{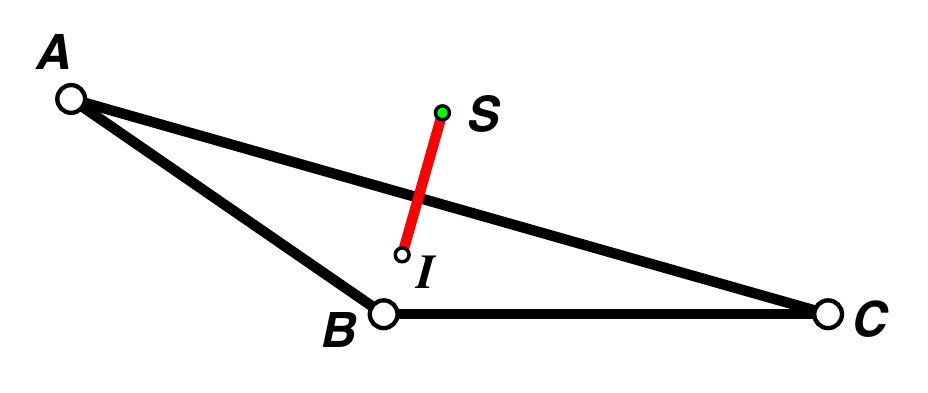}
{}
{$SI<\frac14 R$.}
\extra{The smallest $k$ such that $SI<kR$ is
$k\approx 0.2370406267$, where $k$ is the positive root of
$4 x^6+36 x^5+120 x^4+288 x^3+513 x^2-72 x-16$.}

\pro[Isogonal Conjugate]{\cite{MathWorld}}{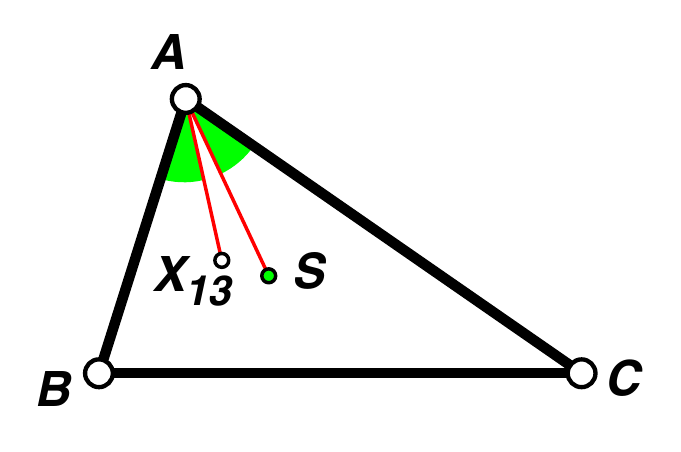}
{}
{$\angle BAX_{13}=\angle SAC$}
\extra{In other words, the isogonal conjugate of $X_{15}$ is $X_{13}$.
See \cite[p.~296]{Johnson}.
The isotomic conjugate of $X_{15}$ is $X_{300}$.
See\cite{ETC15}.
}

\mysubsection{$S$ and $T$ and one notable point}

\nopagebreak
\pro{\cite{Fettis}}{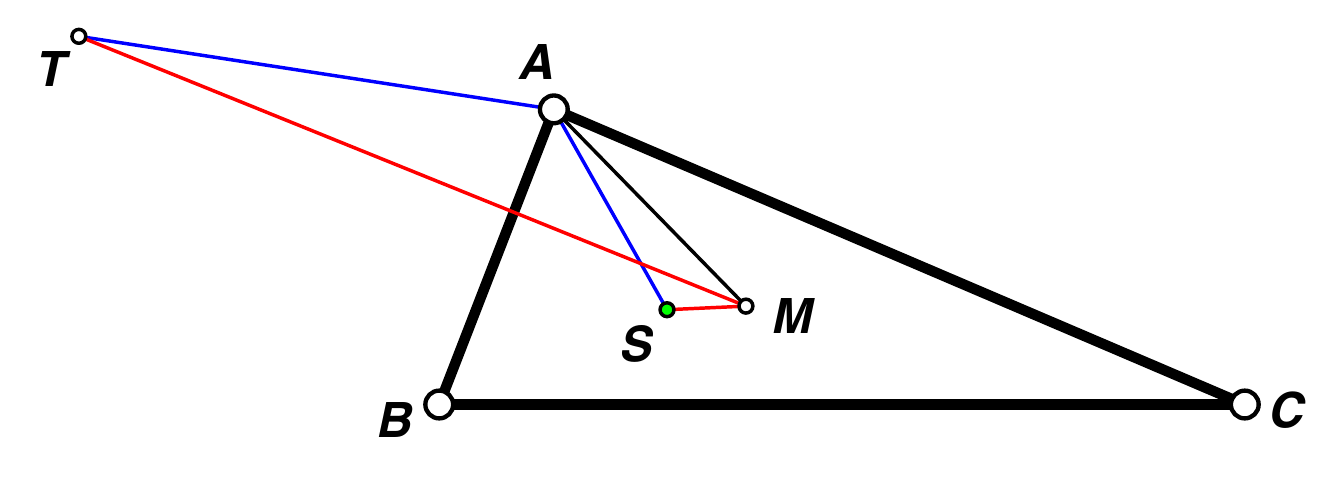}
{}
{$\displaystyle\left(\frac{SA}{TA}\right)^3=\frac{SM}{TM}$}

\new
\pro{\cite{Euclid3368}}{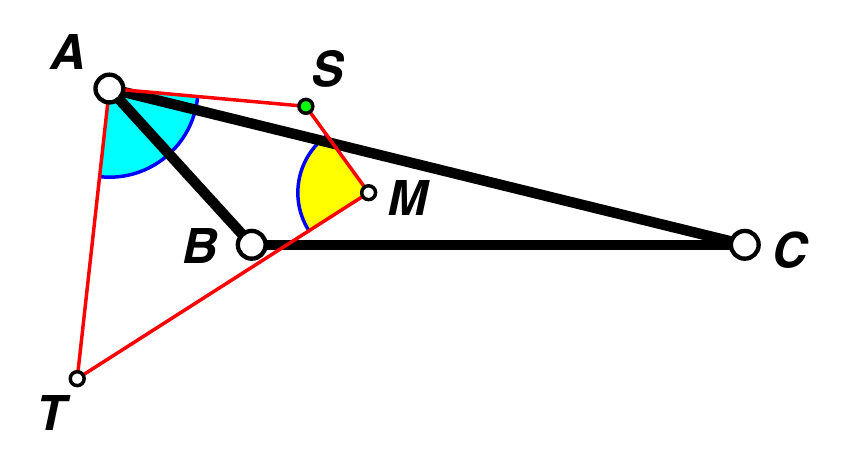}
{$x=\angle TMS$, $\alpha=\angle TAS$}
{$x=3\alpha\pmod{2\pi}$}

\new
\pro{\cite{Euclid3356}}{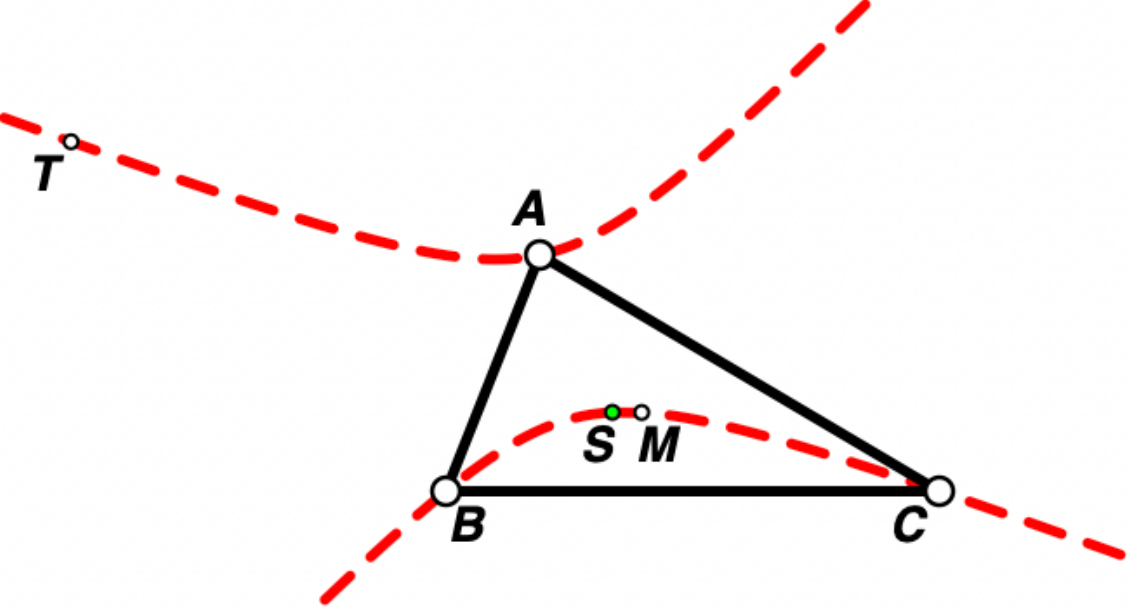}
{}
{$A$, $B$, $C$, $S$, $T$, $M$ lie on a hyperbola.}
\extra{The center of the hyperbola is $X_{18334}$.}

\pro[Inverse Property]{\cite[p.~296]{Johnson}}{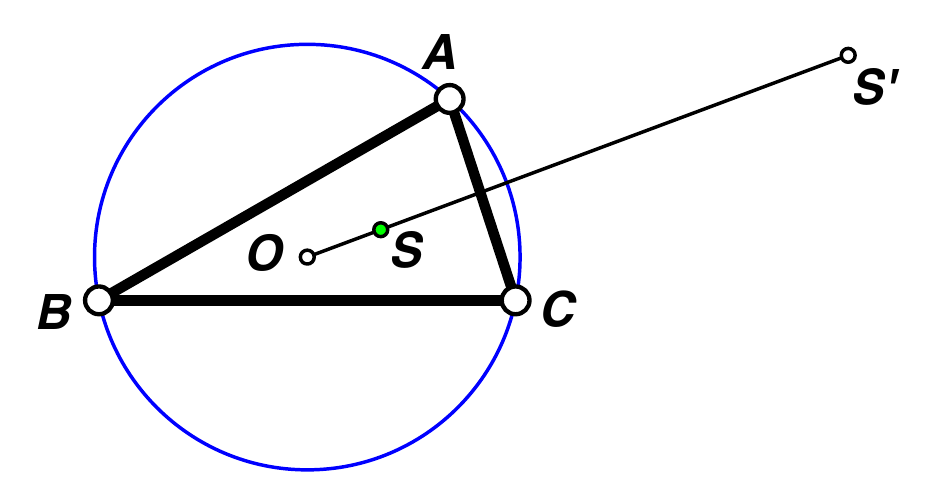}
{}
{$O$ -- $S$ -- $S'$\\
\then {$OS\cdot OS'=R^2$}}
\extra{Inversion of $\triangle ABC$ with respect to an isodynamic point transforms $\triangle ABC$ into an equilateral triangle.}

\new
\pro{}{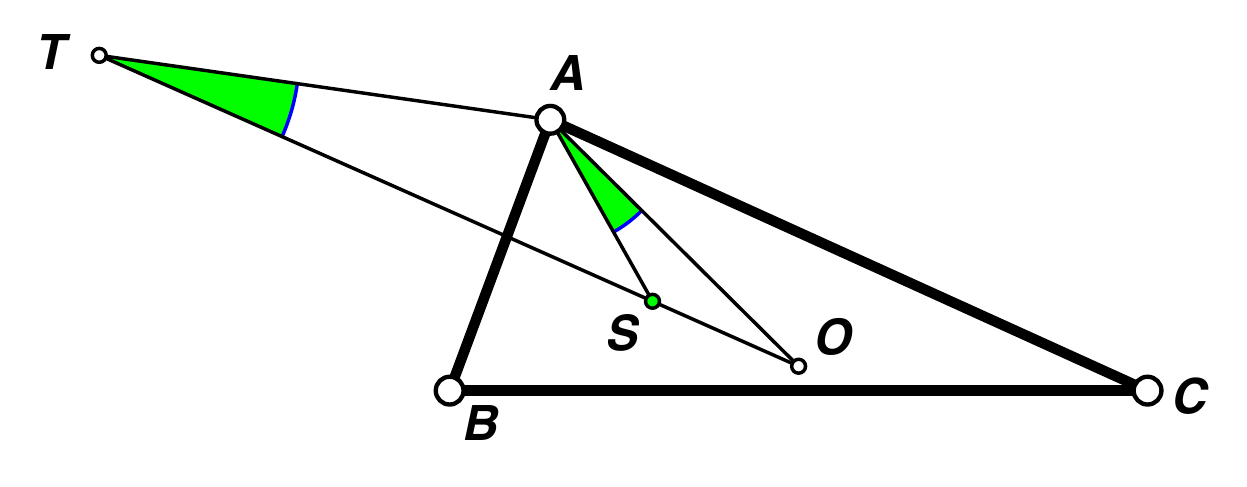}
{}
{$\triangle OSA\sim\triangle TAO$}

\new
\pro{\cite{RG9101}}{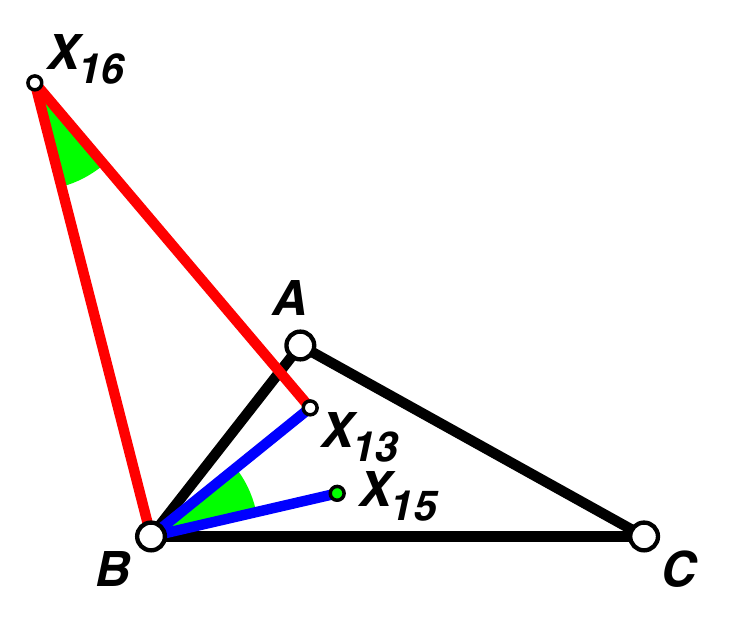}
{}
{$\angle BX_{16}X_{13}=\angle X_{15}BX_{13}$}

\new
\pro{\cite{RG9102}}{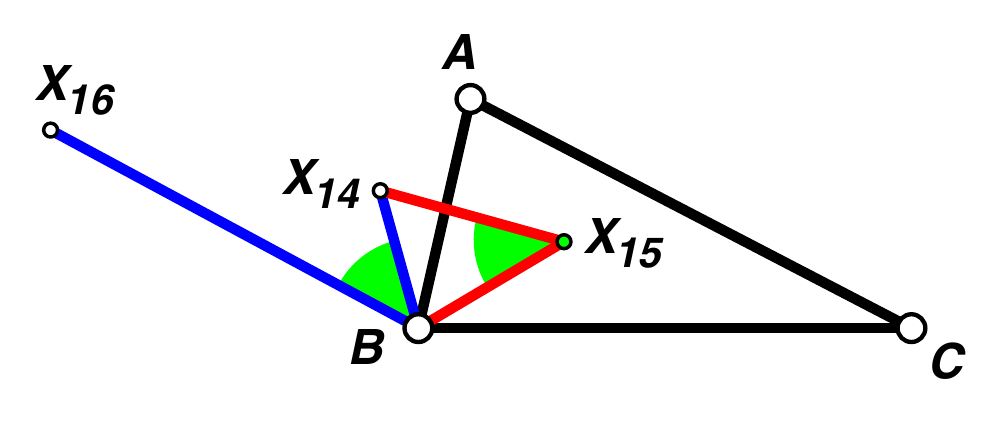}
{}
{$\angle X_{14}BX_{16}=\angle X_{14}X_{15}B$}

\mysubsection{$S$ and $T$ and a point $P$}

\pro{\cite{Fettis}}{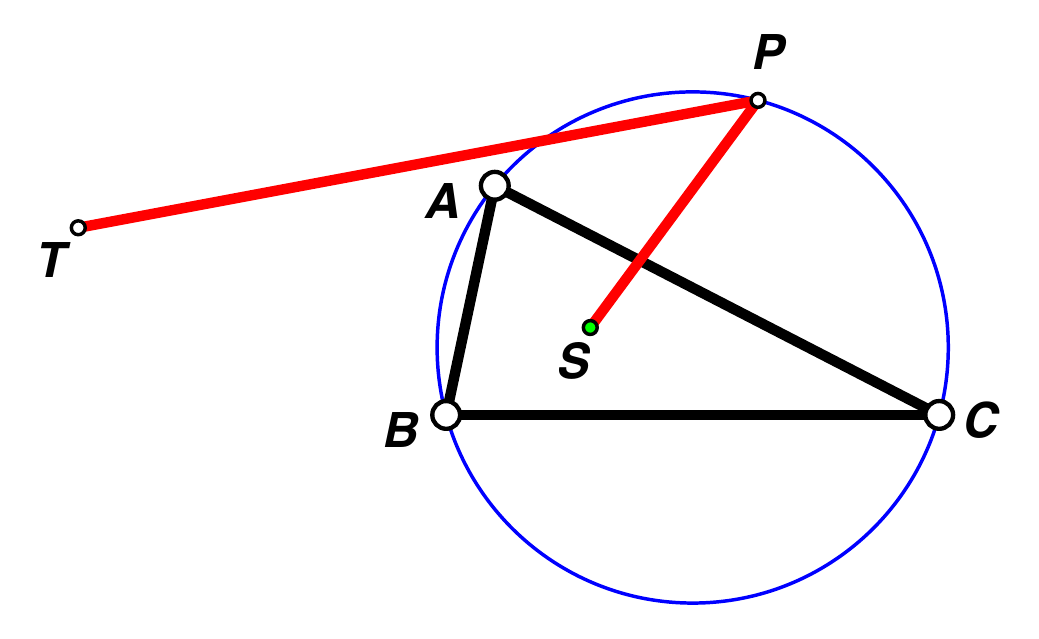}
{}
{$PS/PT$ remains invariant as $P$ moves along the circumcircle.}
\extra{The common ratio is $MX_{13}/MX_{14}$.}

\mysubsection{$S$ and two notable points}

\pro{\cite[Art.~602]{Johnson}}{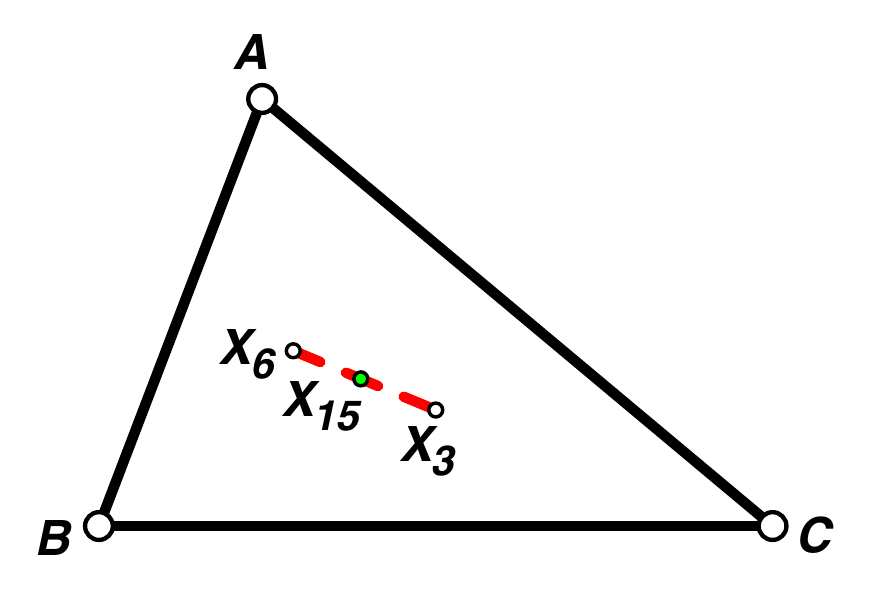}
{}
{$X_3$ -- $X_{15}$ -- $X_6$\\
\then $\displaystyle (X_3X_{15})^2\cdot X_3X_6=R^2(X_3X_{15}-X_6X_{15})$.}
\extra{See also \cite[p.~106]{Gallatly}.
Equation was discovered by computer.
Corollary: $X_3X_{15}> X_6X_{15}$.
}

\pro{\cite{Dekov}}{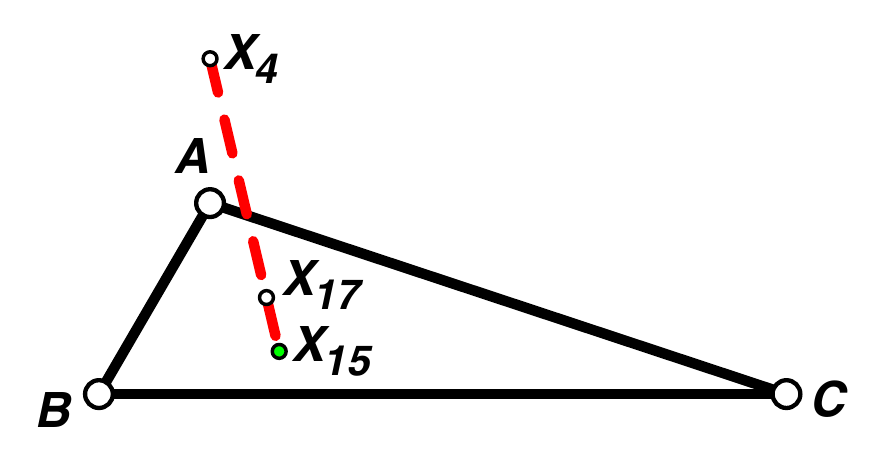}
{}
{$X_4$ -- $X_{17}$ -- $X_{15}$}

\pro{\cite{Dekov}}{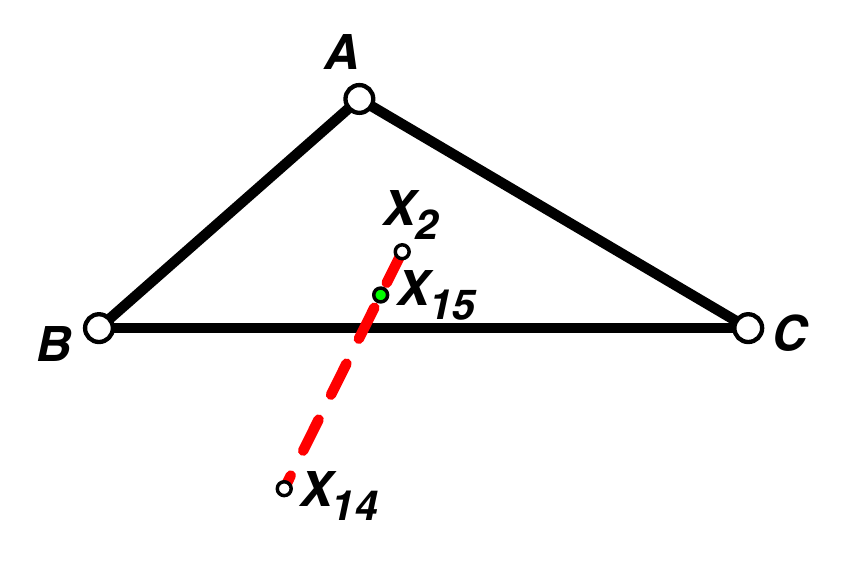}
{}
{$X_2$ -- $X_{15}$ -- $X_{14}$}

\pro{\cite{ETC15}}{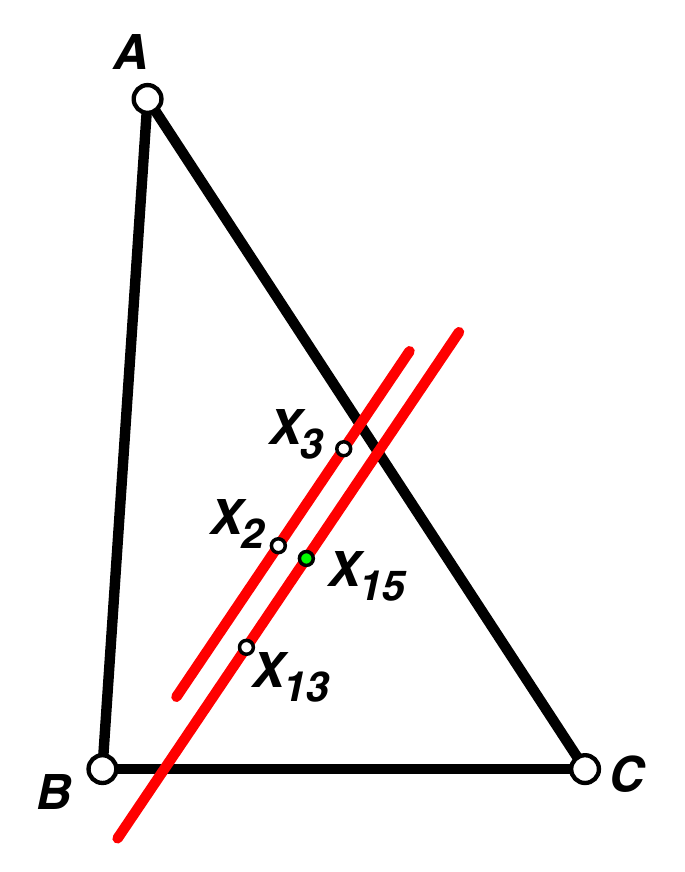}
{}
{$X_{13}X_{15}||X_2X_3$}
\extra{In other words, $X_{13}X_{15}$ is parallel to the Euler line of $\triangle ABC$.
}

\new
\pro{\cite{RG9103}}{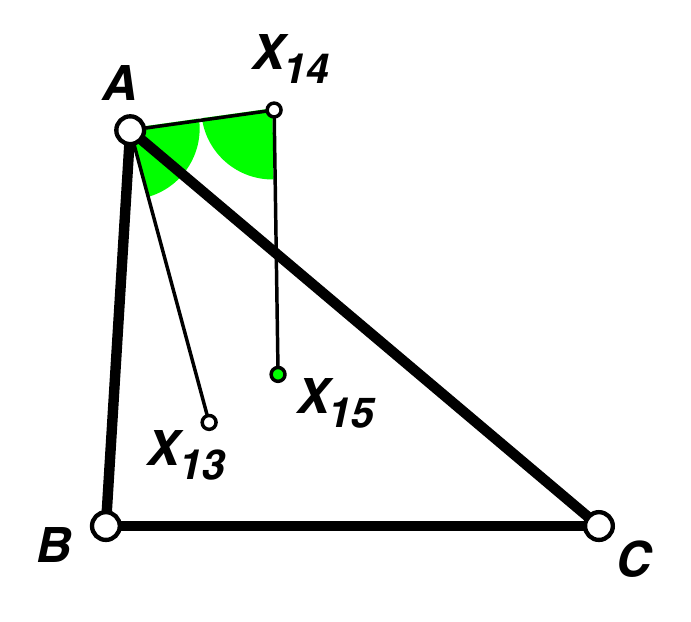}
{}
{$\angle X_{14}AX_{14}=\angle AX_{14}X_{15}$}

\new
\pro{}{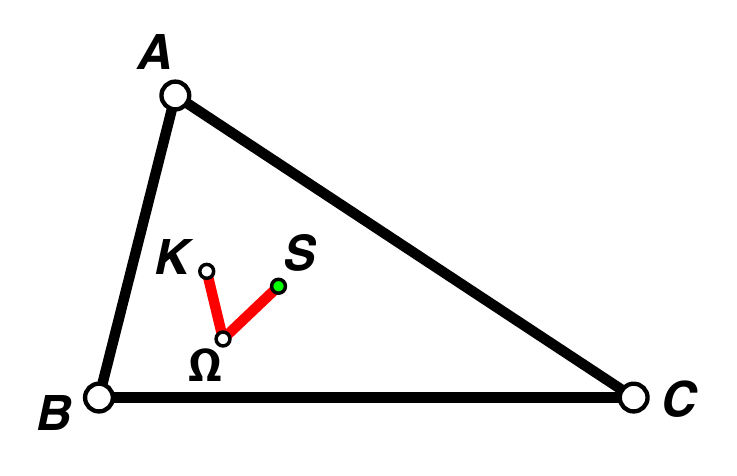}
{}
{$\angle S\Omega K=60\degrees$}

\new
\pro{}{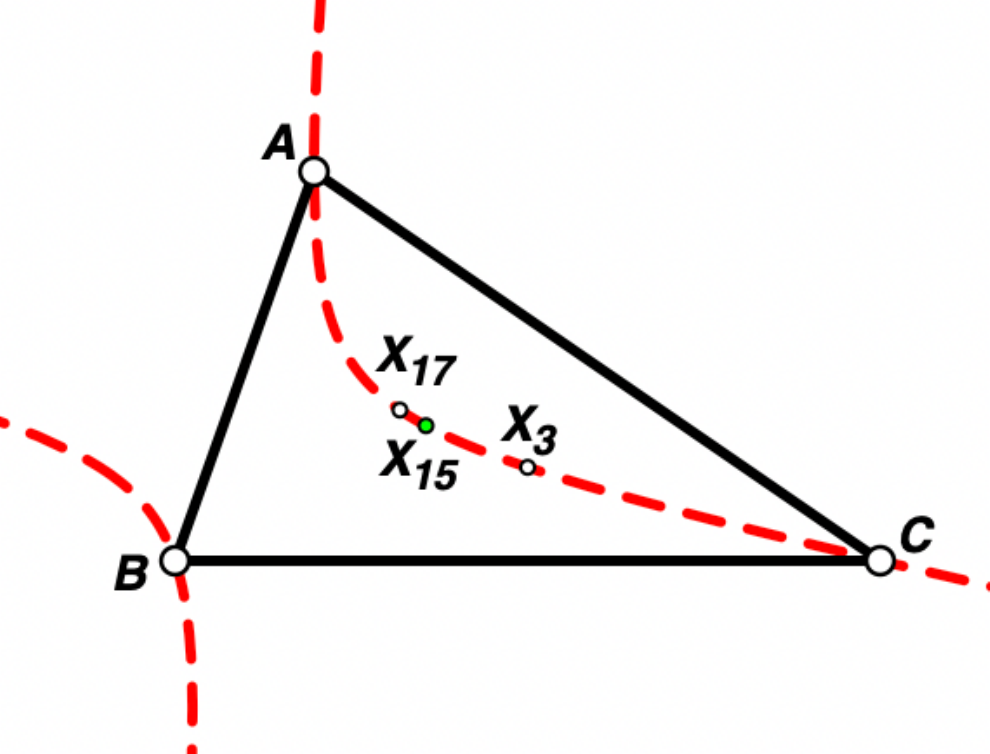}
{}
{$A$, $B$, $C$, $X_3$, $X_{15}$, $X_{17}$ lie on a conic.}

\void{
\new
\pro{}{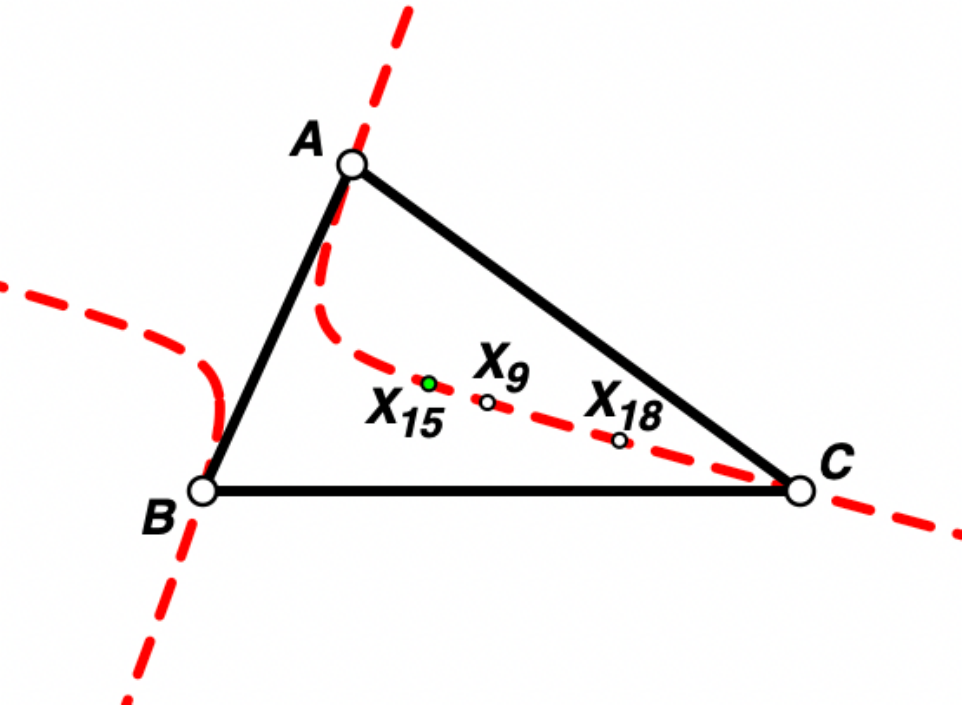}
{}
{$A$, $B$, $C$, $X_9$, $X_{15}$, $X_{18}$ lie on a conic.}
}

\mysubsection{$S$ and three notable points}

\new
\pro{}{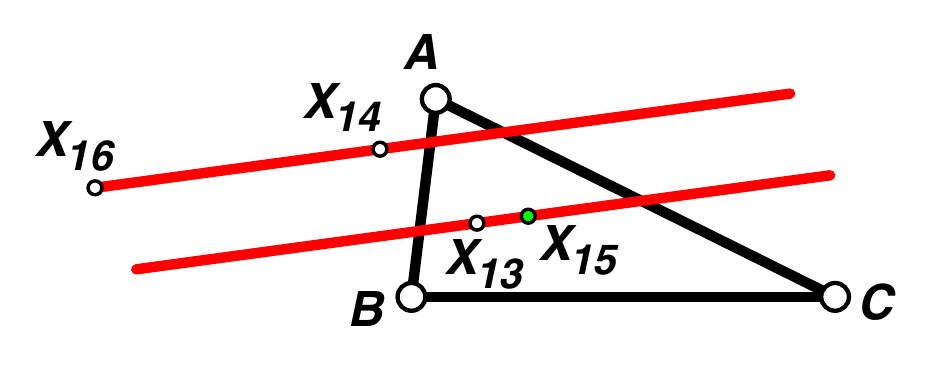}
{}
{$X_{13}X_{15}||X_{14}X_{16}$}

\pro[Brocard Axis]{\cite{MathWorld}}{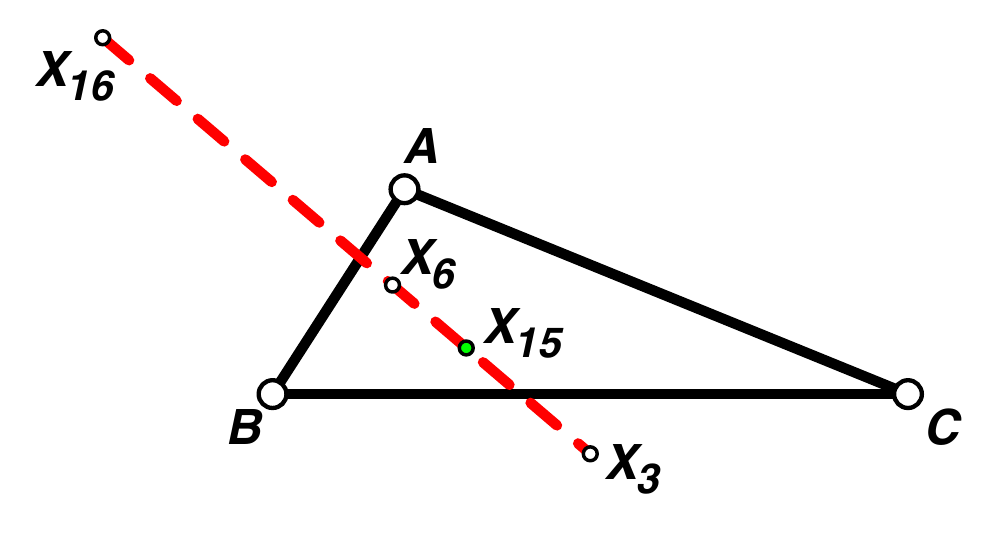}
{}
{$X_3$ -- $X_{15}$ -- $X_6$ -- $X_{16}$}
\extra{A few other named points that lie on the Brocard Axis are the 3rd power point, the Brocard midpoint, 
the Kenmotu point, and the Taylor center.
See \cite{MathWorld-B}.
}

\pro[Brocard Axis Metrics]{\cite{MathWorld}}{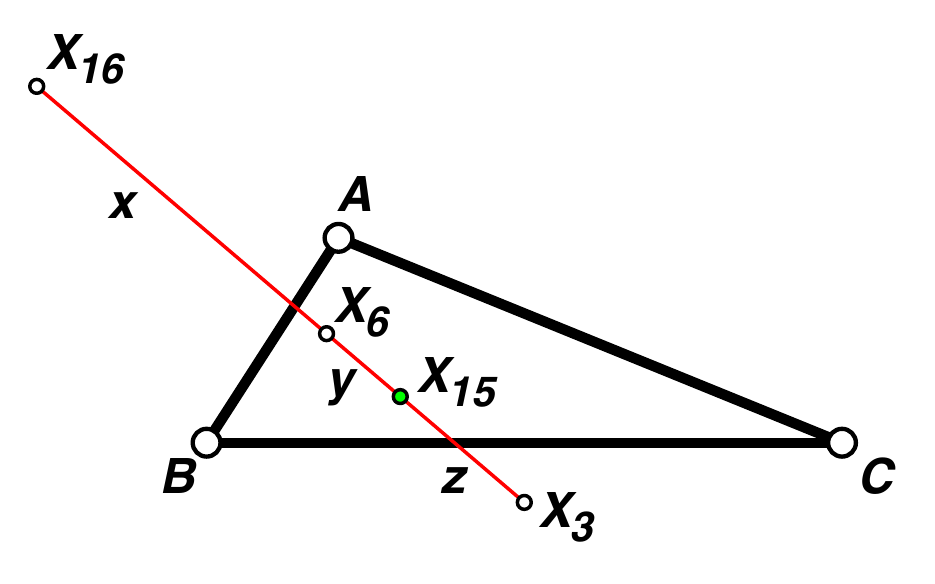}
{}
{$y(x+y+z)=xz$
\\
\then $z(x+y+z)=R^2$
\\
\then $\displaystyle \frac{1}{z}+\frac{1}{x+y+z}=\frac{2}{y+z}$
\\
\then $\displaystyle \frac{1}{x}+\frac{1}{x+y+z}=\frac{2}{x+y}$
\\
\then $(x+y)(y+z)=2xz$
}
\extra{See also \cite[p.~103]{Gallatly} and \cite{Dekov}.
}

\pro{\cite{MathWorld}}{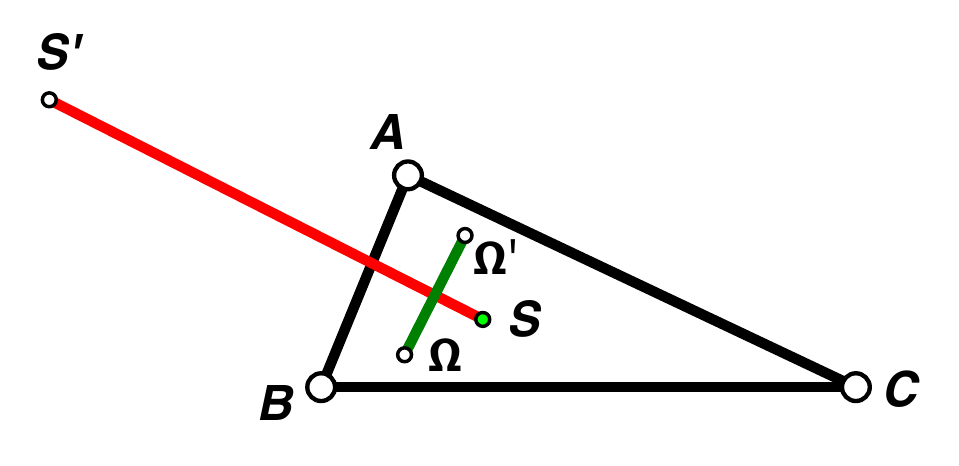}
{}
{$SS'$ is the perpendicular bisector of $\Omega\Omega'$.}
\extra{The point of intersection of the two lines is $X_{39}$, the Brocard midpoint.
}

\new
\pro{}{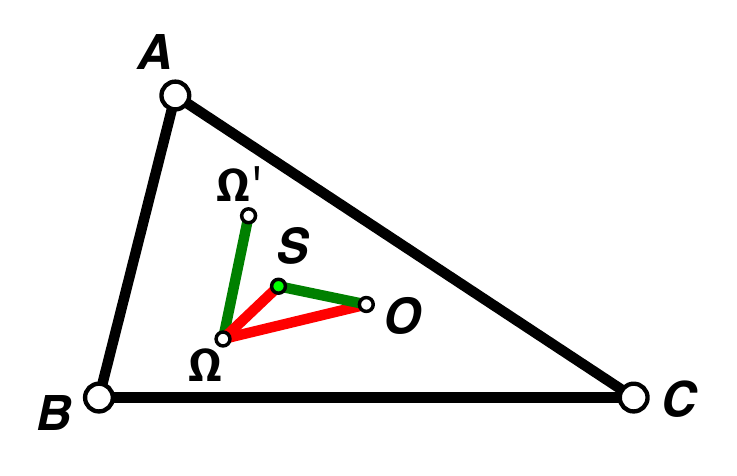}
{}
{$O\Omega\cdot S\Omega=OS\cdot \Omega\Omega'$}

\pro{\cite{ADGEOM5310}}{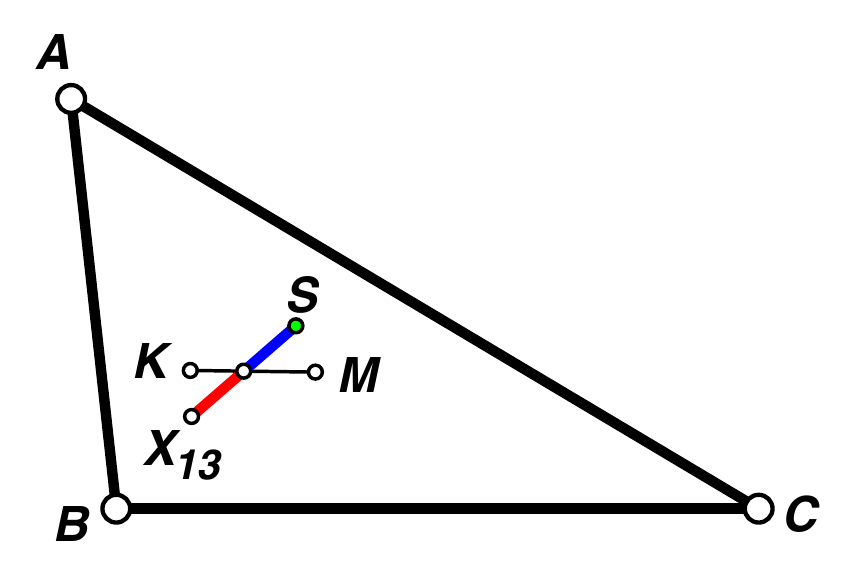}
{}
{$KM$ bisects $SX_{13}$}

\pro{\cite{ADGEOM5285}}{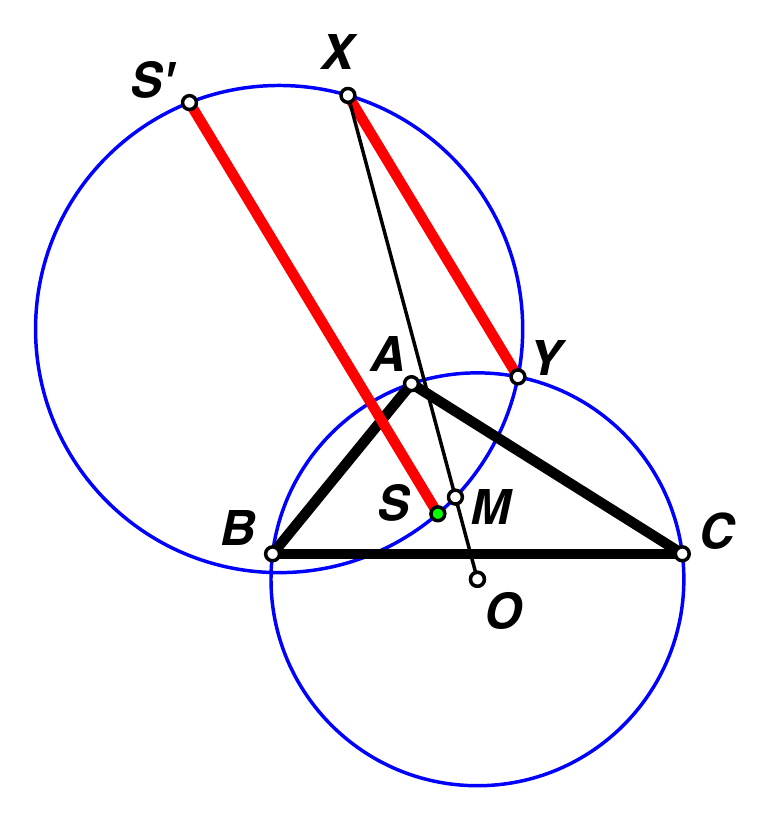}
{}
{$SS'\parallel XY$}
\extra{$Y=X_{110}$.}

\new
\pro{\cite{RG9091}}{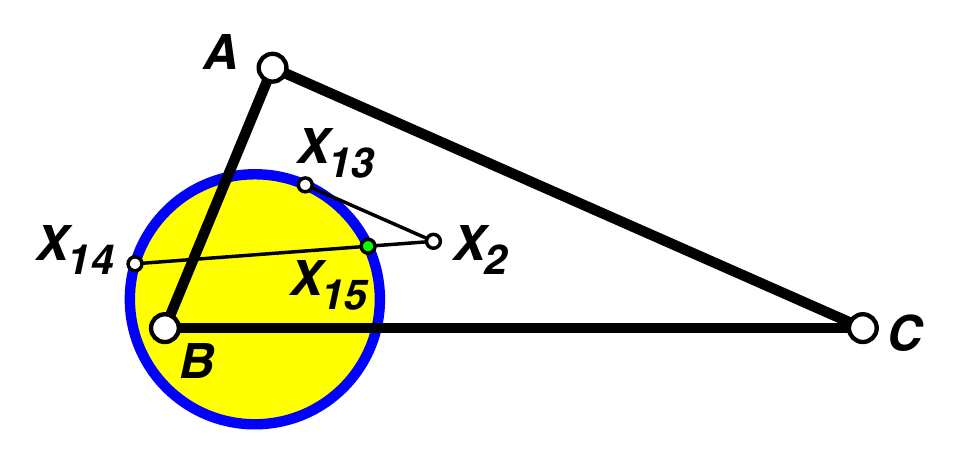}
{}
{$X_2X_{13}$ is tangent to $\odot X_{13}X_{14}X_{15}$}

\pro{\cite{RG8673}}{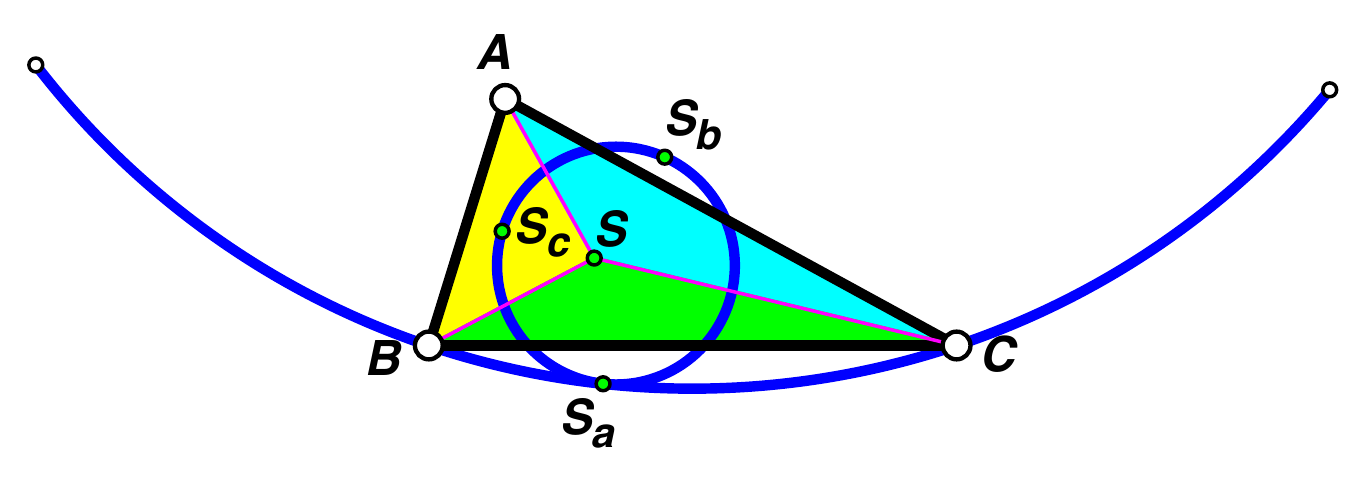}
{}
{$\odot S_aS_bS_c$ is tangent to $\odot SBC$.}
\extra{By symmetry, $\odot S_aS_bS_c$ is also tangent to $\odot SCA$ and $\odot SAB$.
The center of $\odot S_aS_bS_c$ is $X_{5238}$.}

\mysubsection{$S$ and 4 or more notable points}

\nopagebreak
\pro{\cite{Altintas1625}}{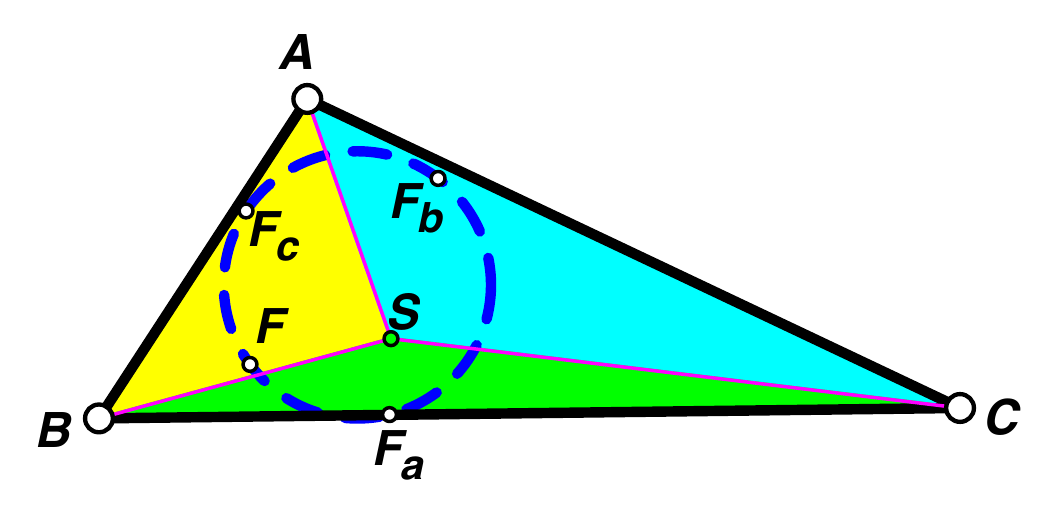}
{}
{$F_a$, $F_b$, $F_c$, $F$ concyclic.}
\extra{The result also works for $S'$. See \cite{RG3673}.}

%\newpage
\new
\prowide{}{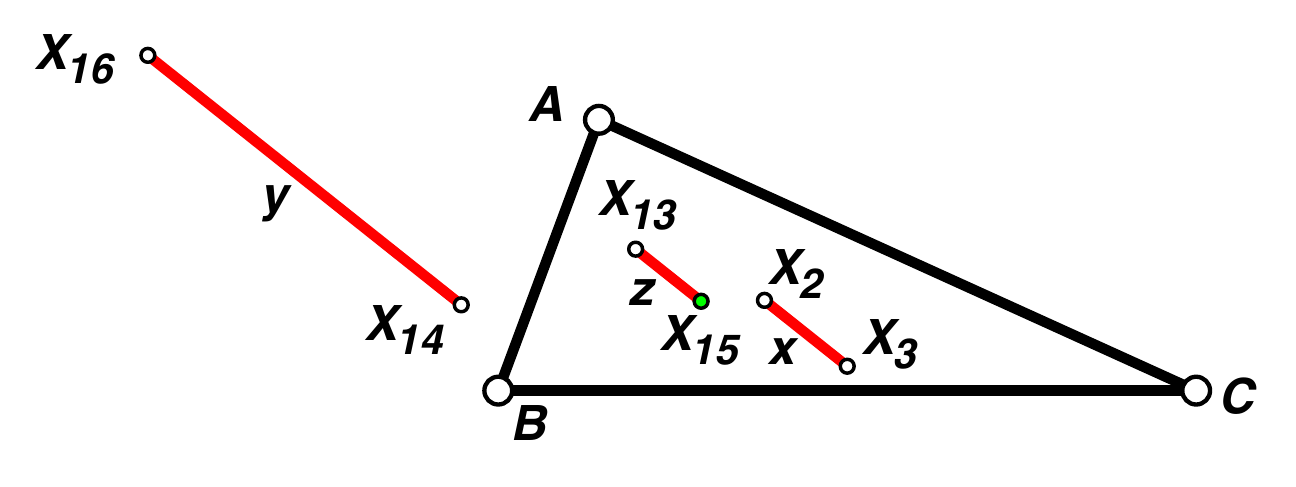}
{}
{$\displaystyle\frac1x+\frac1y=\frac1z$}

\prowide[Evans Conic]{\cite{Evans}}{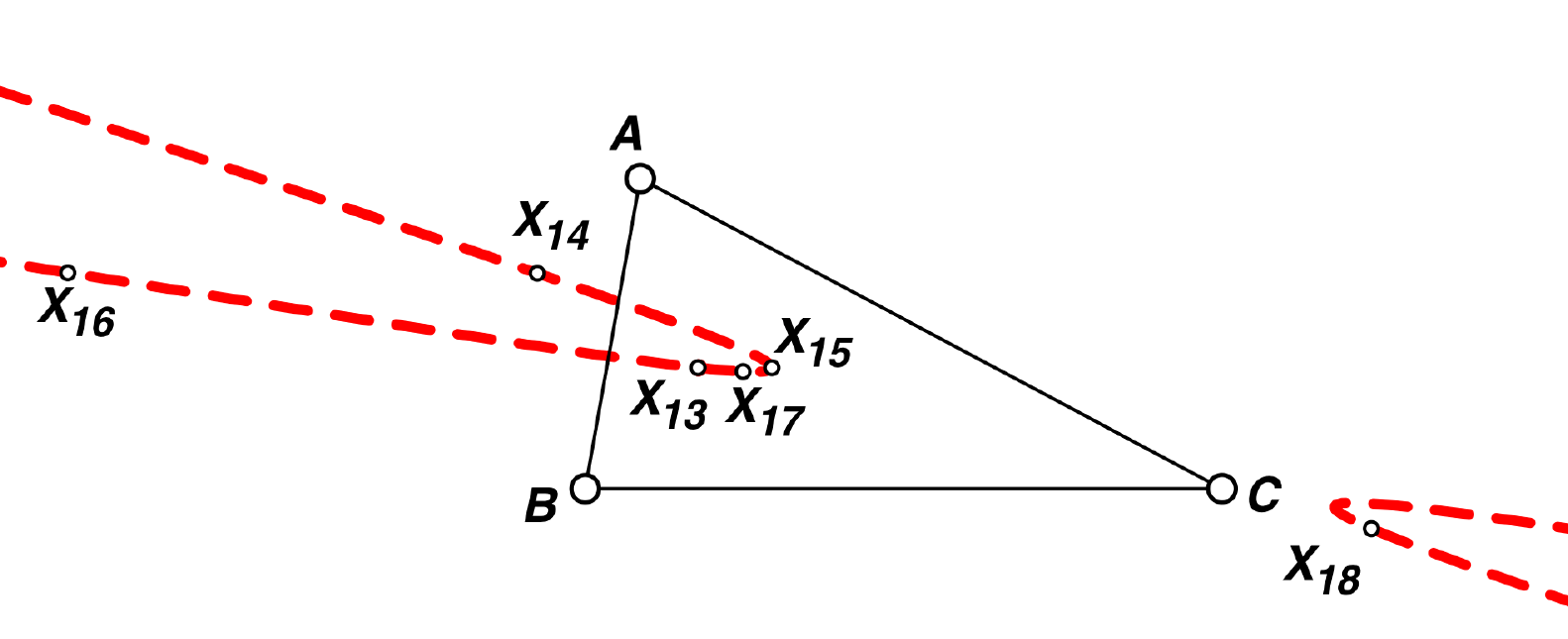}
{}
{$X_{13}$, $X_{14}$, $X_{15}$, $X_{16}$, $X_{17}$, $X_{18}$ lie on a conic.}
\extra{See also \cite{MathWorld-Evans}.
}

\new
\prowide{\cite{Euclid3337}}{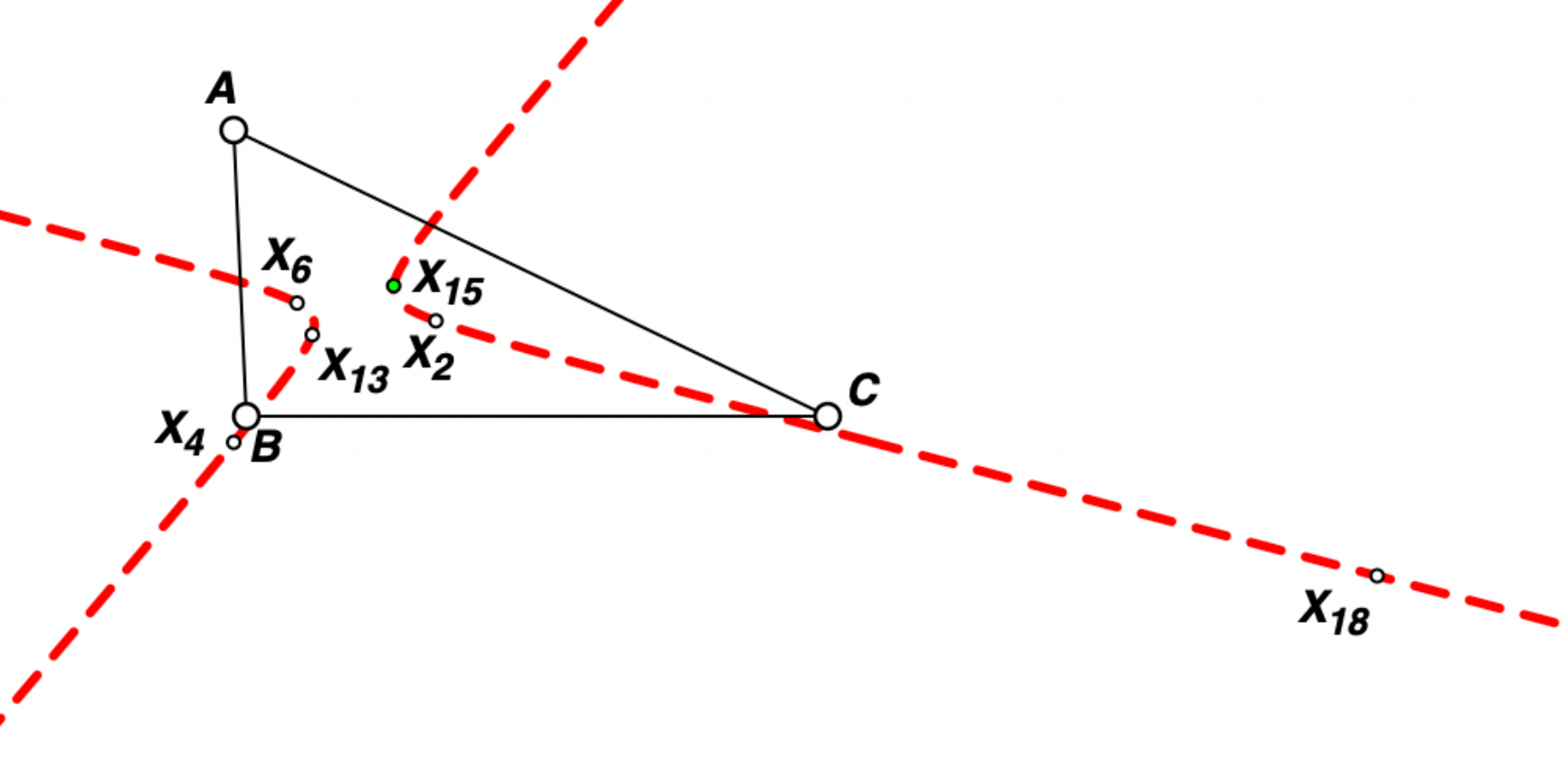}
{}
{$X_{2}$, $X_{4}$, $X_{6}$, $X_{13}$, $X_{15}$, $X_{18}$ lie on a conic.}

\pro[Neuberg Cubic]{\cite{MathWorld-Neuberg}}{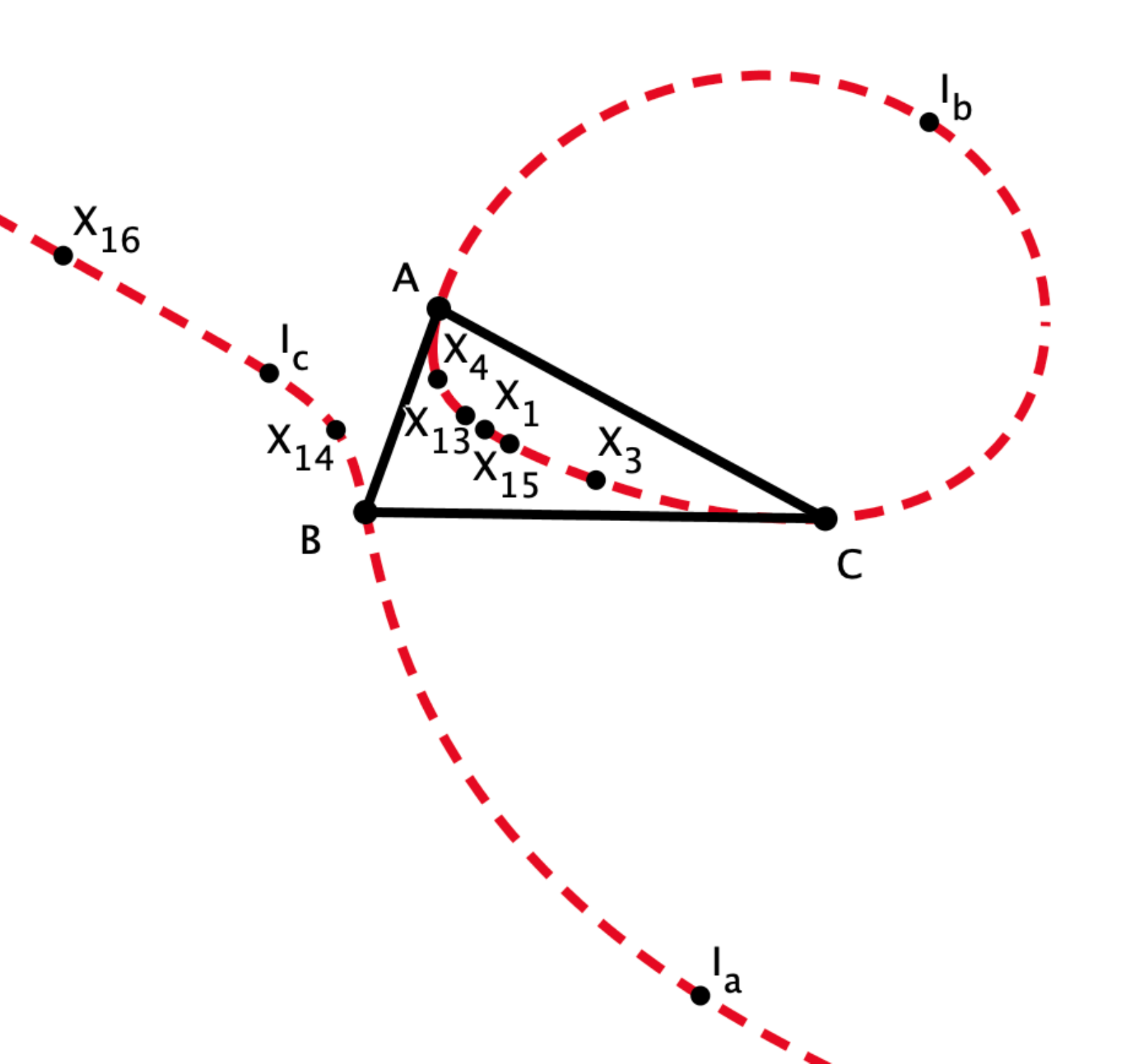}
{}
{The following ten points lie on a cubic curve:
$X_1$, $X_3$, $X_4$, $X_{13}$, $X_{14}$, $X_{15}$, $X_{16}$, $I_a$, $I_b$, $I_c$,
where $I_a$ is the $A$-excenter of $\triangle ABC$.}
\extra{A few other points that lie on the Neuberg cubic are the reflections of the
vertices of the triangle about their opposite sides and the six vertices
of equilateral triangles erected on the sides of $\triangle ABC$.
See \cite{Gibert-Neuberg}.
}

\newpage
\vspace*{-22pt}
\mysection{Triangle plus lines through $S$}

\mysubsection{antiparallels}

\pro{\cite{Altintas1746}}{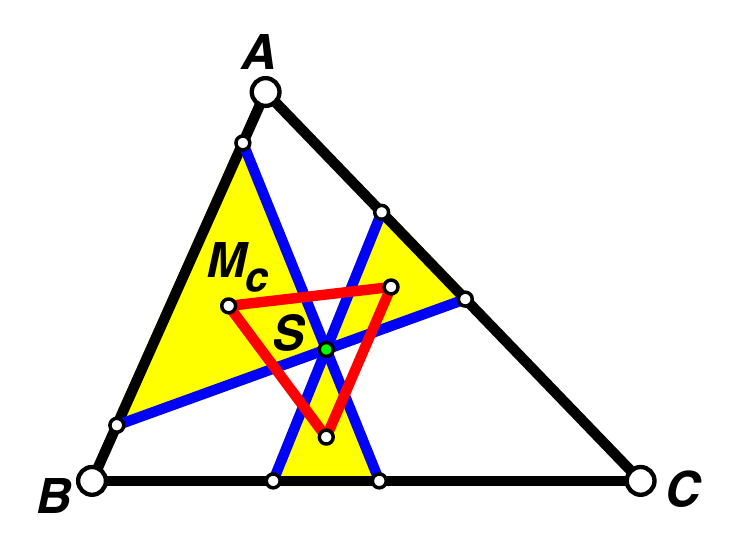}
{}
{$\triangle M_aM_bM_c$ is equilateral.}
\extra{See also \cite{RG8828}.}

\mysubsection{apothems}

\nopagebreak
\pro[Trilinear Coordinates]{\cite{Wiki}}{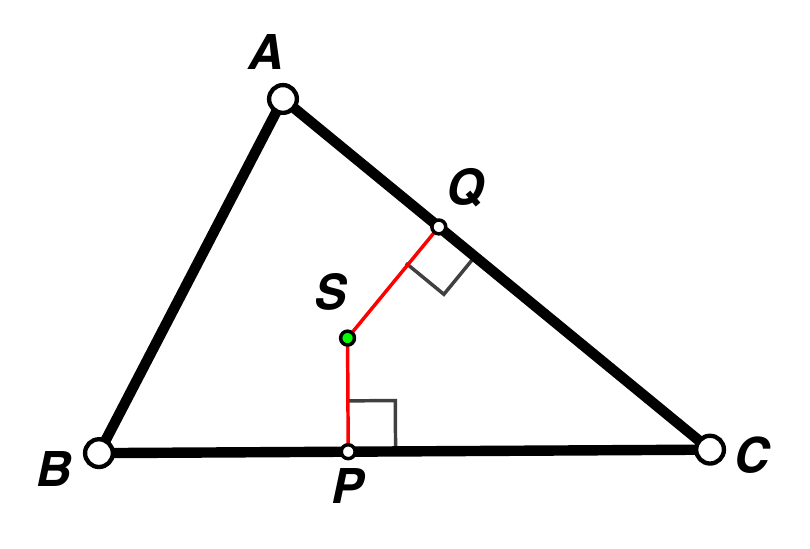}
{\void{Drop perpendiculars from $S$ to two sides of the triangle.}}
{$\displaystyle\frac{SP}{SQ}=\frac{\sin(A+60\degrees)}{\sin(B+60\degrees)}$.}

\pro[Pedal Triangle]{\cite[p.~106]{Gallatly}}{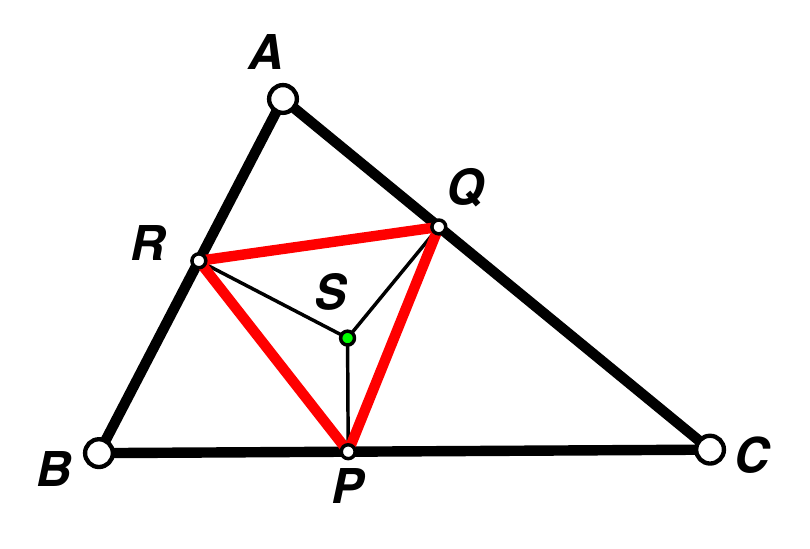}
{\void{Drop perpendiculars from $S$ to the sides of the triangle.}}
{$\triangle PQR$ is an equilateral triangle.\\
\then $\displaystyle [PQR]=\frac{2K^2\sqrt3}{a^2+b^2+c^2+4K\sqrt3}$}
\extra{The center of the equilateral triangle is the midpoint of $SX_{13}$.
See \cite{RG3786} for the area formula.}

\mysubsection{cevians}

\nopagebreak

\new
\pro[Cevian Length]{}{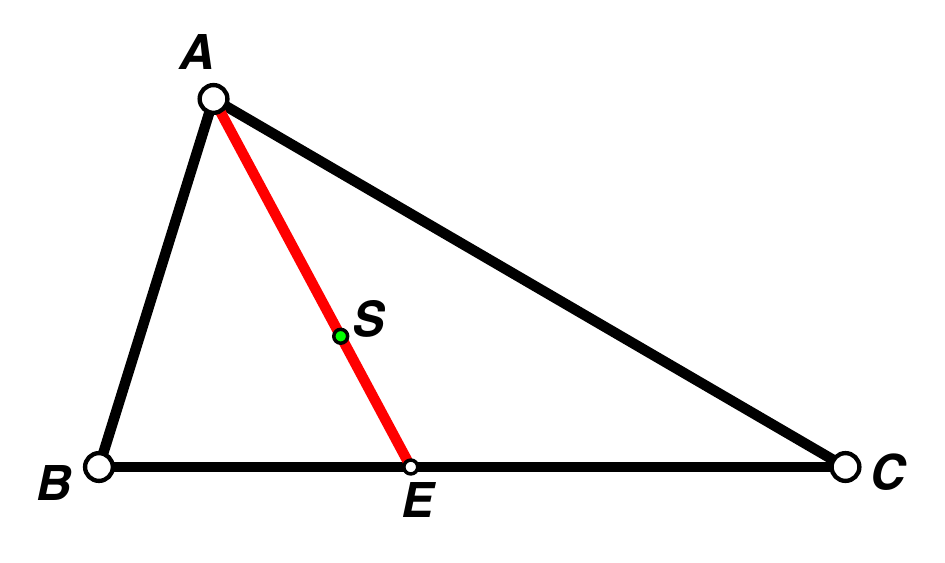}
{}
{$\displaystyle AE=\frac{4\sqrt2bcK\sqrt{a^2+b^2+c^2+4K\sqrt3}}{(b^2+c^2)(4K+a^2\sqrt3)-\sqrt3(b^2-c^2)^2}$
}

%[Orthiac Triangle]
\pro{\cite{Altintas1684}}{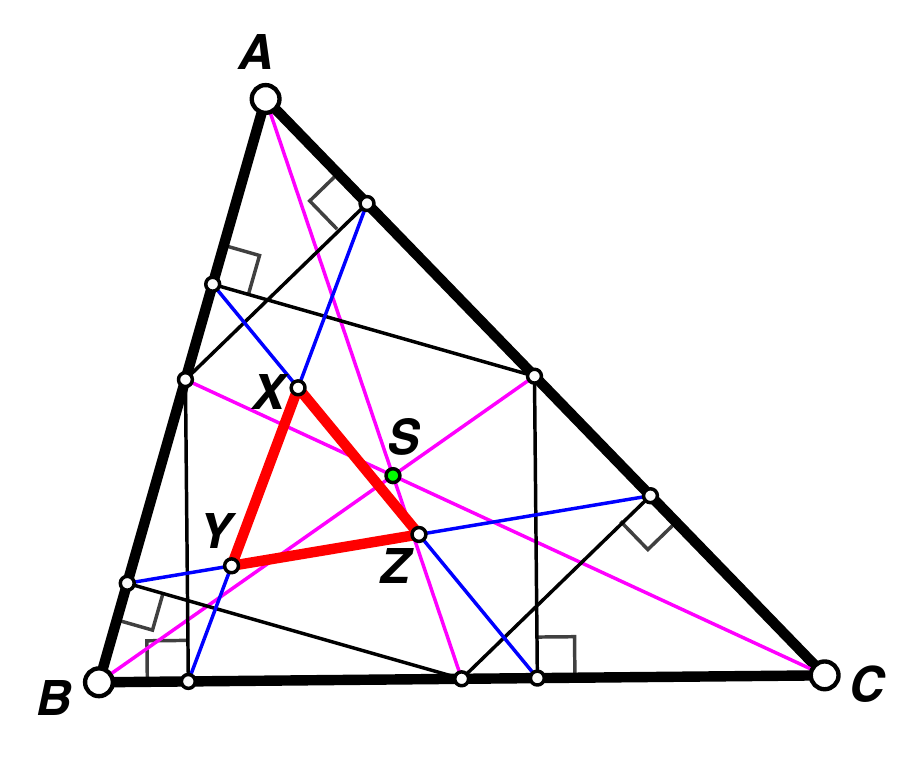}
{}
{$\triangle XYZ$ is equilateral.}
\extra{See also \cite{RG8664}.}

\mysubsection{circlecevians}

\nopagebreak
\new
\pro{\cite{RG4848}}{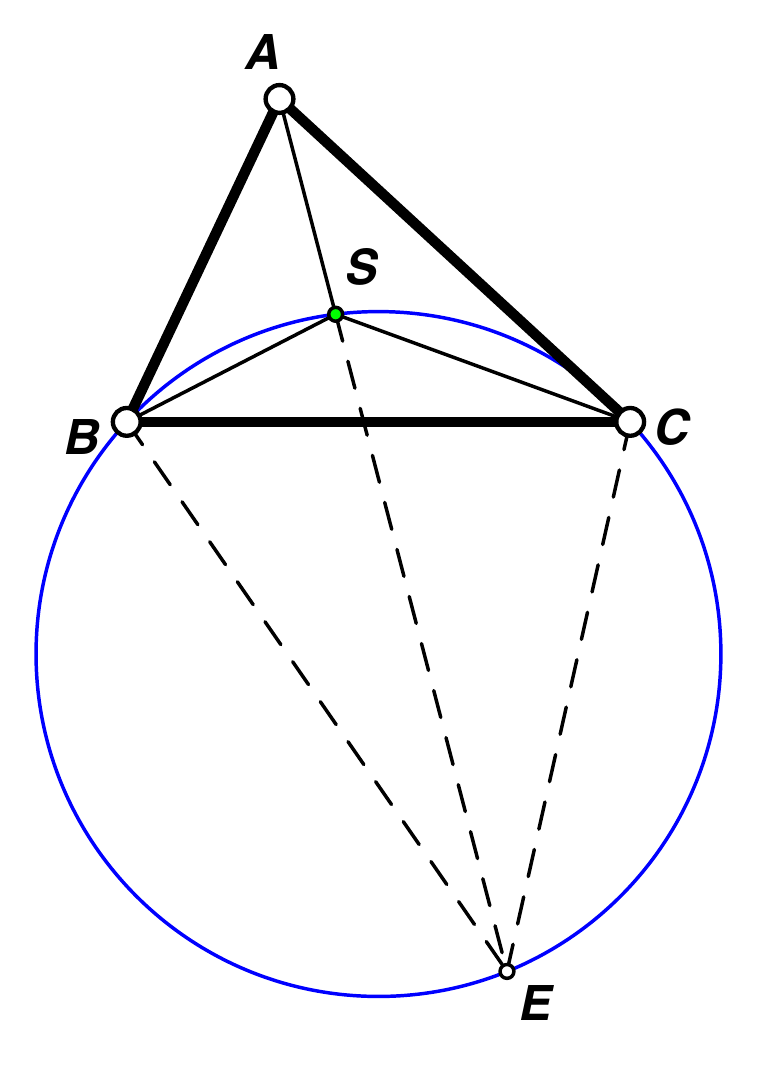}
{\void{$AS$ meets $\odot BCS$ at $E$.}}
{$\angle EBA=\angle ACE=120\degrees$.}

\pro{\cite{Euclid1455}}{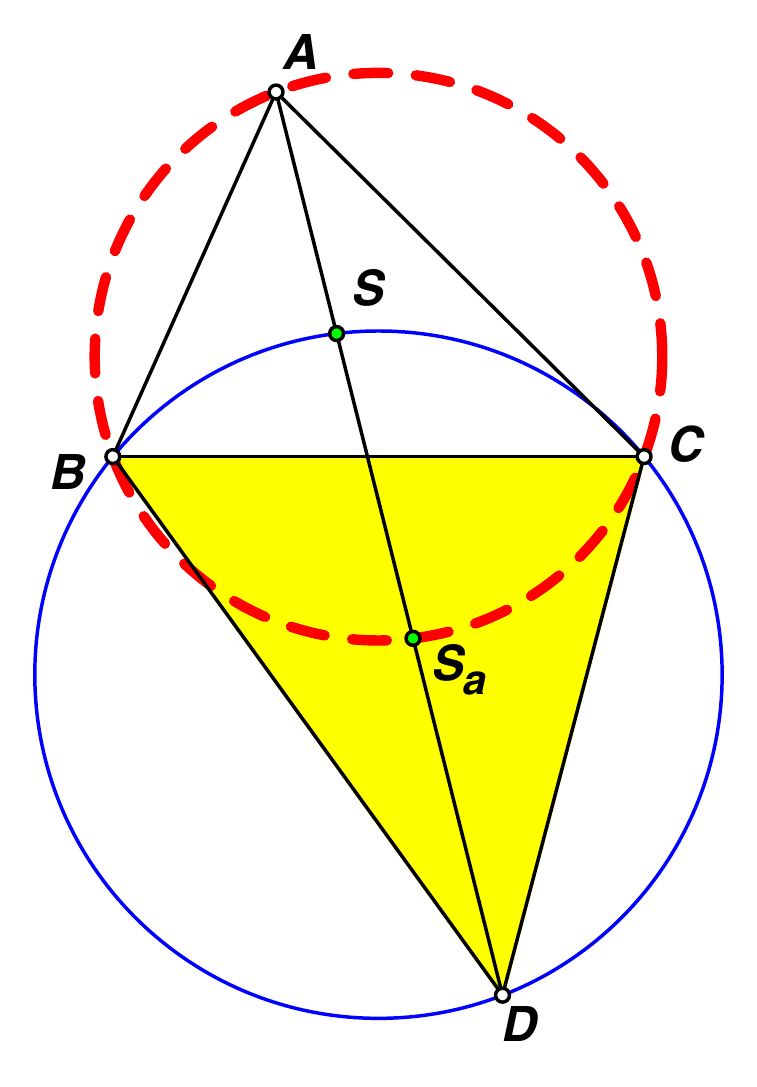}
{}
{$A$ -- $S$ -- $S_a$}

\mysubsection{circumcevians}

\nopagebreak
\new
\pro{\cite{RG4849}}{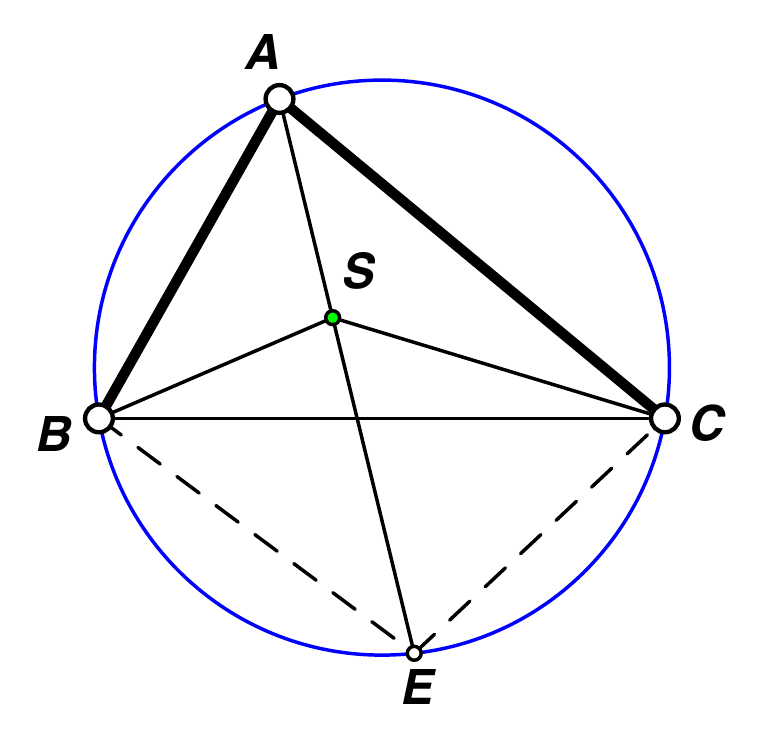}
{\void{$AS$ meets $\odot ABC$ at $E$.}}
{$\angle EBS=\angle SCE=60\degrees$.}

\pro{\cite[p.~296]{Johnson}}{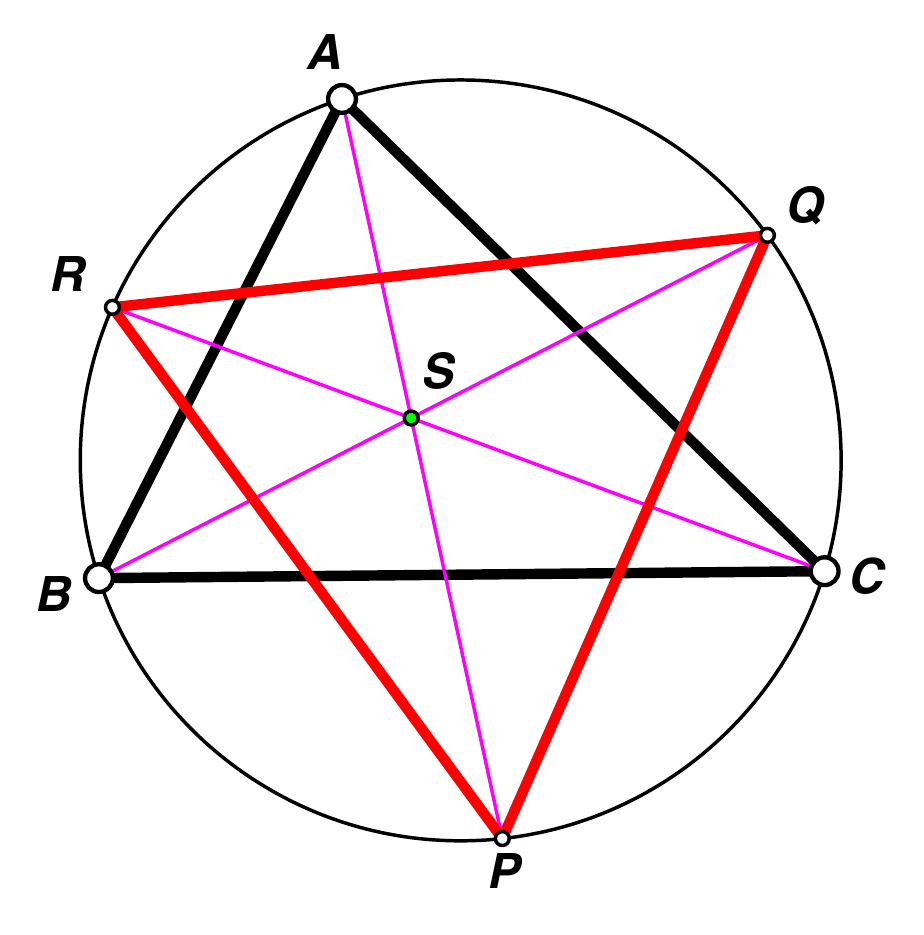}
{\void{Let $AP$, $BQ$, and $CR$ be the circumcevians through $S$.}}
{$\triangle PQR$ is an equilateral triangle.}

\pro{\cite{RG8754}}{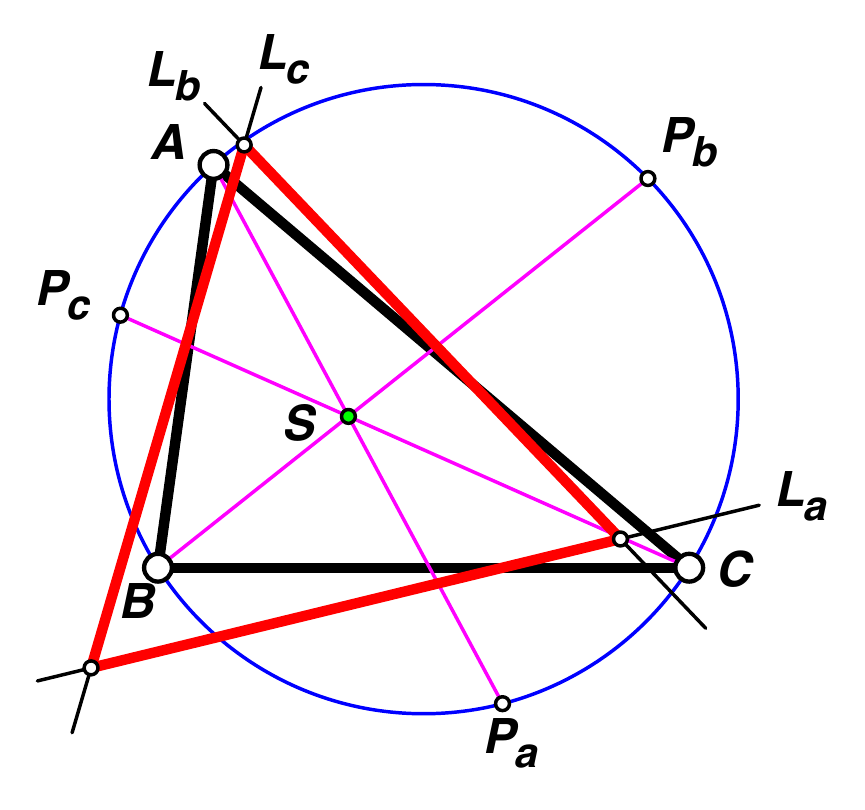}
{$L_a$ is the Simson line of $P_a$.}
{$L_a$, $L_b$, $L_c$ bound an equilateral triangle.}
\extra{The center of the triangle is the nine-point center of $\triangle ABC$.}

\pro{\cite[Thm. 8.1.1]{Cosmology}}{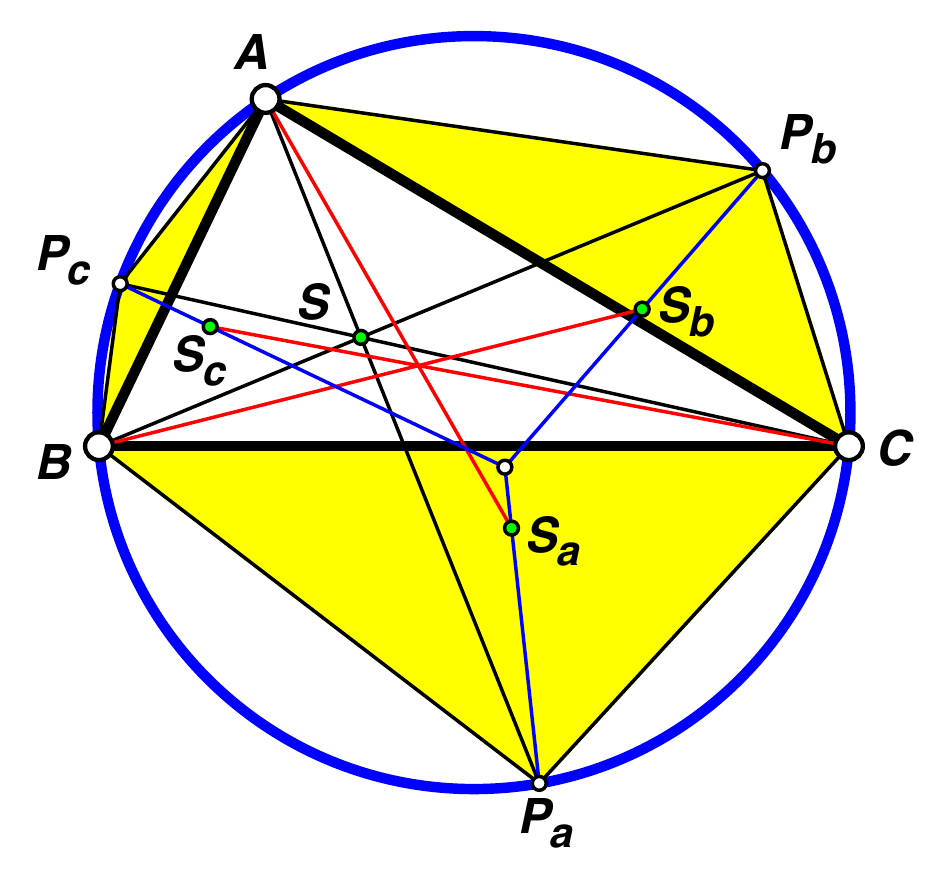}
{}
{$AS_a$, $BS_b$, $CS_c$ concur.\\
\then $P_aS_a$, $P_bS_b$, $P_cS_c$ concur.}

\mysubsection{inclines}

\nopagebreak
\new
\pro{}{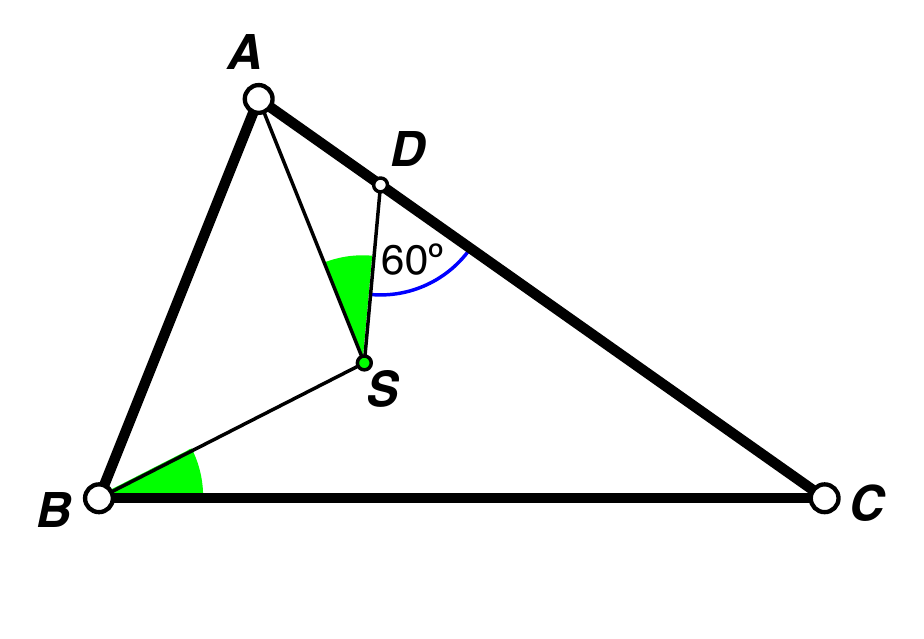}
{}
{$\angle DSA=\angle CBS$.}
\extra{This follows from Property \ref{property-iso-6}.}

%\pagebreak
%\mysubsection{three inclines through $S$}

\new
\pro{}{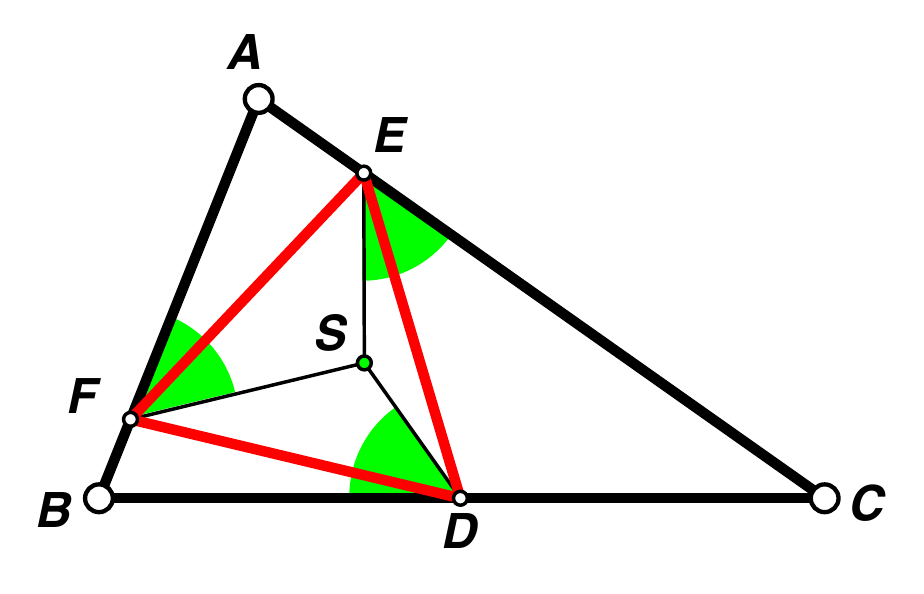}
{}
{$\triangle DEF$ is equilateral.}

\mysubsection{parachords}

\nopagebreak
\new
\pro{\cite{RG9075}}{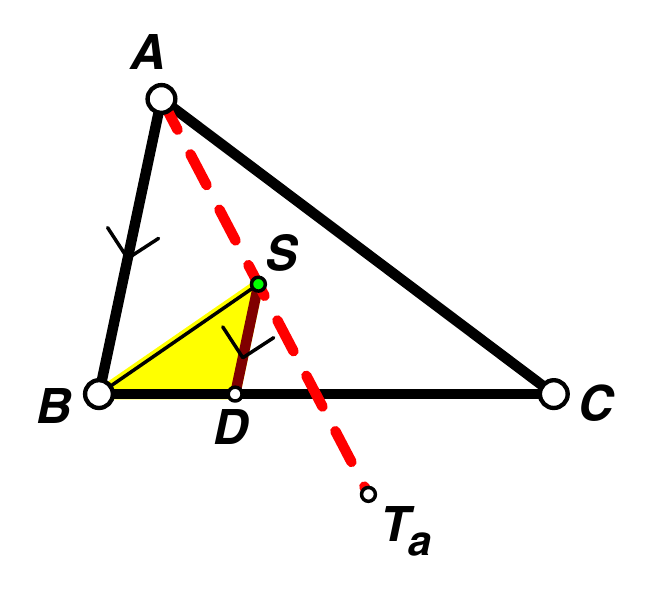}
{}
{$A$ -- $S$ -- $T_a$}

\new
\pro{\cite{PGR43}}{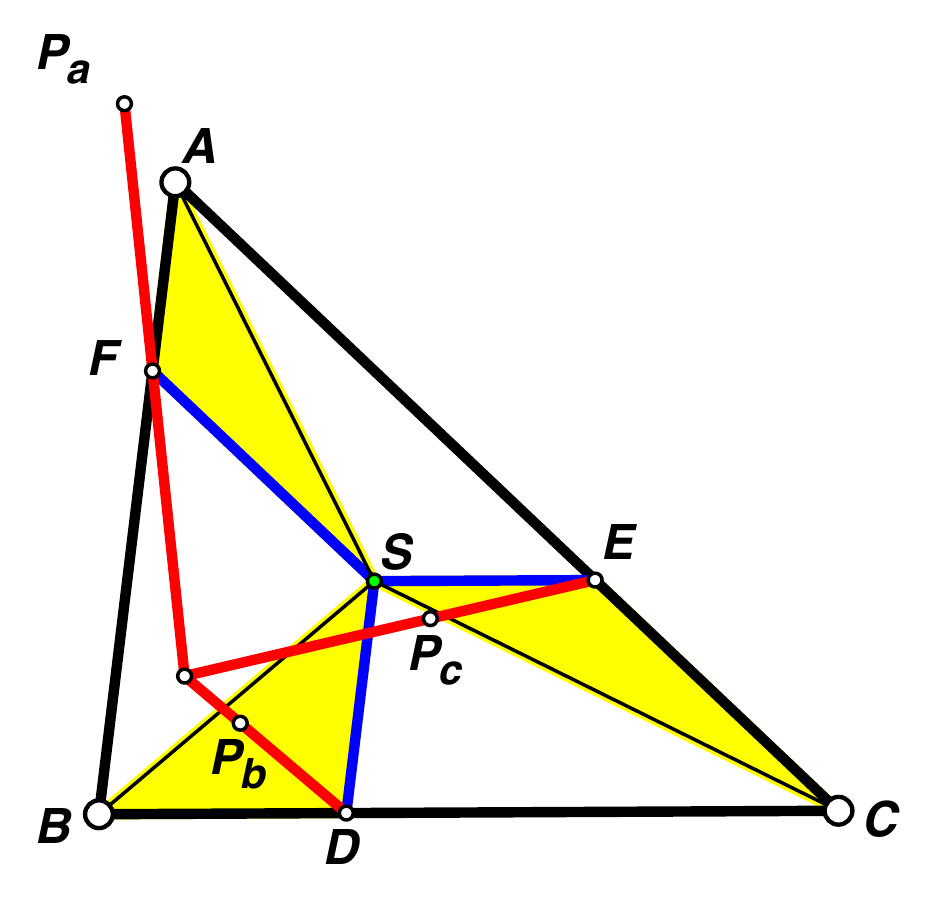}
{$P_a$ is the 2nd Napoleon point of $\triangle ASF$.}
{$DP_b$, $EP_c$, and $FP_a$ are concurrent.}

\new
\pro{\cite{RG8850}}{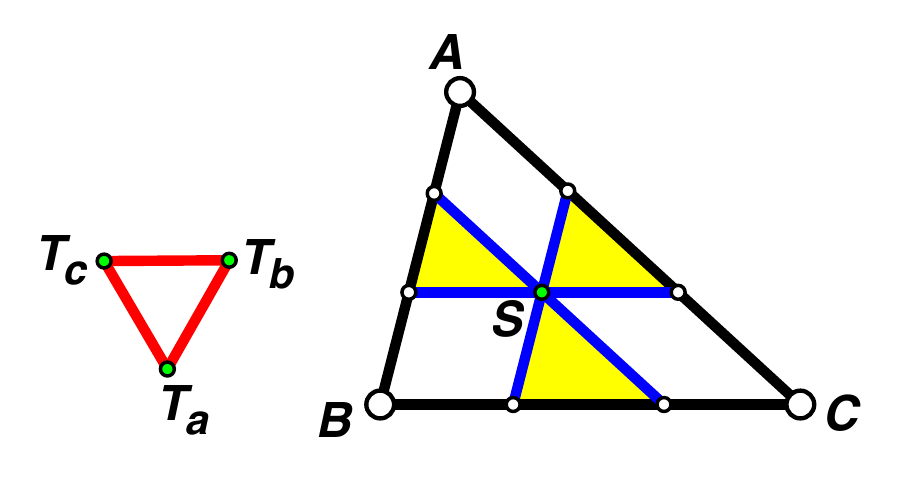}
{}
{$\triangle T_aT_bT_c$ is equilateral.}
\extra{The center of the equilateral triangle is $X_{39555}$.}

%\nopagebreak
\new
\pro{\cite{RG9064}}{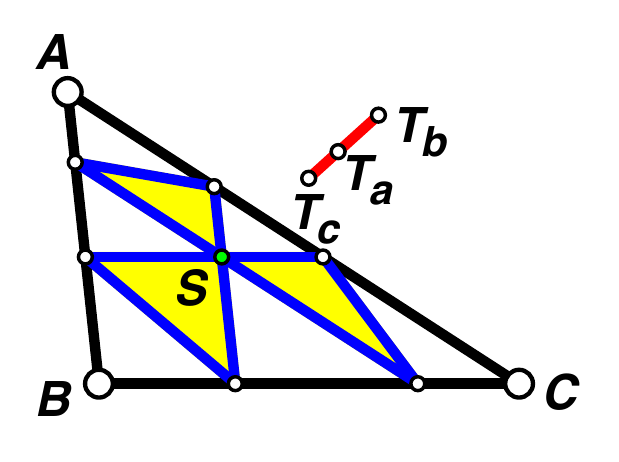}
{}
{$T_c$ -- $T_a$ -- $T_b$.}

\nopagebreak
\new
\pro{\cite{RG9063}}{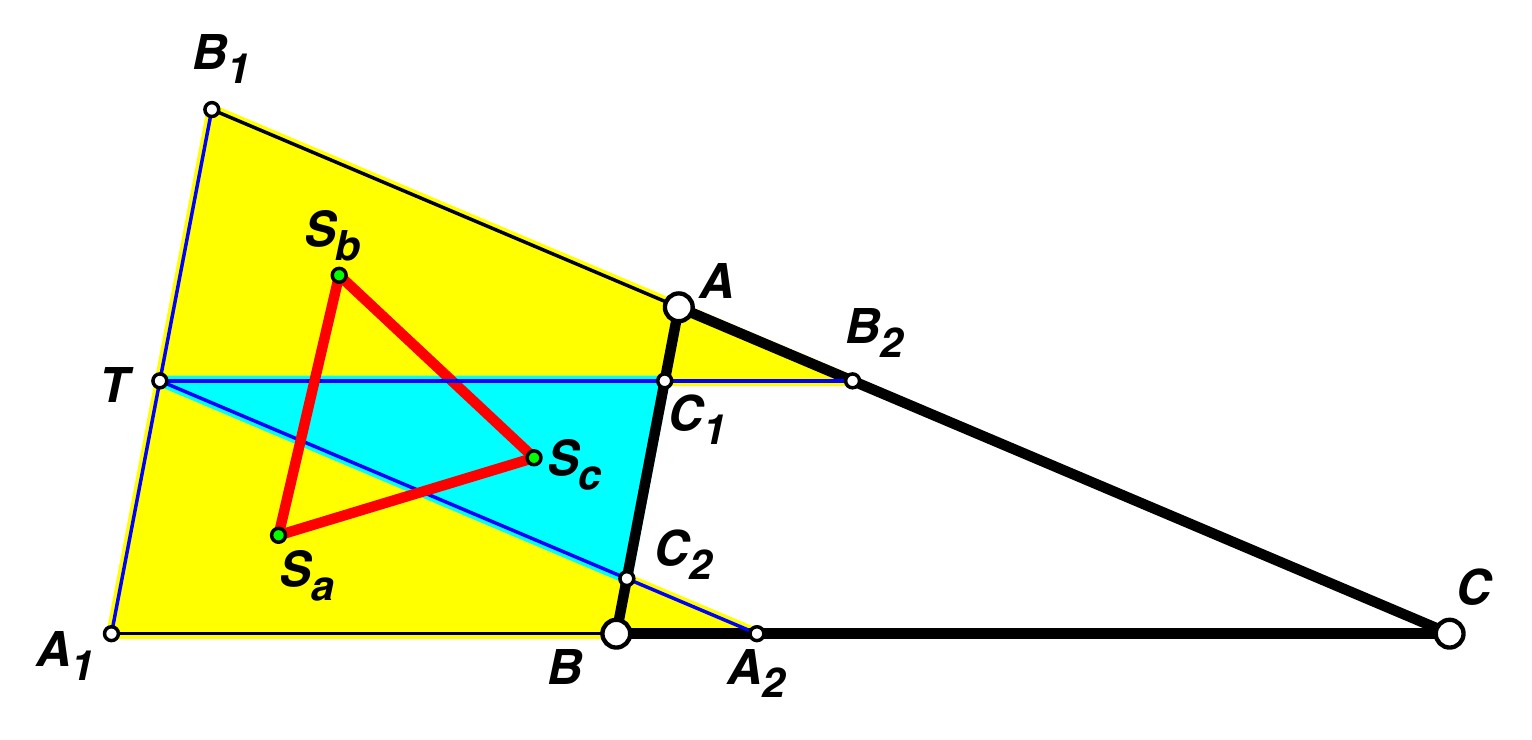}
{}
{$\triangle S_aS_bS_c$ is equilateral.}

\mysubsection{spokes}

\nopagebreak
\new
\pro{\cite{PGR41}}{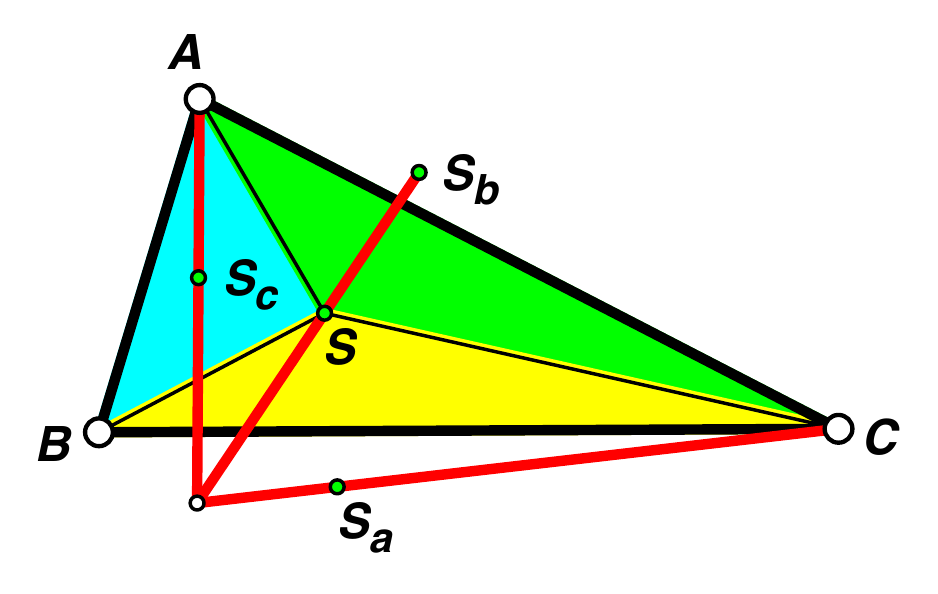}
{}
{$AS_c$, $CS_a$, and $SS_b$ are concurrent.}

\new
\pro{\cite{PGR42}}{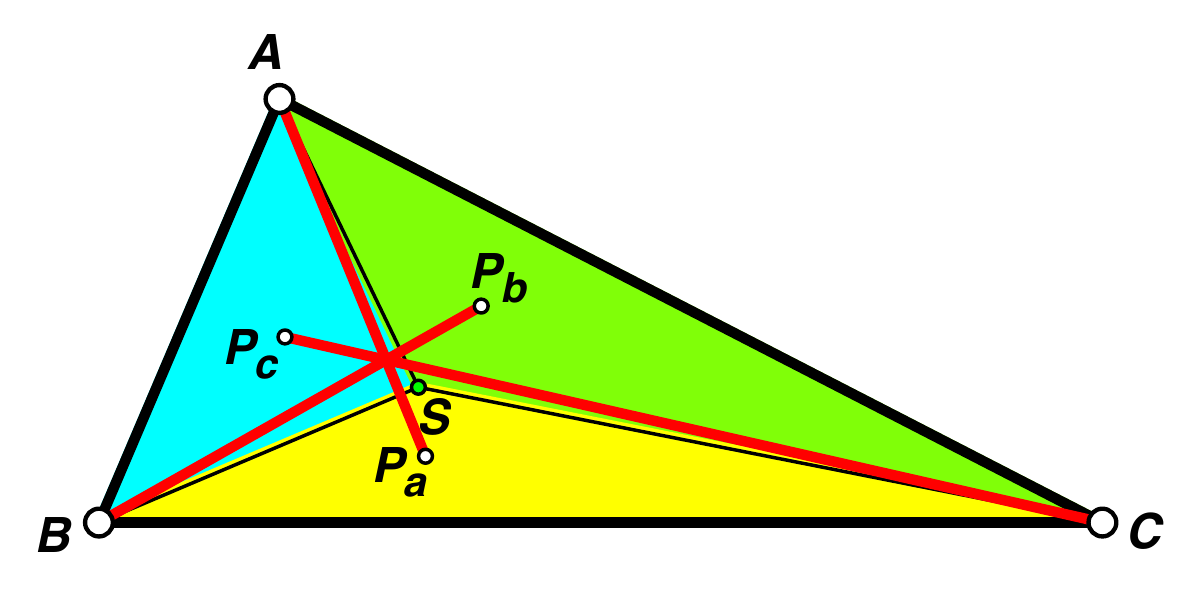}
{$P_a$ is the incenter of $\triangle SBC$.}
{$AP_a$, $BP_b$, and $CP_c$ are concurrent.}
\extra{The result remains true if $P_a$ is replaced by $X_n(SBC)$, for $n=$
2, 3, 6, 13, 15, 31, 32, 36, 39, or 50.}

\newpage
\vspace*{-22pt}
\mysection{Triangle with $S$ in subtriangles}

\mysubsection{formed by $H$}

\pro{\cite{RG8762}}{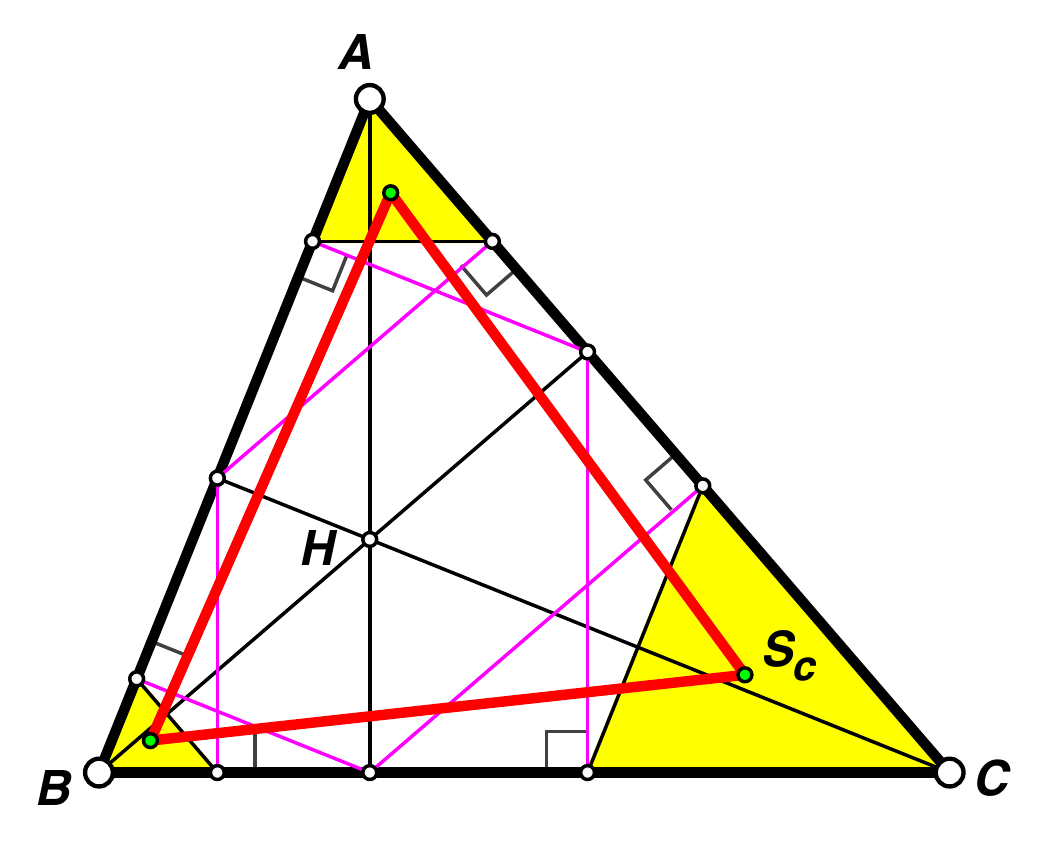}
{}
{$\triangle S_aS_bS_c$ is equilateral.}
\extra{$\triangle T_aT_bT_c$ is equilateral and $\triangle T_aT_bT_c\cong \triangle S_aS_bS_c$.}

\pro[3rd isodynamic-Dao triangle]{\cite{ETC31683}}{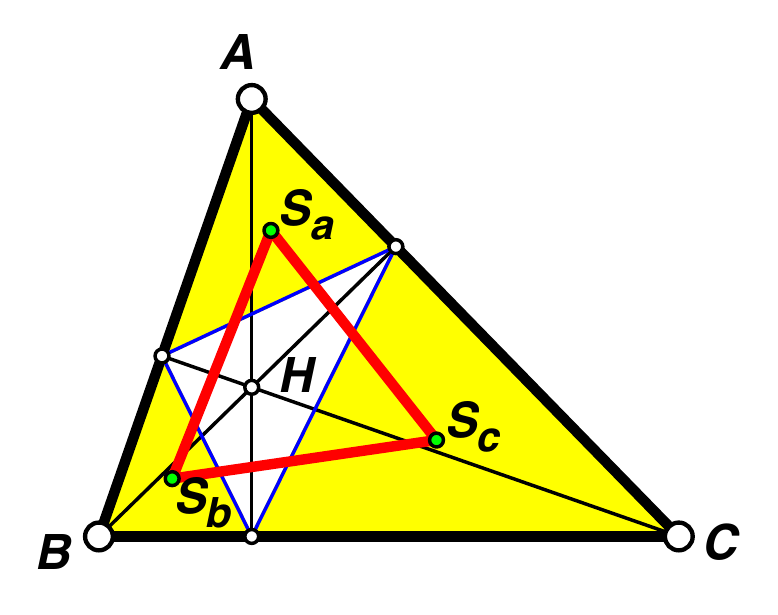}
{}
{$\triangle S_aS_bS_c$ is equilateral.}
\extra{Result is also true when using $S'$ instead of $S$.\\
{\Large $\bullet$}\ \ See also \cite{Altintas1}.}

\pro{\cite{PGR14-Nov-2021}}{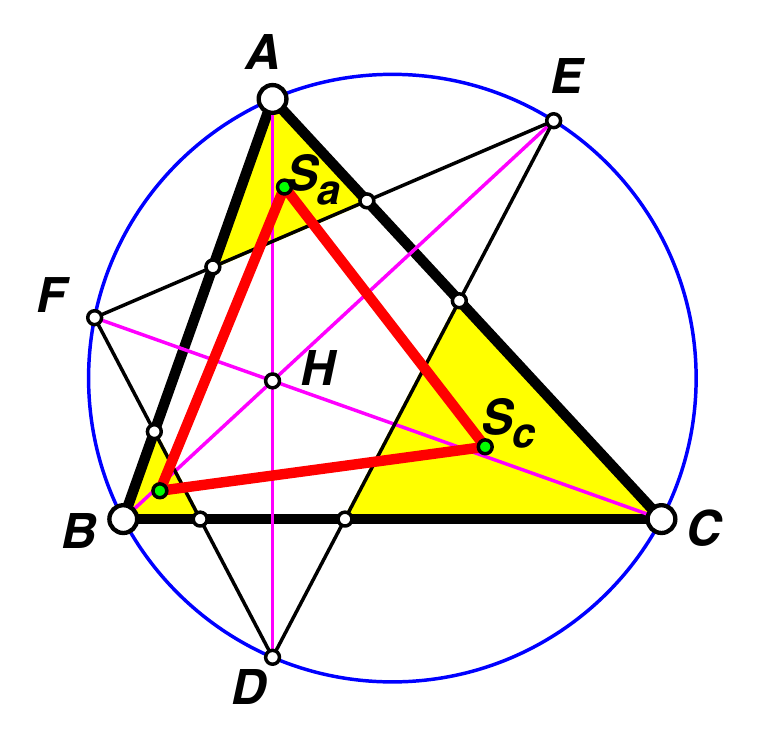}
{}
{$\triangle S_aS_bS_c$ is equilateral.}

\mysubsection{formed by $I$}

\nopagebreak
\pro{\cite{PGR14-Nov-2021}}{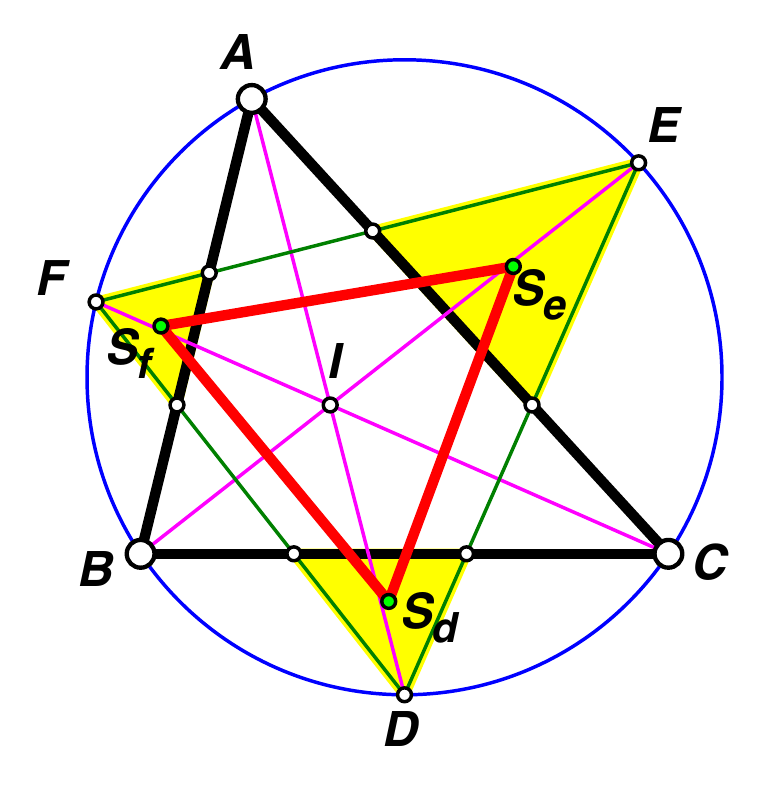}
{}
{$\triangle S_dS_eS_f$ is equilateral.}

\new
\pro{\cite{RG6825}}{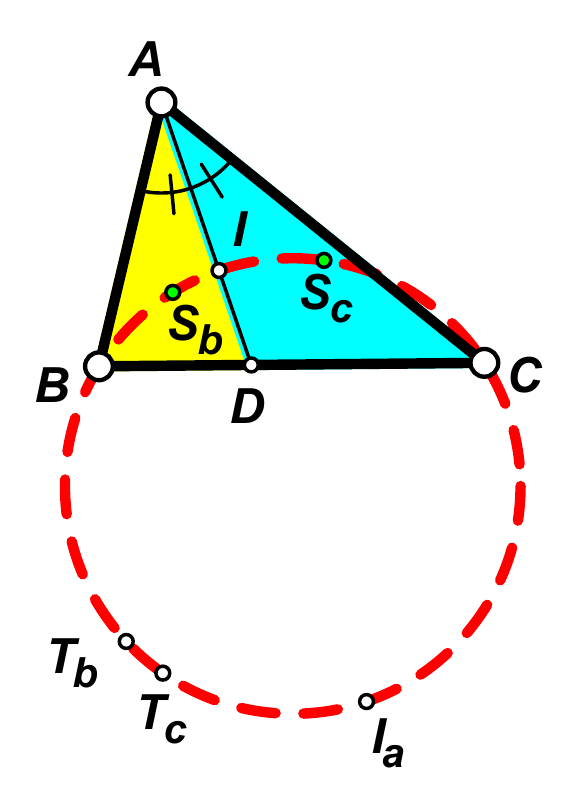}
{}
{$B$, $C$, $S_b$, $S_c$, $T_b$, $T_c$, $I$, $I_a$ concyclic}

\mysubsection{formed by $K$}

\nopagebreak
\pro{\cite{Altintas1}}{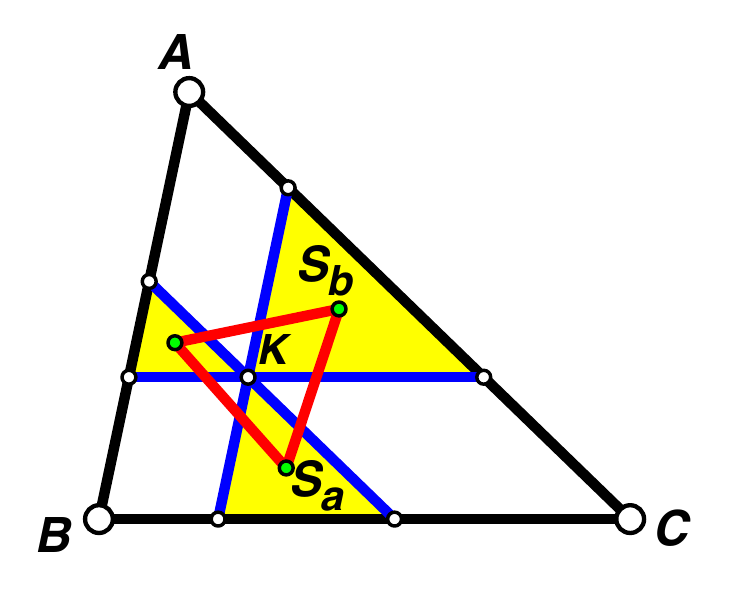}
{}
{$\triangle S_aS_bS_c$ is equilateral.}

\nopagebreak
\pro{\cite{ADGEOM4871}}{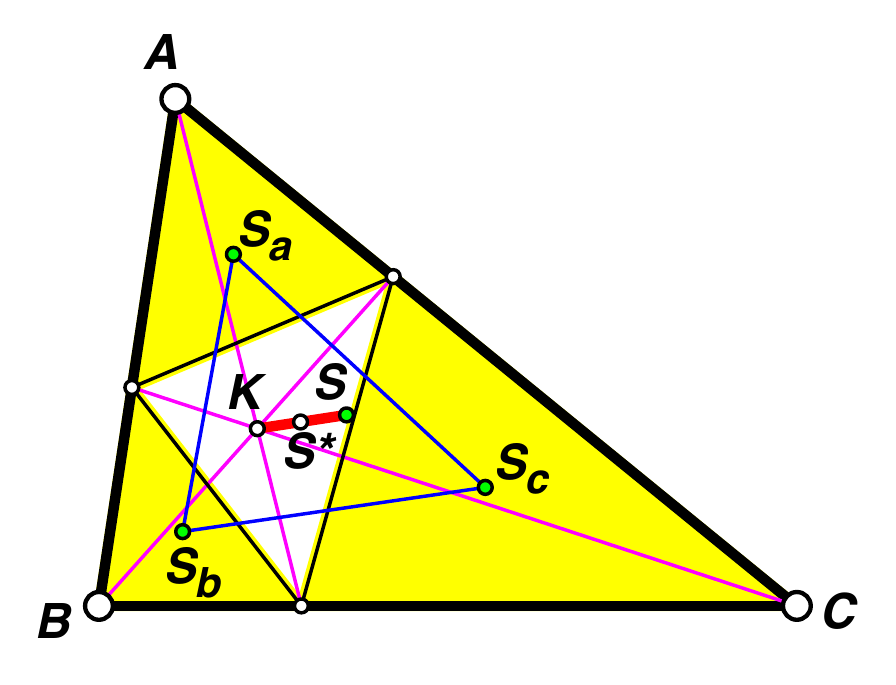}
{$S^*=S(S_a,S_b,S_c)$}
{$K$ -- $S^*$ -- $S$.}
\extra{$S^*=X_{61}$.}

%\newpage
\pro{\cite{PGRAlt}}{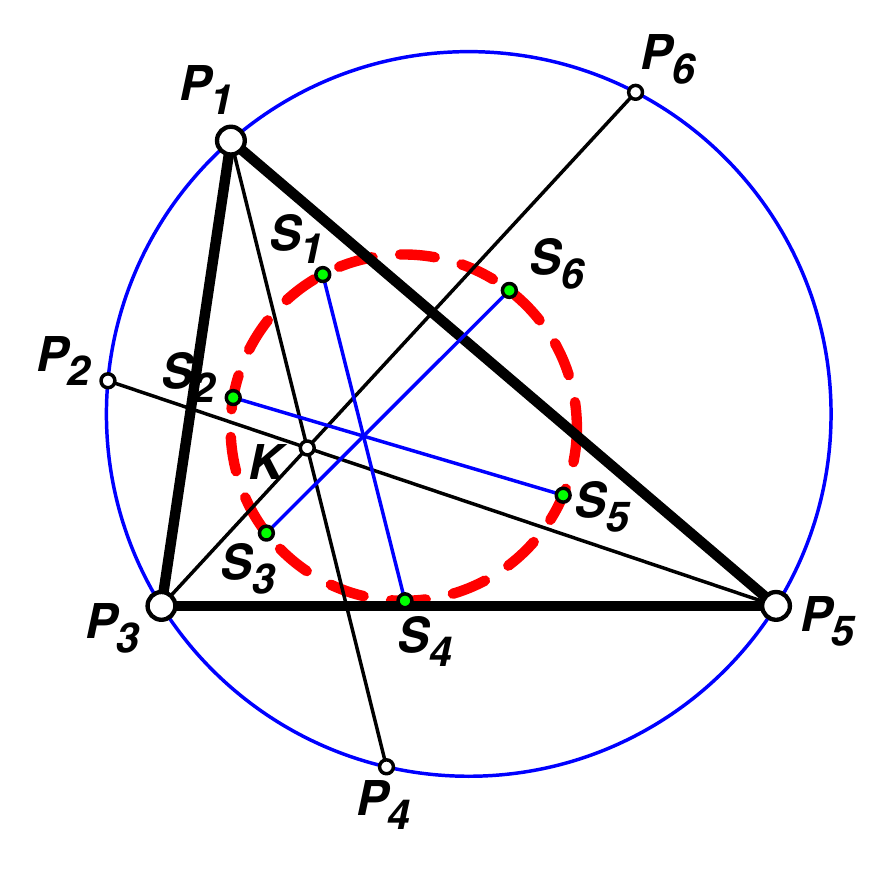}
{$K=K(P_1,P_3,P_5)$ $S_i=S(P_{i-1}P_iP_{i+1})$}
{$S_1$, $S_2$, $S_3$, $S_4$, $S_5$, $S_6$ are concyclic.\\
\then $P_1P_4$, $P_2P_5$, $P_3P_6$ concur.}
\extra{Result is also true using $S'$.}

\mysubsection{formed by $O$}

\nopagebreak
\pro{\cite{Altintas1696}}{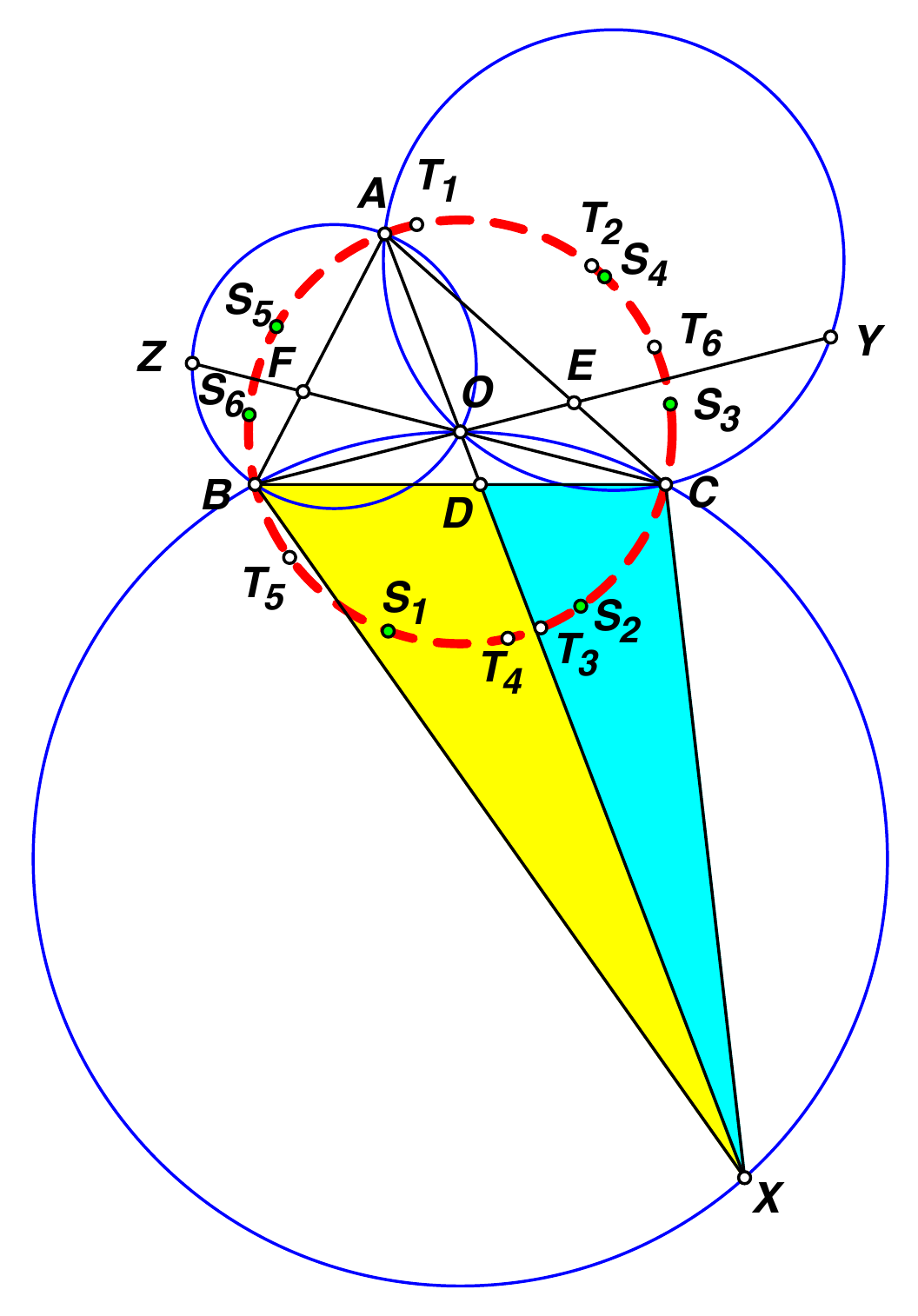}
{}
{$\{S_1,S_2,\ldots S_6, T_1,T_2,\ldots T_6\}$ lie on $\odot ABC$.}

\mysubsection{formed by interior point $P$}

\nopagebreak
\pro{\cite{Altintas1574}}{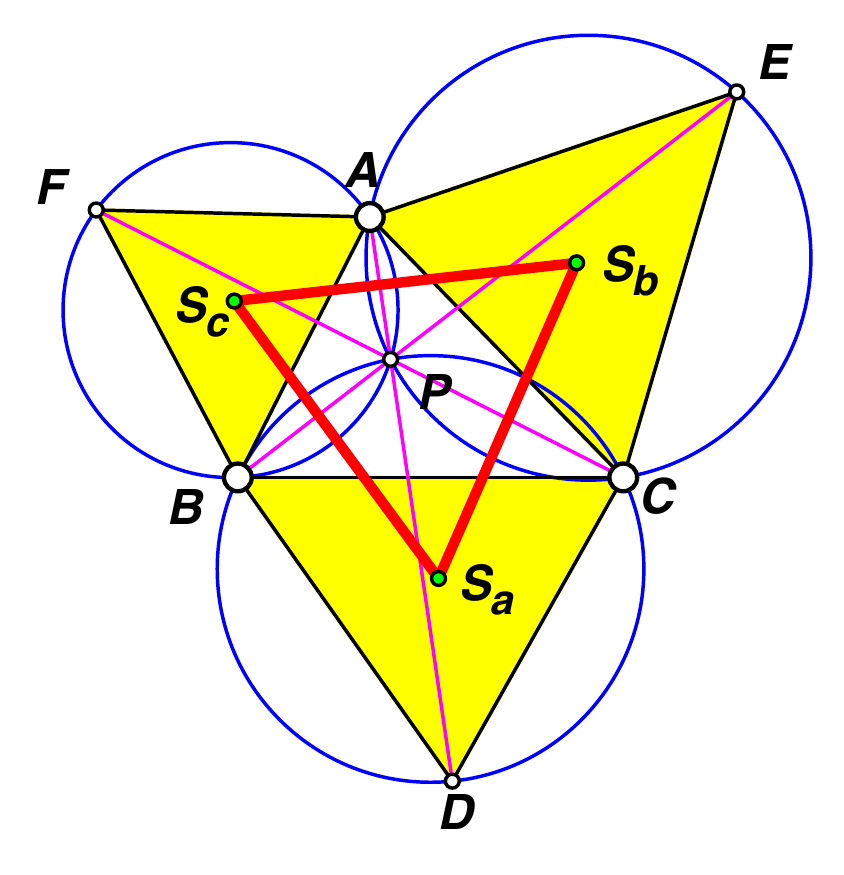}
{}
{$\triangle S_aS_bS_c$ is equilateral.}
\extra{$\triangle T_aT_bT_c$ is equilateral and $\triangle T_aT_bT_c\cong \triangle S_aS_bS_c$.
See \cite{ADGEOM5334}.}

\mysubsection{formed by $X_{13}$}

\nopagebreak
\pro{\cite{ADGEOM4750}}{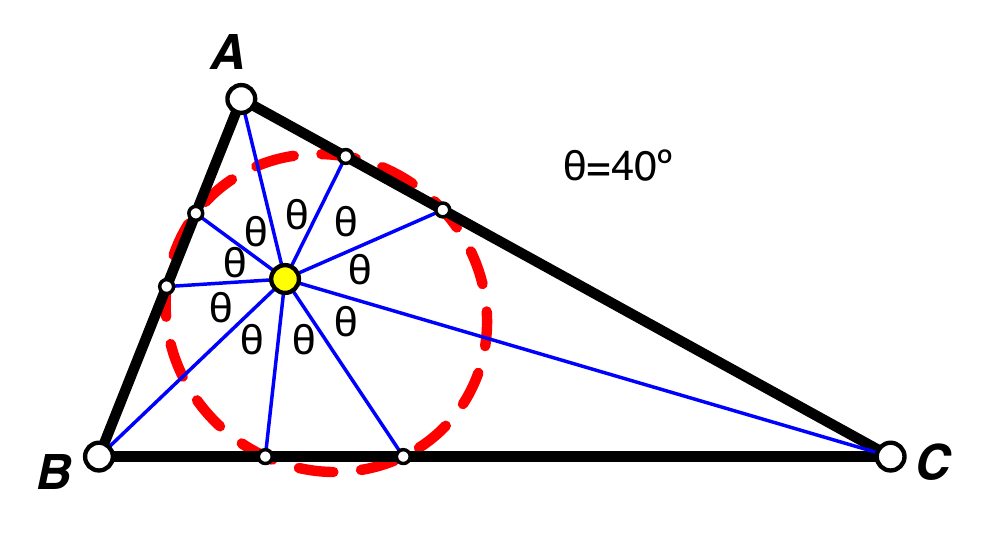}
{Yellow point = $X_{13}$}
{$S$ lies on the axis of an ellipse that passes through the six
points on the perimeter of $\triangle ABC$.}

\pro{\cite[Thm. 7.8.1]{Cosmology}}{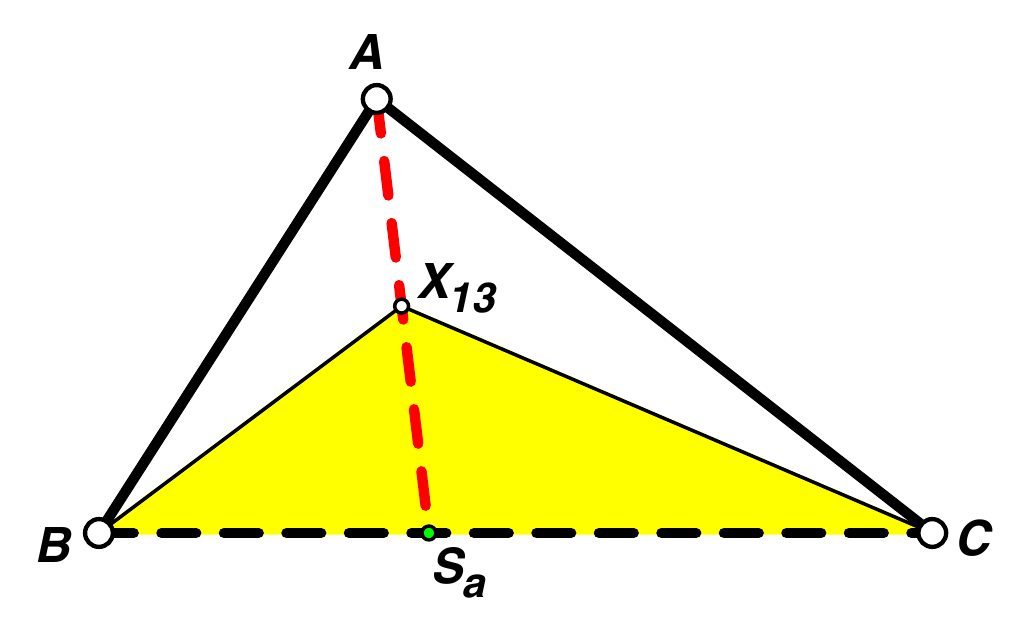}
{}
{$A$ -- $X_{13}$ -- $S_a$\\
\then $B$ -- $S_a$ -- $C$}

\pro{\cite[Thm. 7.8.2]{Cosmology}}{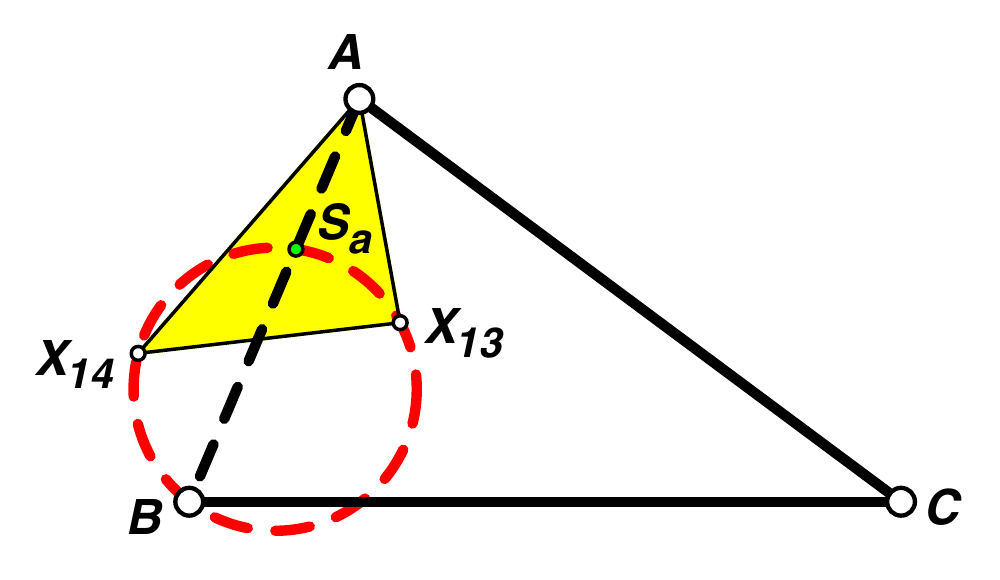}
{$\triangle AX_{14}X_{13}$ has same orientation as $\triangle ABC$}
{$X_{13}$, $S_a$, $X_{14}$, $B$ are concyclic.}

\mysubsection{formed by similar triangles}

\pro{\cite{Rabinowitz/Suppa}}{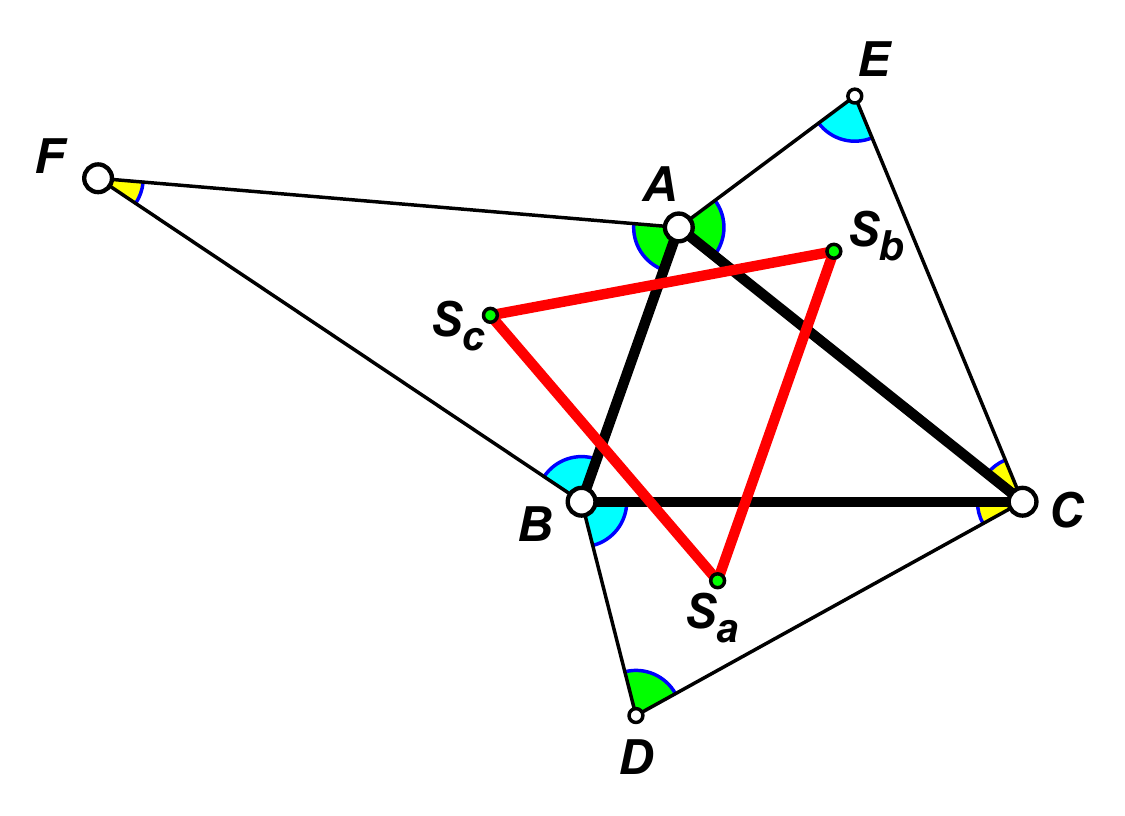}
{}
{$S_aS_bS_c$ form an equilateral triangle.}

\mysubsection{formed by excenters}

\nopagebreak
\pro{\cite{ADGEOM5336}}{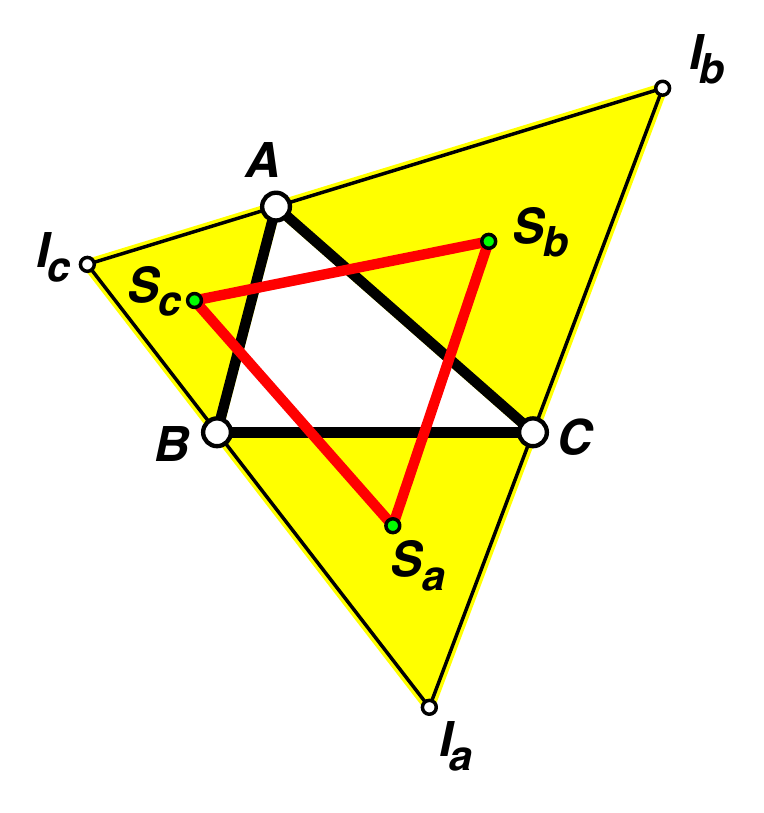}
{}
{$S_aS_bS_c$ form an equilateral triangle.}
\extra{Result is also true using $S'$.}

\mysubsection{formed by six cevians}

\nopagebreak
\pro{\cite{Euclid3140}}{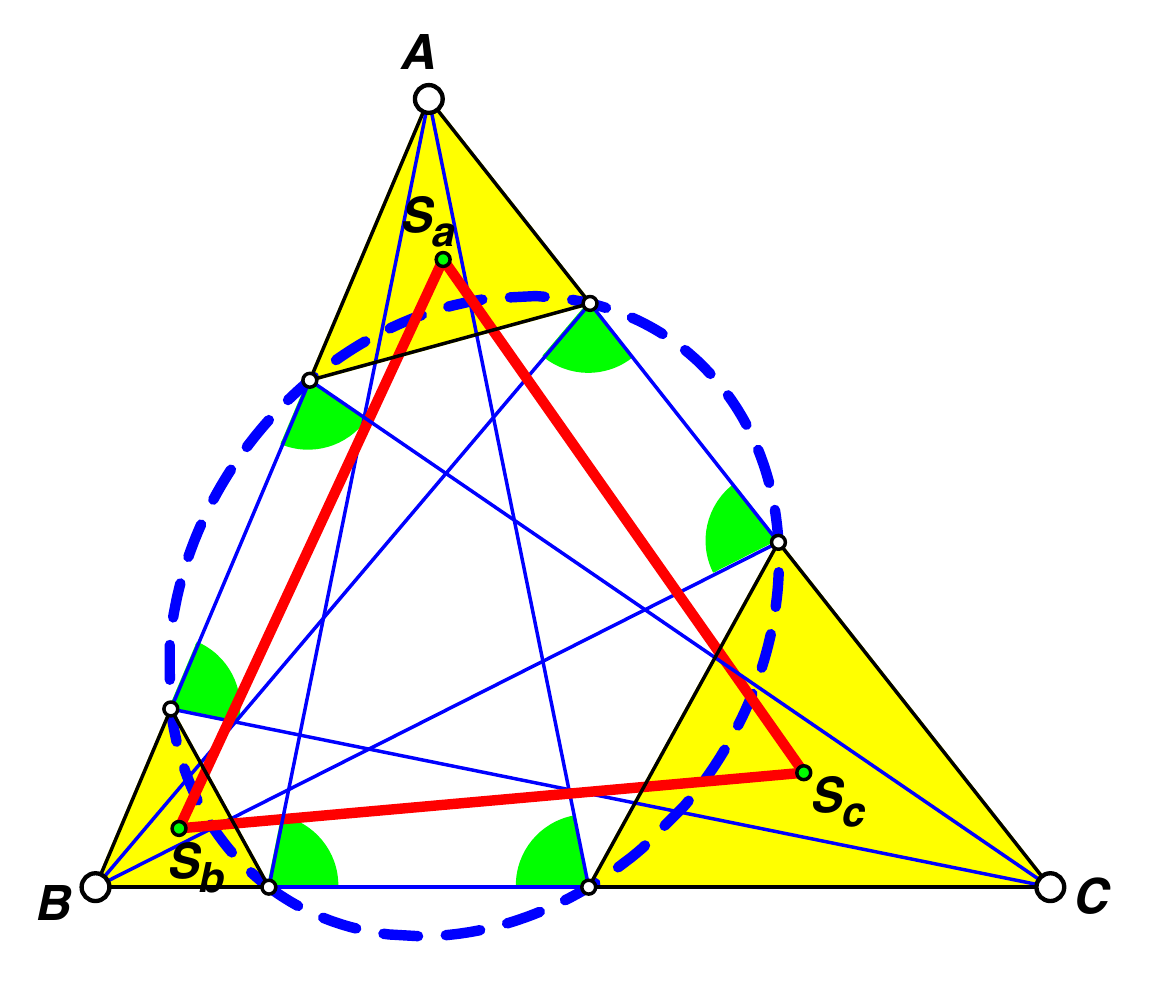}
{}
{$S_aS_bS_c$ is an equilateral triangle.\\
\then Vertices of green angles lie on an ellipse.
}
\extra{Result is also true using $S'$.}

\mysubsection{formed by a Tucker Hexagon}

\nopagebreak
\pro[Tucker Hexagon]{\cite{RG9047}}{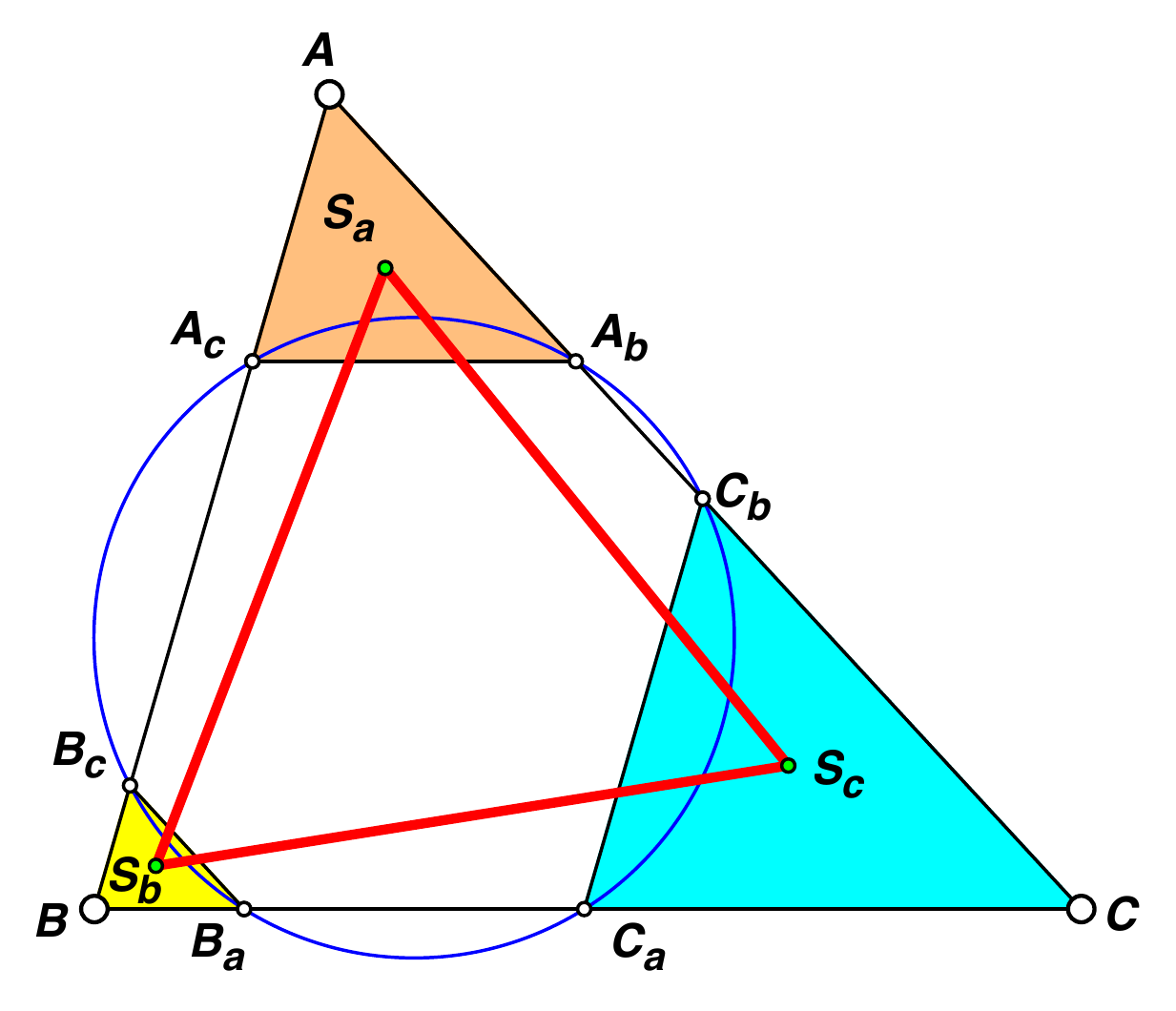}
{$C_aC_b\parallel AB$, $A_bA_c\parallel BC$, $B_cB_a\parallel CA$}
{$\triangle S_aS_bS_c$ is equilateral.}
\extra{The result is true if $S$ is replaced by $S'$.
The two equilateral triangles are congruent.}

\mysubsection{formed by isogonal cevians}

\nopagebreak
\pro{\cite{RG6830}}{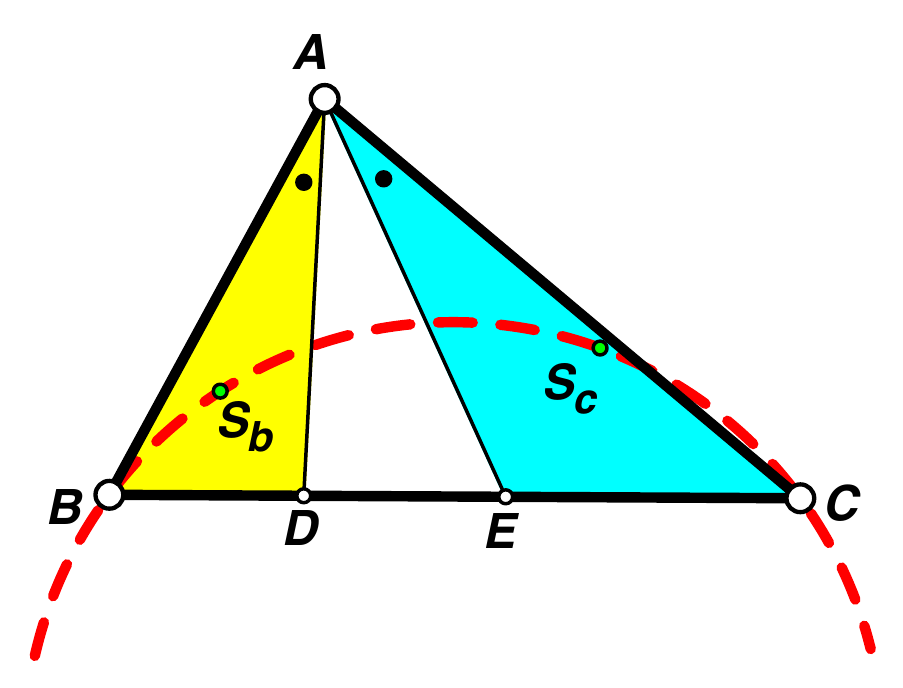}
{}
{$B$, $C$, $S_b$, $S_c$ concyclic}

\pro{\cite{RG6830}}{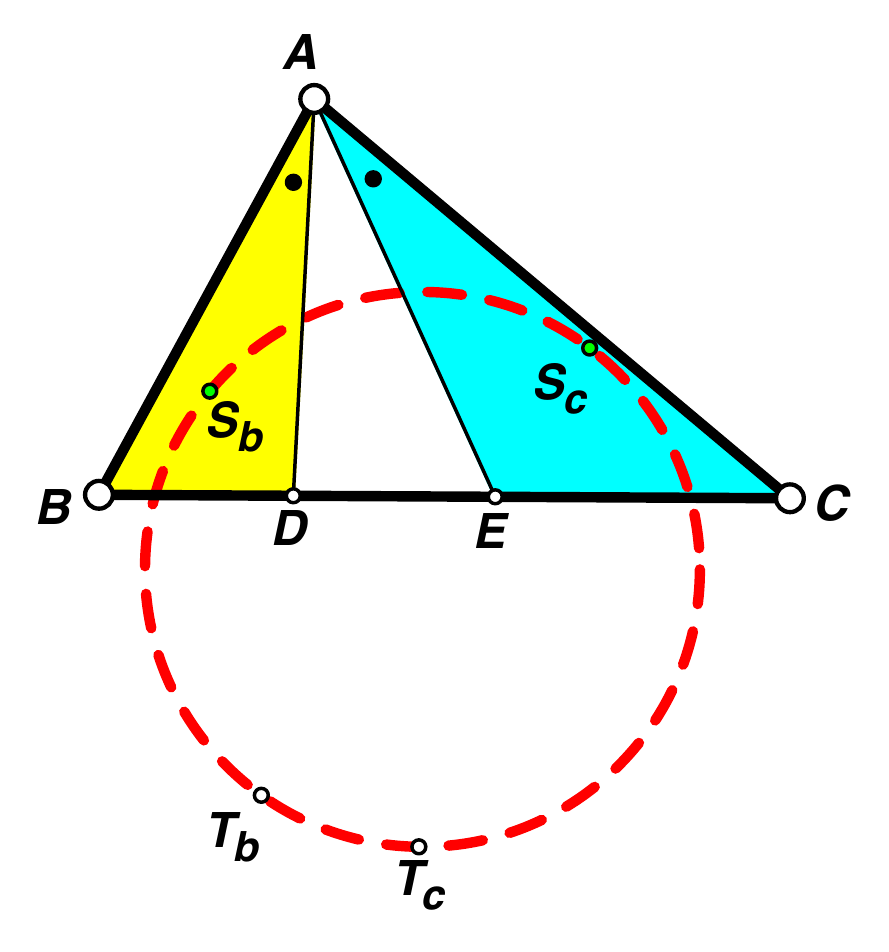}
{}
{$S_b$, $S_c$, $T_b$, $T_c$ concyclic}

\mysubsection{formed by varying vertex $A$}

\nopagebreak
\pro{\cite{RG9079}}{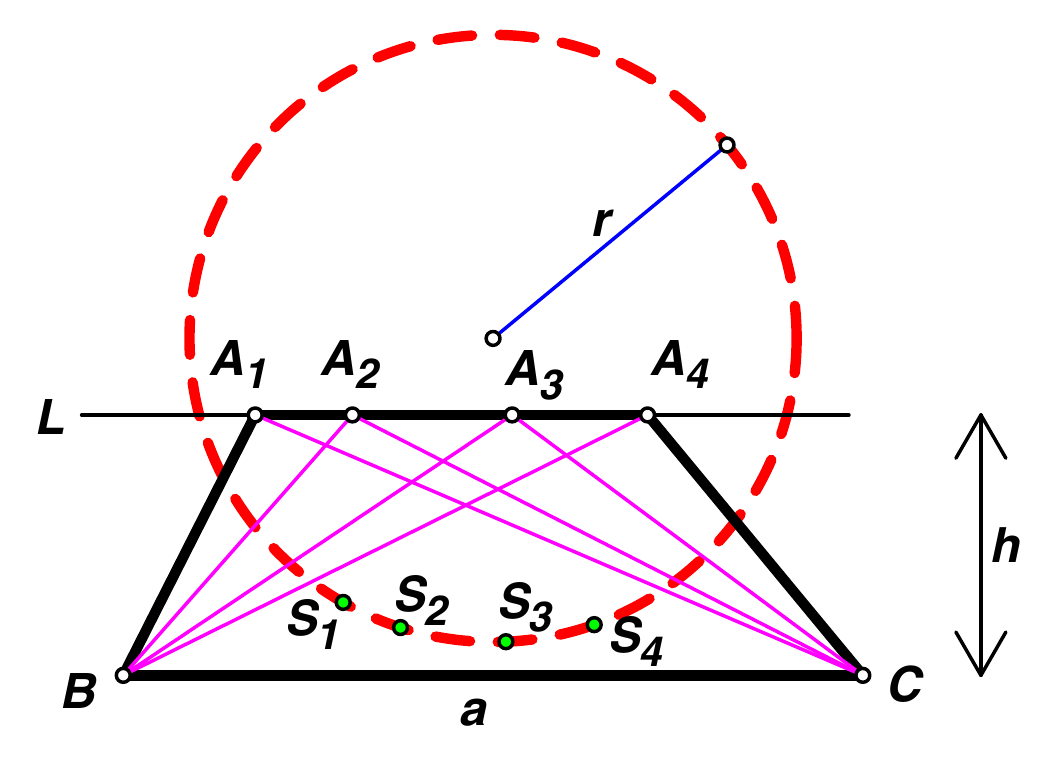}
{$L\parallel BC$, $S_i=S(A_iBC)$}
{$S_1$, $S_2$, $S_3$, $S_4$ are concyclic.\\
\vspace{3pt}
\then $\displaystyle r=\frac{a^2}{2h+a\sqrt3}$}

\pro{\cite{Euclid3253a}}{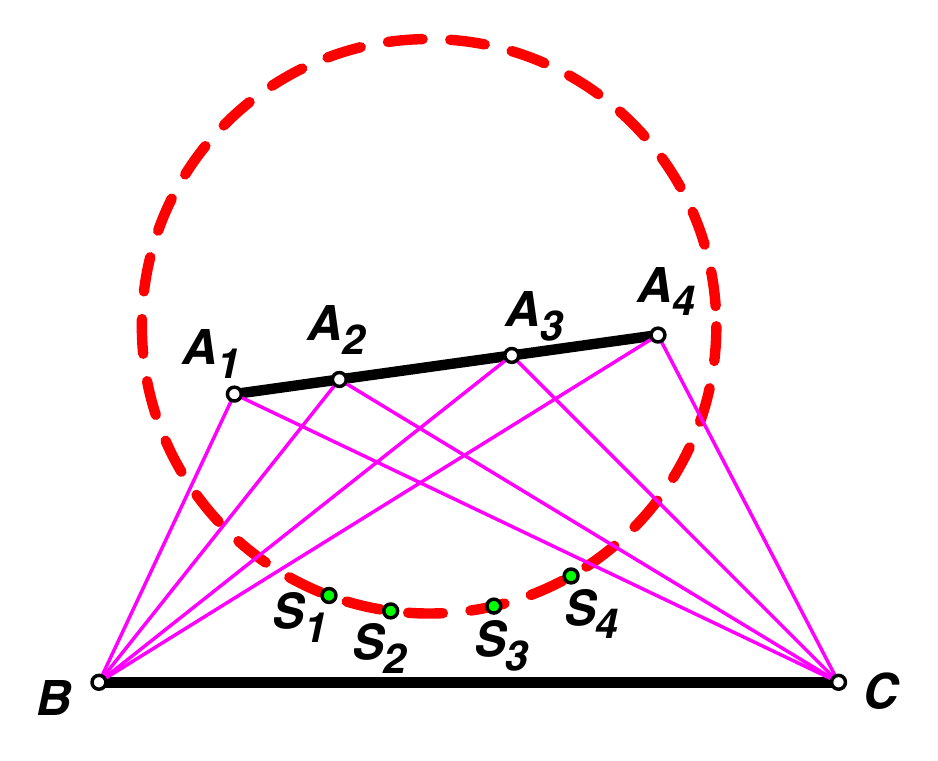}
{$S_i=S(A_iBC)$}
{$S_1$, $S_2$, $S_3$, $S_4$ are concyclic.}

\nopagebreak
\pro{\cite{Euclid3253}}{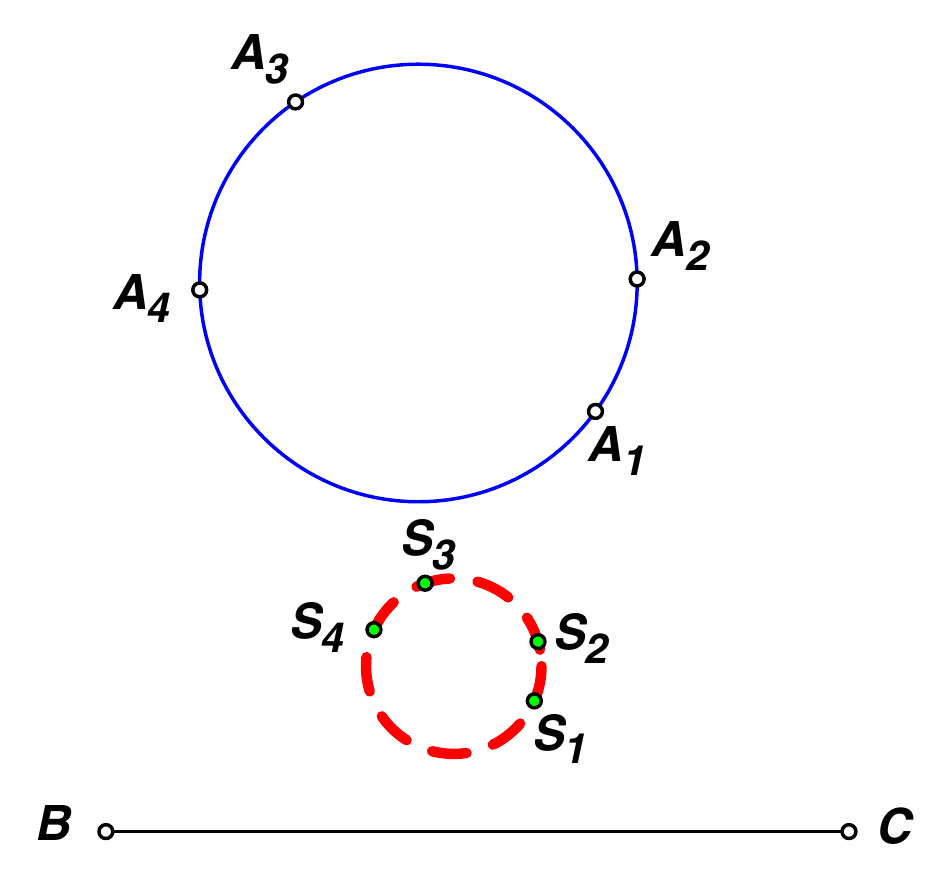}
{$S_i=S(A_iBC)$}
{$S_1$, $S_2$, $S_3$, $S_4$ are concyclic.}

%\pagebreak
%\vspace*{-24pt}
\mysection{Quadrilateral plus $S$}

\mysubsection{cyclic quadrilateral}

\pro{\cite{RG1080}}{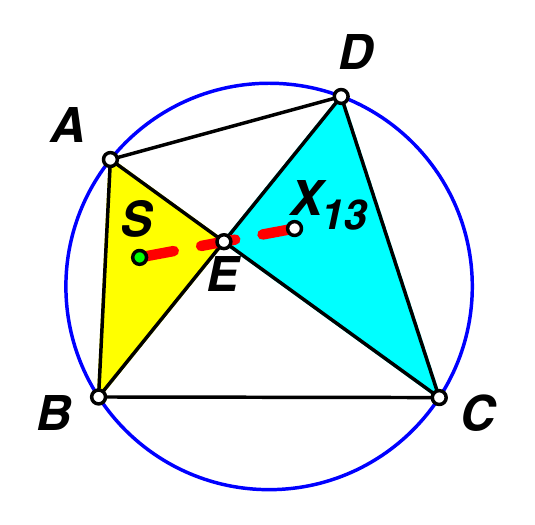}
{}
{$S$ -- $E$ -- $X_{13}$}

\new
\pro{}{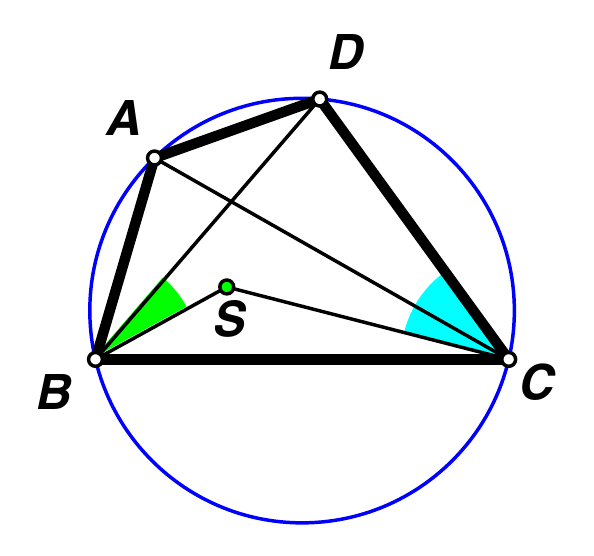}
{}
{$\angle SBD+\angle DCS=60\degrees$}

\new
\pro{}{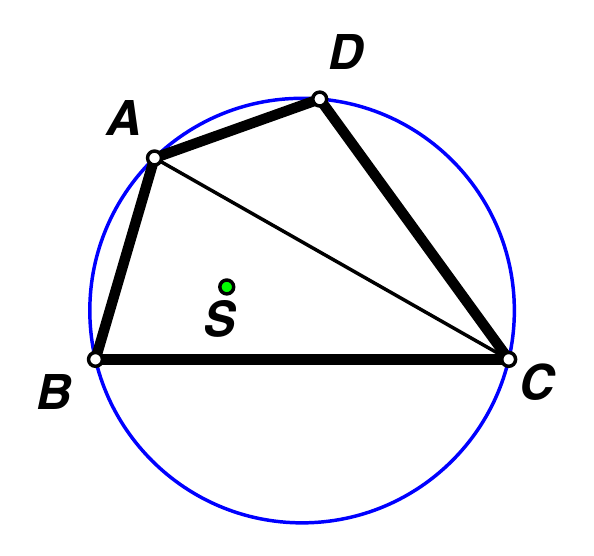}
{}
{$\displaystyle \frac{1}{R(BDS)}+\frac{1}{R(CDS)}=\frac{1}{R(ADS)}$}

\new
\pro{}{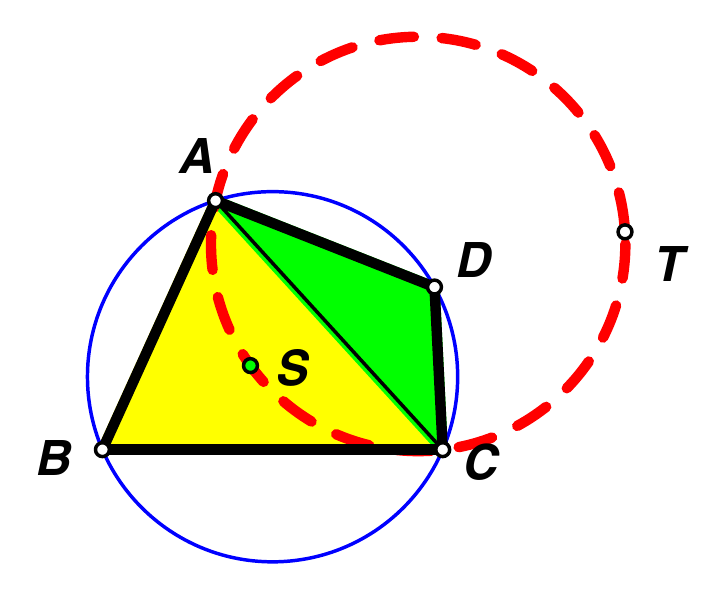}
{$T=X_{16}(ACD)$}
{$A$, $S$, $C$, $T$ concyclic}

\new
\pro{}{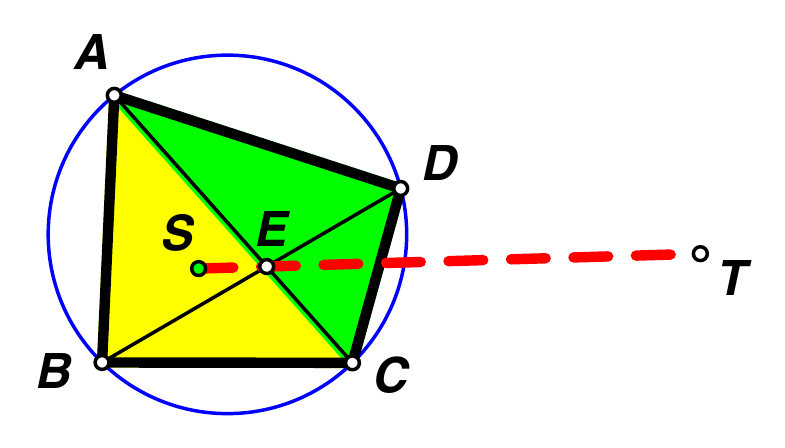}
{$T=X_{16}(ACD)$}
{$S$ -- $E$ -- $T$}

\pro{\cite{Altintas1}}{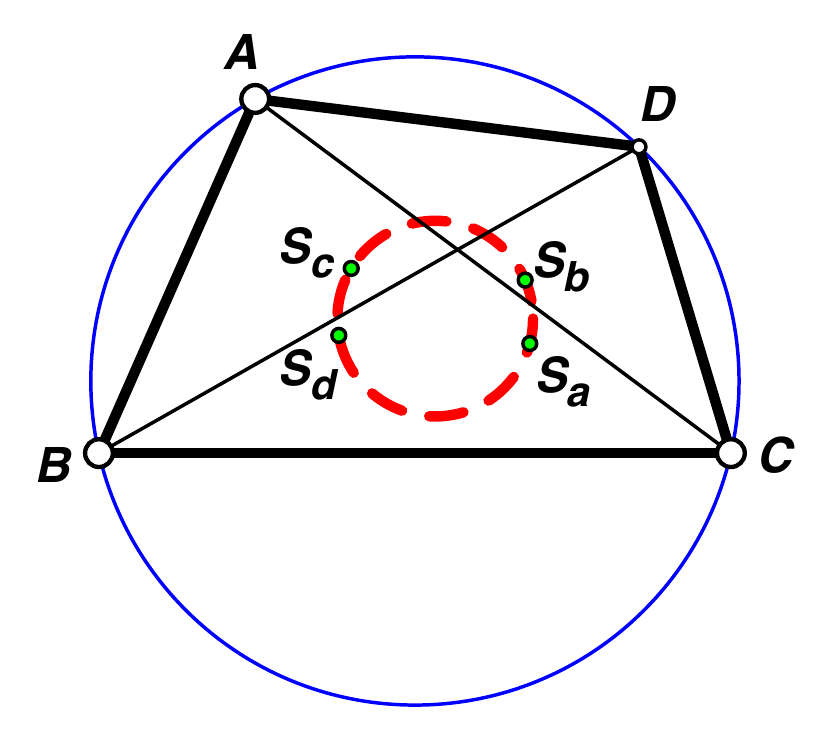}
{$S_a=S(BCD)$}
{$S_a$, $S_b$, $S_c$, $S_d$ are concyclic.}

\mysubsection{square}

\nopagebreak
\new
\pro{\cite{RG6526}}{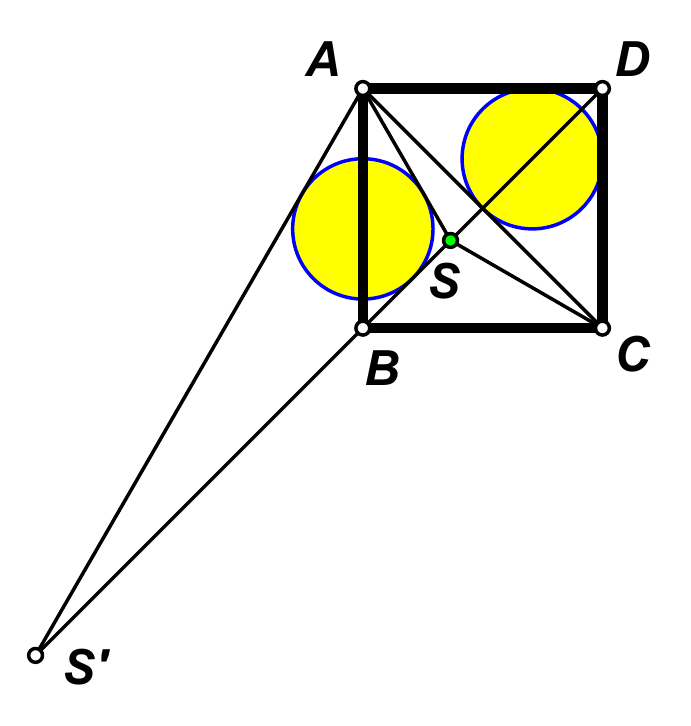}
{}
{Yellow incircles are congruent.}

\new
\pro{\cite{RG6527}}{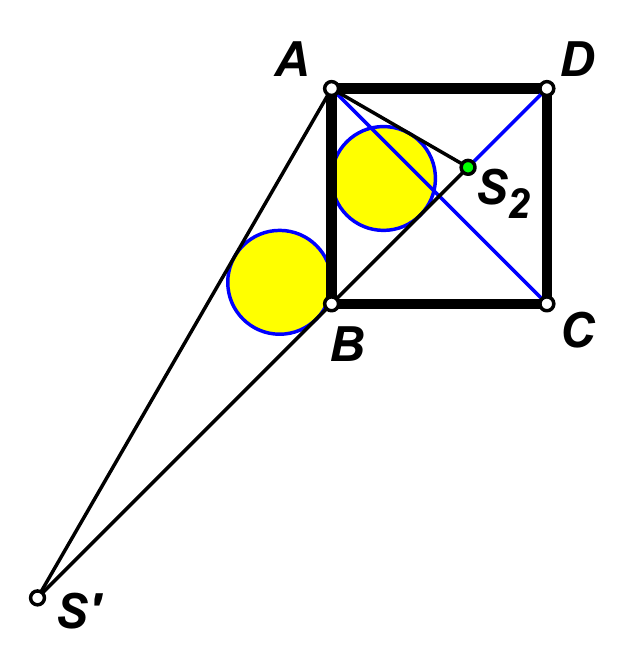}
{$S_2$ is the 1st isodynamic point of $\triangle ACD$.}
{Yellow incircles are congruent.}

\mysubsection{trilateral trapezoid}

\nopagebreak
\new
\pro{}{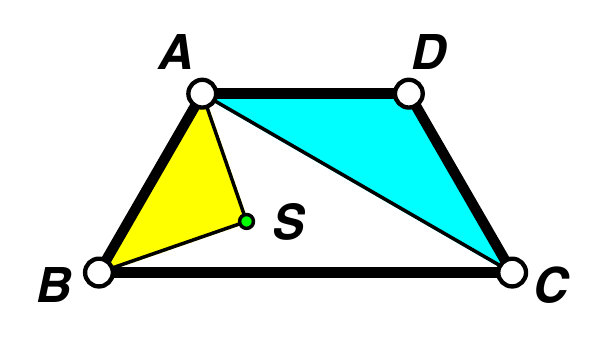}
{}
{$[ADC]=\frac74[BAS]$\\
\then $[CDS]=2[BCS]$\\
\then $[ACS]=3[BDS]$\\
\then $[ABC]=7[BDS]$}

%\newpage
%\vspace*{-24pt}
\mysection{Pentagon plus $S$}

\mysubsection{regular pentagon}

\pro{\cite{RG8232}}{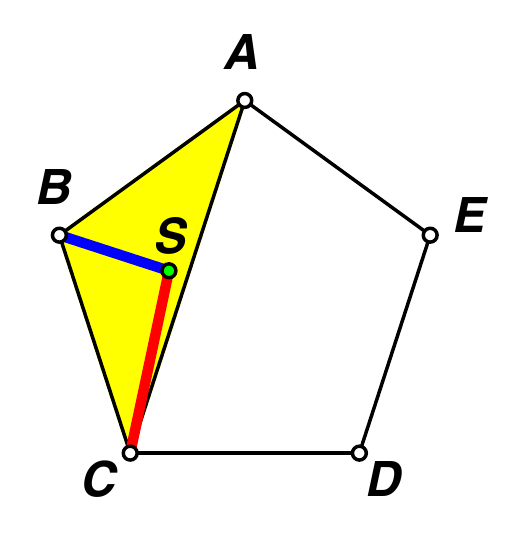}
{}
{$CS/BS=\phi$\\
\then $SD/CS=\sqrt2$}

\pro{\cite{RG8232}}{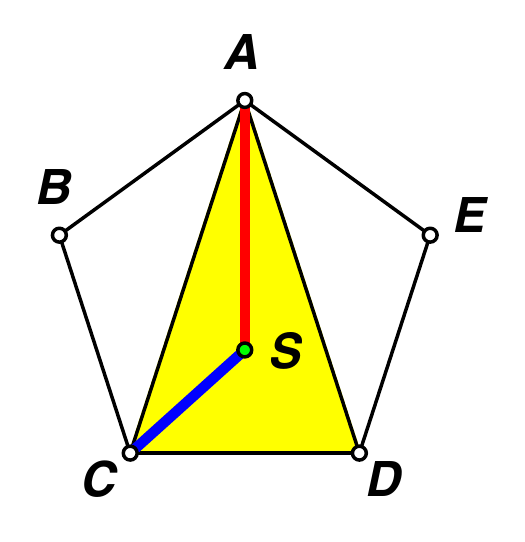}
{}
{$AS/CS=\phi$\\
\then $SE/CS=\sqrt2$}

%\newpage
\vspace*{-15pt}
\mysection{Hexagon plus $S$}

\mysubsection{regular hexagon}

\nopagebreak
\pro{\cite{RG9077}}{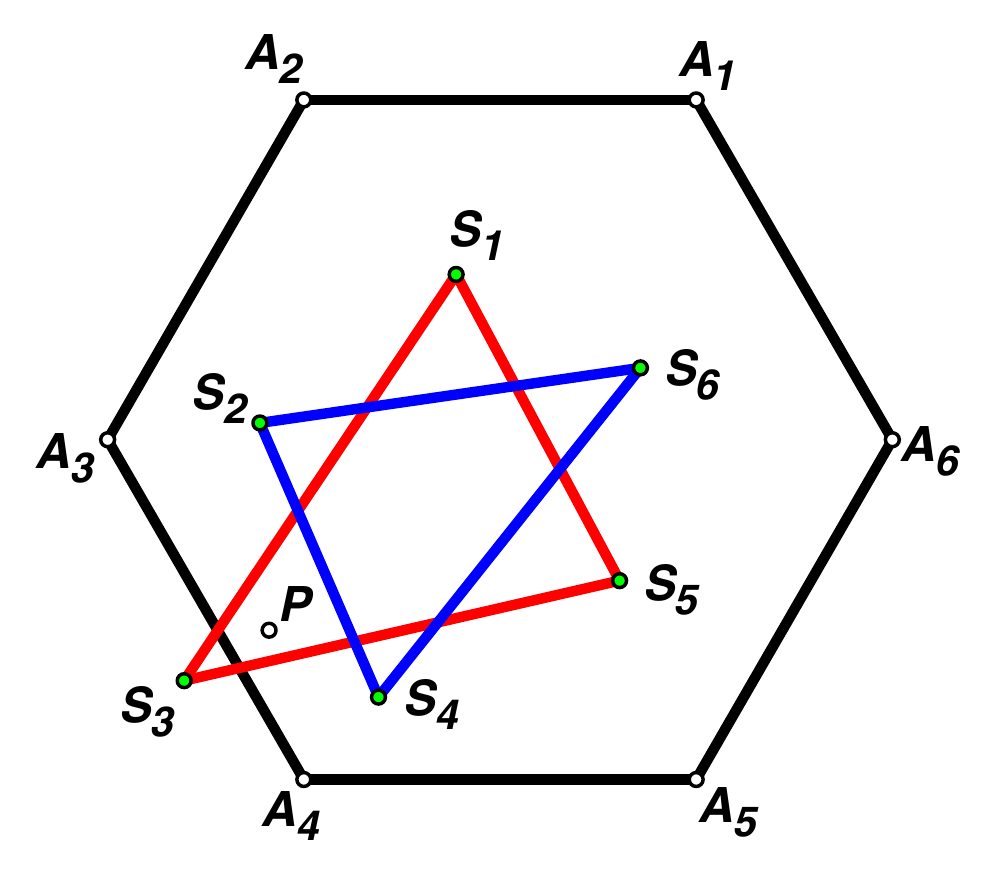}
{$P_i=S(PA_iA_{i+1})$}
{$\triangle S_1S_3S_5\sim\triangle S_4S_6S_2$}

\onecolumn

\section*{Other Properties}

Other properties of the first isodynamic point (discovered by computer) were found by Dekov \cite{Dekov}. A typical result is: the first isodynamic point is the isogonal conjugate
of the inner Fermat point of the anticevian triangle of the outer Fermat point.

Many properties of the first isodynamic point can be found in \cite{ETC15}.
A typical result is: $X_{15}$ is the isogonal conjugate of the isotomic conjugate of $X_{298}$.
Also, lists are given for many of the lines, circles, conics, and cubics that $X_{15}$ lies on.

Many properties of the isodynamic points can be found in the dissertation \cite{Moll}.

Many properties of the isodynamic points can be found in \cite{Cosmology}.

%\newpage
%==========================
% Bibliography
%==========================

% RG results about pentagon, not used: 8108, 8066

\end{document}